\documentclass[11pt]{amsproc}

\usepackage{amssymb,latexsym}
\usepackage{amsmath,amsthm,amsfonts}
\usepackage{psfrag,graphicx}

\newtheorem{thm}{Theorem}[section]
\newtheorem{lemma}[thm]{Lemma}
\newtheorem{cor}[thm]{Corollary}
\newtheorem{prop}[thm]{Proposition}
\newtheorem{conj}[thm]{Conjecture}

\theoremstyle{definition}
\newtheorem{defn}[thm]{Definition}
\newtheorem{example}[thm]{Example}
\newtheorem{xca}[thm]{Exercise}

\theoremstyle{remark}
\newtheorem{obs}[thm]{Observation}
\newtheorem{rmk}[thm]{Remark}
\newtheorem{question}[thm]{Question}

\numberwithin{equation}{section}

\newcommand{\abs}[1]{\lvert#1\rvert}

\newcommand{\bC}{{\mathbb C}}

\newcommand{\bQ}{{\mathbb Q}}
\newcommand{\bR}{{\mathbb R}}
\newcommand{\bZ}{{\mathbb Z}}

\newcommand{\cA}{{\mathcal A}}

\newcommand{\cG}{{\mathcal G}}
\newcommand{\cI}{{\mathcal I}}
\newcommand{\cL}{{\mathcal L}}
\newcommand{\cM}{{\mathcal M}}
\newcommand{\cN}{{\mathcal N}}
\newcommand{\cS}{{\mathcal S}}

\newcommand{\bdy}{\partial}
\newcommand{\del}{\partial}
\newcommand{\wtilde}{\widetilde}

\begin{document}
\bibliographystyle{hplain}

%\date

\title{Four dimensions from two in symplectic topology}

\author{Margaret Symington} 
\address{School of Mathematics, Georgia Inst. of Technology, Atlanta, 
GA 30332}
\email{msyming@math.gatech.edu}
\thanks{The author was supported in part by NSF grant DMS0204368.}

%    General info
\subjclass{Primary 57R15, 55R55; Secondary 53D05, 57M99 }
% 53D05 Sympl. mflds, general
% 55R55 Fiberings with singularities
% 57M99 Low diml topology ``None of above but in this section'' 
% 57R15 Symplectic and contact topology

\date{September 11, 2002}

\keywords{Lagrangian fibrations, almost toric manifolds, moment map}

\begin{abstract}
This paper introduces two-dimensional diagrams
that are slight generalizations of moment map images for toric
four-manifolds and catalogs 
techniques for reading topological and symplectic properties
of a symplectic four-manifold from these diagrams.
The paper offers a purely topological approach to toric manifolds as well as
methods for visualizing other manifolds, including the K3 surface,
and for visualizing certain surgeries.  Most of the techniques extend
to higher dimensions.
\end{abstract}

\maketitle

\tableofcontents

\section{Introduction}

The goal of this paper is to establish techniques for \lq\lq seeing\rq\rq\
certain symplectic four-manifolds: a method for drawing two
dimensional diagrams (base diagrams)
that determine symplectic four-manifolds and 
rules for reading topological and symplectic features
of the manifold from the diagrams.
The inspiration comes from closed toric manifolds, i.e. 
closed symplectic $2n$-manifolds equipped with an effective Hamiltonian
$T^n$ action.  
Delzant's Theorem (Theorem~\ref{Delzant.thm}) says that
the moment map image, a polytope in $\bR^n$, determines the symplectic
manifold and the torus action.
We are concerned here with the symplectic manifold and not the action.
The philosophy therefore is to forget the torus action and then try
to extend Delzant's theorem to a broader class of manifolds using
diagrams that generalize moment map images.
Without the torus action, what remains is a Lagrangian
fibration by tori with singular fibers.
We define the class of almost toric fibrations as
Lagrangian fibrations with certain very simple singularities allowed in
the fibers (Definition~\ref{almosttoric.def}); accordingly, we call
symplectic manifolds equipped with such fibrations almost toric manifolds.
This emphasis on the fibration rather than the torus action gives
a new perspective on toric manifolds in arbitrary dimension 
(including a simple topological
construction from the moment map image) and, generalizing to the almost
toric case, gives ways to visualize four-manifolds such as the K3 surface
and rational balls.

The definition of almost toric fibrations is guided by a classification,
due to Zung~\cite{Zung.II},
of Lagrangian fibrations that have the local structure of an integrable
system with non-degenerate topologically stable singularities.
This classification, up to fiber-preserving symplectomorphism, is in
terms of the base, its geometry, the structure of fibered neighborhoods
of the singular fibers and a Lagrangian Chern class.
The almost toric fibrations are those with no hyperbolic singularities,
and hence in dimension four the topology of the singular fibers is
determined by the geometry of the base. 
Therefore, whenever the Lagrangian Chern class must be trivial (when the
base has the homotopy type of a one-dimensional complex), 
the base determines the manifold.
The two features that distinguish an almost toric fibration from a toric
one are a base that does not immerse, preserving its geometry, in $\bR^2$
and the presence of nodal singularities (the positive transverse
self-intersection points that also arise in Lefschetz fibrations).

While many of the results and techniques of this paper can be generalized
to higher dimensions, unless otherwise specified we restrict attention 
to dimension four so that we can actually draw the base diagrams.
The set of closed toric four-manifolds consists of $\bC P^2$, $S^2\times S^2$
and blow-ups of these manifolds.
The set of closed almost toric manifolds includes the toric ones as well
as sphere bundles over a torus or Mo\"ebius band and blow-ups of them, 
certain torus bundles
over a torus or Klein bottle, the K3 surface and the Enriques surface
(a quotient of the K3 surface by a $\bZ_2$ action). 
See Theorem~\ref{closedalmosttoric.thm}.
However, some of the most interesting almost toric manifolds are
manifolds with boundary.
For instance rational balls that are useful in surgery 
constructions (Section~\ref{rationalballs.sec}) admit almost toric
fibrations that show these surgeries can be done symplectically.
Furthermore, an exotic symplectic four-ball 
(Example~\ref{exoticR4.ex}) -- diffeomorphic to the standard four-ball
but not symplectically embeddable in $(\bR^4,dx\wedge dy)$ --
admits a simple almost toric (even toric) fibration.

The topological and symplectic features of an almost toric manifold that
we explain how to see include 
\begin{itemize}
\item certain surfaces (called visible surfaces), their
symplectic area and self-intersections, and how 
to detect if they are symplectic or Lagrangian (Section~\ref{readingI.sec});
\item the Euler characteristic and first Chern class of the manifold
(Section~\ref{readingII.sec}); and
\item embedded lens spaces, including 
when they admit induced contact structures
and when said contact structures are tight or overtwisted 
(Section~\ref{contact.sec}).
\end{itemize}

We also explain how to blow up and down in the almost toric category
(which is more general than the toric blowing up and down) and how
toric and almost toric fibrations assist in understanding symplectic
surgeries such as generalized rational blowdowns and symplectic
sums.  These are all discussed
in Section~\ref{surgeries.sec}.
As an application of visible surfaces we explain in Section~\ref{K3.sec}
how to \lq\lq see\rq\rq\ surfaces in a K3 surface that reveal the intersection
form 
$E8\oplus E8\oplus 3\left(\begin{smallmatrix} 0 & 1 \\  1 & -2
\end{smallmatrix}\right)$.

In addition to an introduction to almost toric manifolds, 
the material presented here includes some well known facts from toric geometry,
but explained from a non-standard perspective with no reference to torus
actions and moment maps.
For instance, Section~\ref{bdyreduc.sec} 
includes a purely topological construction of toric manifolds 
that is more direct and more general than the usual quotient construction
(cf.~\cite{Audin.torus}).
In Section~\ref{closedtoric.sec} we provide a simple proof of the
well-known classification of closed toric (symplectic) four-manifolds.  Most
symplectic proofs appeal to the classification of complete toric (complex)
surfaces (cf.~\cite{Fulton.toric}); 
ours is purely symplectic, relying only on the notions in this paper.
In Section~\ref{closedalmosttoric.sec} we outline the classification 
of almost toric four-manifolds~\cite{LeungSymington.almosttoric}.

\subsection{Conventions} \label{conventions.sec}
Throughout this paper we assume manifolds are smooth, except possibly
on their boundaries where they may have corners.
We assume that the fibers of a fibration are connected, compact and
without boundary.

We use $\cA$ to denote an integral affine structure (Section~\ref{affine.sec})
and $\cS$ to denote a stratification.  
We write $\phi:(B_1,\cA_1,\cS_1)\rightarrow (B_2,\cA_2,\cS_2)$
to indicate that the map $\phi$ preserves the affine structures and
stratifications, while we  write
$\phi:(B_1,\cA_1)\rightarrow (B_2,\cA_2)$ if only the affine structure
is necessarily preserved by $\phi$.

We use two sets of {\it standard coordinates} and {\it standard 
symplectic forms} on two different spaces:
$(x,y)$ and $\omega_0=dx\wedge dy$ on $\bR^{2n}$ and $(p,q)$ and 
$\omega_0=dp\wedge dq$ on $\bR^n\times T^n$ where we are employing
vector and summation notation so $(p,q)=(p_1,q_1\ldots p_n,q_n)$,
$dp\wedge dq=\sum_i dp_i\wedge dq_i$ and 
$v\frac{\del}{\del p}=\sum_i v_i\frac{\del}{\del p_i}$, $v\in\bR^n$.
We assume the $q$-coordinates are cyclic with period $2\pi$.
By $a[\frac{\del}{\del q}]$ we mean the homology class in
$H_1(\bR^n\times T^n)$ represented by integral curves of the vector
field $\sum_i a_i \frac{\del}{\del q_i}$.

The meanings of light vs. heavy lines and light lines with heavy dashes
are explained in Sections~\ref{bdyreduc.sec} and~\ref{affinenode.sec}.

The following common notations are employed: $A^{-T}$ denotes the inverse
transpose of the matrix $A$; $I$ is the unit interval $[0,1]$; and given
a fibration $\pi:M\rightarrow B$, the vertical bundle of $TM$ means the
subbundle of $TM$ that is the kernel of $\pi_*$.
By a lattice in a vector bundle of rank $n$ we mean a smoothly varying 
lattice (isomorphic to $\bZ^n$) in each fiber of the bundle.

\noindent{\bf Acknowledgments.\,\,}
I am grateful to John Etnyre for encouraging me to write this paper and
to Eugene Lerman for bringing the
work of Nguyen Tien Zung to my attention.
I am also thankful for the helpful comments of John Etnyre, Andy Wand
and the referee on early versions of the paper.

\section{Background}

\subsection{Lagrangian fibrations} \label{basics.sec}

A smooth $2n$ dimensional manifold $M$ is symplectic if it is
equipped with a closed non-degenerate two-form $\omega$ (i.e.
$d\omega=0$ and $\omega^n$ is
nowhere-vanishing).
A submanifold $L$ of a symplectic manifold $(M,\omega)$
is {\it Lagrangian} if ${\rm dim} L=\frac{1}{2}{\rm dim} M$ and
$\omega|_L=0$.
\begin{defn}  
A locally trivial fibration of a symplectic manifold is a
{\it regular Lagrangian fibration} if the fibers are Lagrangian.
More generally, a projection $\pi:(M^{2n},\omega)\rightarrow B^n$, is a 
{\it Lagrangian fibration} if it restricts to a regular Lagrangian
fibration over an open dense set $B_0\subset B$, the set of regular
values.
(We assume throughout this paper that the fibers  of a Lagrangian
fibration are connected, compact and without boundary.)
\end{defn}

A basic feature of symplectic manifolds is expressed in Darboux's theorem
which states that in a neighborhood of any point of a symplectic
$2n$-manifold one can find {\it Darboux coordinates} 
$(x,y)$ such that the symplectic form is standard:  
$\omega_o=dx\wedge dy$.
Given a regular Lagrangian fibration we can instead demand that
the coordinates be adapted to the fibration, as in the following
basic example.

\begin{example} \label{standard.ex}  
Consider standard coordinates $(p,q)$ on $\bR^n\times T^n$ where
the $q$ coordinates have period $2\pi$.
Then projection on the first factor 
$\pi:(\bR^n\times T^n,\omega_0)\rightarrow \bR^n$ 
is a Lagrangian fibration with respect to the standard symplectic
structure $\omega_0= dp\wedge dq$.
\end{example}

Example~\ref{standard.ex} captures the local behavior near any 
fiber of a regular Lagrangian fibration.
This is a consequence of the Arnold-Liouville 
theorem~\cite{Arnold.classicalmech} about integrable systems 
(where the coordinates $(p,q)$ are called {\it action-angle coordinates}).
Rephrased in terms of Lagrangian fibrations we have:

\begin{thm}{  (Arnold-Liouville)} \label{ArnoldLiouville.thm}
Let $F$ be a fiber of a regular Lagrangian fibration
$\pi: (M^{2n},\omega)\rightarrow B^n$.  
Then there is a open neighborhood of $F$ that is fiber-preserving 
symplectomorphic to
$(V\times T^n,\omega_0)\rightarrow V$ where $V\subset \bR^n$.
\end{thm}
Notice in particular that Theorem~\ref{ArnoldLiouville.thm} 
implies that the regular fibers of a Lagrangian fibration are tori.
For a proof of Theorem~\ref{ArnoldLiouville.thm} the reader can 
consult~\cite{Duistermaat.actionangle} or~\cite{Arnold.classicalmech}.
The key steps of the proof are the construction of a natural
torus action in the neighborhood of the regular fiber.

The same argument can be carried out globally 
(see Duistermaat~\cite{Duistermaat.actionangle}) in which case one constructs
a natural fiberwise action of $T^*B$ on $M$, namely
$\alpha\cdot x=\phi_\alpha (x)$ where $\phi_\alpha$ is the time-one
flow of the vector field $X_\alpha$ defined by 
\begin{equation}\label{lattice.eqn}
\omega(\cdot,X_\alpha)=\pi^*\alpha.
\end{equation}
Following this
line of reasoning the isotropy subgroup of this action is a lattice
$\Lambda^*$ in $T^*B$ such that $(T^*B/\Lambda^*,d\lambda)$ 
and $(M,\omega)$ are
locally fiberwise symplectomorphic.  (Here, $\lambda$ is the canonical
$1$-form on $T^*B$, the Liouville form.)
This construction yields another lattice: $\Lambda^{\rm vert}$ belonging 
to the vertical bundle in $TM$ defined by the fibration, namely
$\{X_\alpha|\alpha\in\Lambda^*\}$ where $X_\alpha$ is defined by
Equation~\ref{lattice.eqn}.
The lattice $\Lambda^*\subset T^*B$ also
defines a dual lattice $\Lambda\subset TB$ according to
$u\in \Lambda_b$ if $v^*u\in\bZ$ for every $v^*\in \Lambda_b^*$.  

Here we take a different approach to the globalization in which the 
lattice $\Lambda\subset TB$ appears directly.  Our emphasis on this
lattice rather than its dual follows from a na\"{\i}ve preference for
visualizing tangent vectors to the base rather than covectors.
To find the lattice $\Lambda$ we appeal to the Arnold-Liouville theorem, 
identify $\Lambda$ in local standard coordinate charts 
and then show that this lattice is well defined globally.

The following example covers the local case and serves in part to set notation.
\begin{example} \label{bases.ex} 
Given the fibration 
$\pi:(\bR^n\times T^n,\omega_0)\rightarrow \bR^n$ defined by projection
onto the first factor,
the vectors $\{\frac{\del}{\del p_i}\}$
and $\{\frac{\del}{\del q_i}\}$
provide bases for the {\it standard lattice} $\Lambda_0$ on $T \bR^n$ and the
{\it standard vertical lattice} $\Lambda_0^{\rm vert}$ 
in the vertical bundle of
$T(\bR^n\times T^n)$ induced by the projection $\pi$.
Indeed, the dual lattice $\Lambda^*_0$ on $T^* \bR^n$ is spanned by
$\{dp_i\}$ and each lattice covector $v dp\in T_p\bR^n$ 
defines (by Equation~\ref{lattice.eqn}) the vertical vector field
$v\frac{\del}{\del q}\subset\{p\}\times T (T^n)\subset T(\bR^n\times T^n)$
since $dp\wedge dq(u\frac{\del}{\del p}+
u'\frac{\del}{\del q},v\frac{\del}{\del q})= u\cdot v = 
v dp( u\frac{\del}{\del p}+ u'\frac{\del}{\del q})$ for all $u,u'$.
(Here $u,u'$ are vectors in $\bR^n$ and $v\in\bZ^n$.)
\end{example}

To globalize we rely on the following lemma
about the local structure of fiber-preserving symplectomorphisms.
\begin{lemma} \label{rigidity.lem} 
Let $\varphi$ be a fiber-preserving symplectomorphism of 
$\pi:(\bR^n\times T^n,\omega_0)\rightarrow \bR^n$ where $\pi$ is projection
on the first factor.
Then $\varphi(p,q)=(A^{-T}p+c,Aq+f(p))$ for some $A\in GL(n,\bZ)$,
$c\in \bR^n$ and smooth map $f:\bR^n\rightarrow \bR^n$ such that
$A^{-1}\frac{\del f}{\del p}$ is symmetric.
\end{lemma}

\begin{proof}
Since $\varphi$ is fiber-preserving we can write
$\varphi(p,q)=(\varphi_1(p),\varphi_2(p,q))$.
Requiring $\varphi$ to be a symplectomorphism amounts to imposing
\begin{equation} \label{sympl.eqn}
\varphi^{T}_* J \varphi_* = J \qquad\quad {\rm where} \ 
J= \left(
\begin{smallmatrix}
0 & -I\\
I & 0
\end{smallmatrix}
\right) .
\end{equation}  

In terms of $\varphi_1, \varphi_2$ this is equivalent to 
\begin{equation}
\varphi_1,_p  =  \varphi_2,_{q}^{-T} \quad {\rm and} \qquad
\varphi^T_1,_{p}\varphi_2,_p\ \   {\rm is \ symmetric.}
\label{fconstr.eq}
\end{equation}
Since $\varphi_1,_p$ does not depend on $q$, neither does $\varphi_2,_q$.
Therefore, $\varphi_2= A(p)q + f(p)$ where $A(p)$, for each $p$, is
an element of $GL(2,\bZ)$ and $f(p)$ gives rise to a translation in each
fiber.
Finally, $GL(2,\bZ)$ being a discrete group implies that 
$\varphi_1,p=A^{-T}(p)$ is constant.
Thus we can write $\varphi(p,q)=(A^{-T}p+c,Aq+f(p))$ for $c\in R^2$ where the
condition on $f$ follows from (\ref{fconstr.eq}).
\end{proof}

The above simple lemma captures the rigidity of a Lagrangian fibration and
hints at why, at least in the toric case, one might expect to see a
linear object such as a polytope appear as the base of a Lagrangian
fibration.  Specifically, it leads to:
\begin{thm} \label{affine.thm}
If $\pi:(M,\omega)\rightarrow B$ is a regular
Lagrangian fibration then there are lattices $\Lambda\subset TB$,
$\Lambda^*\subset T^*B$ 
and $\Lambda^{\rm vert}$ in the vertical bundle of $TM$ (induced by $\pi$) 
that, with respect
to standard local coordinates, are the standard lattice, its dual, 
and the standard vertical lattice (of the previous example).
\end{thm}

\begin{proof}
The Arnold-Liouville theorem implies that we can find a cover
$\{U_i\}$ of $B$ and fiber-preserving symplectomorphisms 
$\Phi_i:\pi^{-1}(U_i)\rightarrow V_i\times T^n \subset 
(\bR^n\times T^n,\omega_0)$.
The $\Phi_i$ cover smooth maps $\phi_i:U_i\rightarrow V_i$ which,
by Lemma~\ref{rigidity.lem}, satisfy $(\phi_j\circ \phi_i^{-1})_*\in
GL(n,\bZ)$ on any $U_j\cap U_i$.
Therefore the pullbacks via the $(\phi_i)_*$ of the standard lattices in
$T\bR^n$ and $T^*\bR^n$ give well-defined lattices $\Lambda\subset TB$ and
$\Lambda^*\subset T^*B$.
Since the pairing that relates the lattices 
$\Lambda^*$ and $\Lambda^{\rm vert}$ is defined locally by the symplectic
form $\omega$ it is automatically defined globally.
Therefore the vertical 
lattice $\Lambda^{\rm vert}$ is well-defined globally.
\end{proof}

The rigidity of a regular Lagrangian fibration is
in stark contrast with a symplectic fibration where, for instance,
the fiber can be any symplectic manifold and there is no interdependence
of diffeomorphisms of the base and fiber components of a local
fiber-preserving symplectomorphism. 

\subsection{Integral affine structures} \label{affine.sec}

\begin{defn}   \label{affine.def}
An {\it integral affine structure} $\cA$ on an $n$-manifold $B$ is 
a lattice $\Lambda$ in its tangent bundle, in which case $(B,\cA)$ is
called an {\it integral affine manifold}.
\end{defn}

\begin{example}\label{lattice.ex} 
The standard lattice $\Lambda_0$ generated by the unit vectors tangent
to the coordinate axes in $\bR^n$ 
defines the {\it standard integral affine structure} $\cA_0$ on $\bR^n$.
\end{example}

\begin{example}
The induced lattice $\Lambda$ on the tangent bundle of
the base $B$ of a regular Lagrangian
fibration provides $B$ with an integral affine structure $\cA$.
\end{example}

\begin{prop} \label{charts.prop}
An $n$-manifold $B$ admits an integral affine structure if and only if
it can be
covered by coordinate charts $\{U_i,h_i\}$, $h_i:U_i\rightarrow \bR^n$
such that the map $h_j\circ h_i^{-1}$, wherever defined, is an element of
$AGL(n,\bZ)$, i.e.
a map of the form $\Phi(x)=Ax+b$ where $A\in GL(n,\bZ)$ and $b\in \bR^n$.
\end{prop}
\begin{proof}
Given an affine structure $\cA$ on $B$, affine coordinate charts 
$h_i:(U_i,\cA)\rightarrow (\bR^n,\cA_0)$ exist since all affine structures
(in a given dimension) are isomorphic.
Then $(h_j\circ h_i^{-1})_*\in GL(n,\bZ)$ so $h_j\circ h_i^{-1}\in AGL(n,\bZ)$.
Given such a cover, the integral affine structure $\cA$ is defined
as the pull-back (via the $h_i$) of the standard affine structure 
$\cA_0$ on $R^n$.
\end{proof}

\begin{rmk} \label{flat.rmk}   
Note that the tangent bundle of an integral affine manifold has
a flat connection, namely the pull-back via the coordinate charts
of the flat connection on $T\bR^n$ (since it is preserved under the
transition maps).
\end{rmk}

\begin{defn}   \label{rational.def}
A $k$-dimensional submanifold $P\subset (B,\cA)$ is 
{\it affine planar} if at each point $b\in P$ the subspace 
$T_bP$ is spanned by $k$ vectors in $\Lambda$, or
equivalently, the subspace in $T^*_b B$ consisting of covectors
that vanish on $T_bP$ is spanned by $n-k$ covectors in $\Lambda^*$.
Accordingly, a one dimensional submanifold is {\it affine linear} 
if at every point it has a tangent vector in $\Lambda$.
\end{defn}

In general there is no metric compatible with an integral affine
structure.  However, we can still define the following $AGL(n,\bZ)$-invariant
quantity:
\begin{defn}  
Consider a parameterized curve $\gamma:[0,1]\rightarrow (B,\cA)$ and 
define $s(t),v(t)$ so that $s(t)>0$ is a scaler, 
$v(t)\in\Lambda_{\gamma(t)}$ is a primitive tangent vector and 
$\gamma'(t)=s(t)v(t)$.
Then the {\it affine length} of $\gamma$ is 
\begin{equation}
\int_{0}^1 s(t) dt.
\end{equation}
\end{defn}

\subsection{Monodromy -- topological and affine} 
\label{monodromy.sec}

An important measure of the non-triviality of a torus bundle is its
monodromy.  
Interestingly, any topological monodromy of a regular Lagrangian fibration 
is reflected in the global geometry of the base.
In this section we review the definition of (topological) monodromy,
define the affine monodromy of an integral affine structure and
then give the relation between the topological 
monodromy of a regular Lagrangian fibration and the affine monodromy
of its base.

Let $\cM(T^n)$ denote the {\it mapping class group} of the $n$-torus, 
namely the set of isotopy classes of self-diffeomorphisms of the torus.
Let $\zeta_M$ be an $n$-torus bundle given by $\pi:M\rightarrow B$.
Then a {\it monodromy representation} of $\zeta_M$ is
a homomorphism $\Psi:\pi_1(B,b)\rightarrow \cM(T^n)$ such that for 
any $[\gamma]\in \pi_1(B,b)$ represented by a based 
loop $\gamma:S^1\rightarrow B$,
the total space of $\gamma^* \zeta_M$ is diffeomorphic to 
$I\times T^n/((0,q)\sim (1,\Psi[\gamma](q))$.
We can identify each isotopy class $\Psi[\gamma]$ with 
the induced automorphism on $H_1(F_b,\bZ)$ where $F_b$ is the fiber over the 
base point $b$.
Choosing a basis for $H_1(F_b,\bZ)$ we then get an induced map from
$\pi_1(B,b)$ into $GL(n,\bZ)$.
In general, by the {\it monodromy} of a regular torus fibration we mean this
induced map for some choice of monodromy representation and some
choice of basis for $H_1(F_b,\bZ)$.

\begin{rmk} \label{orient.rmk}  
Note that the base of a regular Lagrangian fibration need not be orientable
and hence the fiber need not be oriented.
\end{rmk}

Similar to a torus fibration, 
an integral affine structure $\cA$ on a manifold $B$ has an
{\it affine monodromy representation} 
$\Psi_B:\pi_1(B,b)\rightarrow \rm{Aut(\Lambda_b)}$
where $\Lambda_b$ is the restriction to $T_b B$ of the lattice $\Lambda$ that
defines $\cA$.
Specifically,
if we identify $(T_bB,\Lambda_b)$ with $(\bR^n,\bZ^n)$ and let
$\gamma$ be a based loop,
$\Psi_B[\gamma]$ is the automorphism of $(\bR^n,\bZ^n)$ such that 
$\gamma^* (TB,\Lambda)$ is isomorphic 
to $I\times (\bR^n,\bZ^n)/((0,p)\sim (1,\Psi_B[\gamma](p))$.

Now recall that a regular Lagrangian fibration 
$\pi:(M,\omega)\rightarrow B$ induces lattices $\Lambda\subset TB$,
$\Lambda^*\subset T^*B$ and $\Lambda^{\rm vert}\subset TM$ (see
Section~\ref{basics.sec}).  
This allows us to relate the topological monodromy of the underlying
torus fibration to the affine monodromy of the induced affine structure
on the base.

\begin{prop} \label{baselink.prop}
Suppose that, with respect to local standard coordinates near a regular fiber,
the topological monodromy around a loop in the base of a regular
Lagrangian fibration is given by $A\in GL(n,\bZ)$.  Then the affine monodromy
is given by its inverse transpose $A^{-T}$.
\end{prop}

\begin{xca}   \label{baselink.ex} Prove this.  Note that the 
monodromy of the lattice $\Lambda^*$ would be $A$.
\end{xca}

\begin{prop} \label{link.prop}
Suppose $\Phi[\gamma],\Phi_B[\gamma]$ are the topological and affine
monodromies along $\gamma$ with respect to standard local coordinates
near the regular fiber over the base point $b$.
Then $\Phi[\gamma]$ has an invariant sublattice spanned by 
$v\frac{\del}{\del q}$ for some $v\in \bZ^n$
if and only if $\Phi_B[\gamma]$ has an invariant sublattice
spanned by $u\frac{\del}{\del p}$, $u\in\bZ^n$ such that $v\cdot u=0$.
\end{prop}

\begin{xca}   \label{link.ex} Prove this.
\end{xca}

\section{Toric manifolds} \label{toric.sec}

\subsection{More background} \label{morebackground.sec}

If a group $G$ acts on a symplectic manifold $(M,\omega)$ then each element
$\xi$ of the Lie algebra $\cG$ induces a vector field $X_\xi$,
the infinitesimal action.
A smooth function $f:M\rightarrow \bR$ also defines a (Hamiltonian) vector
field $X_f$ by $\omega(\cdot,X_f)=df$.
\begin{defn}  
A group action on a symplectic manifold $(M,\omega)$ is 
{\it Hamiltonian} if there is
a Lie algebra homomorphism from $\cG$
to the smooth functions on $M$ (with the Poisson bracket)
taking $\xi$ to $f_\xi$ so that $X_\xi=X_{f_\xi}$.
\end{defn}
If a group action is Hamiltonian then there is an associated
{\it moment map} $\mu:M\rightarrow \cG^*$
defined implicitly by $<\xi,\mu(x)>=f_\xi(x)$ for $\xi\in\cG$ and $x\in M$.
Thus for $k$-torus actions the image of the moment map is in $\bR^k$.

A fundamental theorem about Hamiltonian torus actions is:
\begin{thm}[Atiyah~\cite{Atiyah.convexity}, 
Guillemin and Sternberg~\cite{GuilleminSternberg.convexity}]
\label{convexity.thm} 
The moment map image $\Delta:=\mu(M)$ for a Hamiltonian
$k$-torus action on a closed symplectic manifold $(M,\omega)$ 
is a convex polytope.
\end{thm}
When the dimension of the torus acting on $(M^{2n},\omega)$ is maximal,
i.e. $k=n$, then $M$ is called {\it toric}.  In the toric case
the regular fibers of $\mu$ are tori of dimension $n$ 
and singular fibers are tori of lower dimension.
In this case
the moment map has a simple interpretation as the orbit map.
However, the geometry of $\Delta$ reveals much more than just the
orbit structure:
\begin{thm}[Delzant~\cite{Delzant.moment}] \label{Delzant.thm}
If a closed symplectic manifold $(M^{2n},\omega)$ is equipped with
an effective Hamiltonian $n$-torus action, then the image of the moment
map $\Delta$ determines the manifold $M$, its symplectic structure $\omega$ 
and the torus action.
\end{thm}

\begin{rmk}   \label{proper.rmk}
Both the Atiyah, Guillemin and Sternberg convexity theorem and Delzant's
theorem remain true if the hypothesis that the manifold be closed is
replaced by the moment map being proper and the fibers of the moment
map being connected (cf.~\cite{LermanTolman.orbifolds}).
\end{rmk}

Not all polytopes can be the moment map image for a closed manifold; the
following is standard:
\begin{prop} \label{toricbasics.prop}
Suppose $\Delta$ is the moment map image of a closed symplectic manifold.
Then each $(n-1)$-dimensional face of $\Delta$ has a primitive
integral normal vector and each 
$(n-k)$-dimensional face, $k\ge 2$, lies at the intersection of
$k$ faces of dimension $n-1$.
Furthermore, smoothness of the manifold implies that the $k$ 
primitive integral vectors normal to the $k$ faces intersecting
at an $(n-k)$-dimensional face span a $k$-dimensional 
sublattice of the standard integral lattice $\Lambda_0$.
(If not, $M$ has orbifold singularities.)
\end{prop}
For further information on Hamiltonian torus actions 
consult~\cite{Audin.torus}.

Now forget the torus action but remember the
projection to $\bR^n$ defined by the moment map $\mu$.
All fibers of the moment map are tori, but their dimensions vary:
the preimage of a point on the interior of $\Delta$ is an $n$-torus
while the preimage of a point in the interior of a $k$-dimensional face
of $\Delta$ is a $k$-torus. 
A standard fact from toric geometry is that the symplectic form 
vanishes on fibers of the moment map.
Therefore regular fibers are Lagrangian submanifolds and
we can think of the moment map $\mu$ as providing a Lagrangian fibration.

If $(M,\omega)$ is toric and either is non-compact or has non-empty  
boundary then preimages under
the moment map may be disconnected, a union of tori
of possibly varying dimensions (as in Example~\ref{exoticR4.ex}).
In this case $\Delta$ is no longer the orbit space; rather
we have $\mu=\Phi\circ \pi$ where $\pi:(M,\omega)\rightarrow B$ is
the induced Lagrangian fibration and $\Phi:B\rightarrow \bR^n$ is
an immersion.
The immersion $\Phi$ respects the geometry of the base in the
following sense:  the induced affine structure on the regular values of $\pi$
extends in a natural way an affine structure $\cA$ on $B$
and we have $\Phi:(B,\cA)\rightarrow (\bR^2,\cA_0)$.
(This will become clear in the next section.) 
Furthermore, the base $B$ has a natural stratification $\cS$:
the $l$-stratum of $(B,\cS)$ is the set of points in 
$b\in B$ such that $\pi^{-1}(b)$ is a torus of dimension $l$.
(That this defines a stratification follows from the normal forms
for the fibered neighborhoods of $l$-dimensional fibers.)
Notice that when $M$ is a closed manifold the $l$-stratum of
$\Delta$ is the union of the 
interiors of the $l$-dimensional faces, $l>0$, and
the $0$-stratum is the union of the vertices.

\begin{defn}   \label{toricfibr.def}
A Lagrangian fibration $\pi:(M^{2n},\omega)\rightarrow (B,\cA,\cS)$ is
a {\it toric fibration} if there is a Hamiltonian $n$-torus action and
an immersion
$\Phi:(B,\cA)\rightarrow (\bR^n,\cA_0)$ such that $\Phi\circ \pi$ is the
corresponding moment map and $\cS$ is the induced stratification.
In this case we call $(B,\cA,\cS)$ a {\it toric base}.
\end{defn}

Our goal to generalize part of Delzant's theorem can then be summarized
by:
\begin{question} \label{Lagfibr2.quest}  
What assumptions about a triple $(B,\cA,\cS)$ ensure that
it determines a unique symplectic manifold?
\end{question}
An answer is provided, for two-dimensional bases,  
by Theorem~\ref{almosttoricbase.thm} and  
Corollary~\ref{almosttoricbase.cor}.  
As a warm-up, in the next section we give criteria for recognizing
whether or not $(B,\cA,\cS)$ is the base of a toric fibration
(Theorem~\ref{topologicaltoric.thm}).
 
\subsection{Construction via boundary reduction} \label{bdyreduc.sec}

In this section we explain how to construct a toric fibration
from a toric base $(B,\cA,\cS)$.
The idea is to start with a regular Lagrangian fibration over the
base $(B,\cA)$ and then collapse certain fibers to lower dimensional tori 
so as to obtain the stratification $\cS$.

The collapsing of fibers is achieved by boundary reduction, an operation
originally defined in~\cite{Symington.blowdowns}:
\begin{prop} \label{bdyreduc.prop} 
Let $(M,\omega)$ be a symplectic manifold with boundary such
that a smooth component $Y$ of $\bdy M$ is a circle bundle
over a manifold $\Sigma$.
Suppose also that the tangent vectors to the circle fibers lie
in  the kernel of $\omega|_{Y}$.
Then there is a projection $\rho:(M,\omega)\rightarrow (M',\omega')$
and an embedding $\phi:\Sigma\rightarrow M'$
such that $\rho(Y)=\phi(\Sigma)$, $\rho|_{M-Y}$ is a symplectomorphism
onto $M'-\phi(\Sigma)$ and $\phi(\Sigma)$ is a symplectic submanifold.
\end{prop}

\begin{proof}
The circle bundle structure of $Y$ gives the existence of the manifold
$M'$ and the embedding $\phi$.
Using a local model for the symplectic neighborhood of $Y\subset M$
it is not hard to see the induced symplectic structure on a neighborhood
of $\rho(Y)\subset M'$.  Essentially, it arises from a fiberwise
application of Example~\ref{cylinder.ex}.
That  $\rho(Y)$ is a symplectic submanifold can be seen in the local
models or else by noting
that, for points $x\in Y$,  $\rho_*(T_xY)$ is the quotient of
$T_xY$ by the kernel of $\omega|_{T_xY}$, thereby showing that $\omega'$ is 
non-degenerate on $\rho(Y)$.  
\end{proof}

\begin{defn} \label{bdyreduc.def}   
With the notation of Proposition~\ref{bdyreduc.prop},
$(M',\omega')=\rho(M,\omega)$ is the
{\it (symplectic) boundary reduction} of $(M,\omega)$ along $Y$.
\end{defn}

Boundary reduction can be viewed as \lq\lq half\rq\rq of a symplectic
cut (introduced by Lerman~\cite{Lerman.symplecticcuts}).
Its essence is captured in the following two-dimensional
example:
\begin{example} \label{cylinder.ex}  
Consider a cylinder with linear coordinate $0\le p\le 1$,
cyclic coordinate $q$ of period $2\pi$ and symplectic form $dp\wedge dq$.
Let $D$ be a disk of radius $\sqrt{2}$ in $(\bR^2,dx\wedge dy)$.
Then the change of coordinates
\begin{equation} \label{coord.eq}
x=\sqrt{2p} \cos(q), 
\qquad y=\sqrt{2p} \sin(q)
\end{equation}
defines a symplectomorphism from the cylinder minus the boundary component
defined by $p=0$ to the punctured disk.
Boundary reduction of the cylinder along the $p=0$ boundary component 
yields the disk $(D,dx\wedge dy)$.
\end{example}

Notice that a hypersurface in $(\bR^n\times T^n,\omega_0)$ that projects
to a hypersurface in $\bR^n$ is fibered by circles in the kernel of the
symplectic form only if the hypersurface in $\bR^n$ is affine planar.
In higher dimensions one can perform boundary reduction along
different smooth components of a boundary simultaneously.
The following example reveals the local convexity of (toric) moment 
map images as well as the properties listed in 
Proposition~\ref{toricbasics.prop}.

\begin{example} \label{bdyreduc.ex}  
For $0\le k\le n$, let
$U_k=\{(p_1,\ldots, p_n) | p_j\ge0, j\le k\} \subset (\bR^n,\cA_0)$ and 
let $P_j$ be the hyperplane through the origin defined by $dp_j=0$.
Let $(M,\omega) = (U_k\times T^n,dp\wedge dq)$ 
so we have a Lagrangian fibration $\pi:(M,\omega)\rightarrow (U_k,\cA_0,\cS)$
in which the stratification is trivial -- all fibers are top dimensional.
Then the boundary reduction of $(M,\omega)$ along
$\pi^{-1}(\cup_{j=1}^k P_j\cap U_k)$ is 
$(M',\omega')$ where
$M' =\bR^{2k}\times (\bR\times S^1)^{n-k}$
with symplectic form $\omega'=\sum_{i=1}^k dx_i\wedge dy_i+
\sum_{i=k+1}^n dp_i\wedge dq_i$.
The torus action on $(M,\omega)$ given by $t\cdot(p,q)=(p,q+t)$ descends to 
the action
$t\cdot(x,y,p,q)= (x \cos( t)-y\sin( t), x\sin( t)+y\cos(t),p,q)$ on 
$(M',\omega')$.
This yields the toric fibration
$\pi':(M',\omega')\rightarrow (U_k,\cA_0,\cS')$ where $\cS'$ is the 
stratification such that $p$ belongs to the $k$-stratum if it is in
the intersection of precisely $n-k$ of the $P_j$, $1\le j\le k$.
\end{example}

\begin{obs} \label{normalform.obs}  
A neighborhood of $\pi^{-1}(0)$ in
Example~\ref{bdyreduc.ex} provides a normal form for a neighborhood of
a $k$-dimensional fiber of a toric fibration.  (Compare~\cite{Audin.torus}).
\end{obs}

In working with manifolds that arise via boundary reduction it is helpful
to have the following two definitions:
\begin{defn}   \label{reducedbdy.def}
The  {\it reduced boundary} of a toric base
$(B,\cA,\cS)$, denoted by $\bdy_R B$,
is the set of points in $B$ that belong to lower-dimensional strata, i.e. 
$k$-strata for $k<n$.  (Notice that a point can belong to a lower-dimensional stratum only if it belongs to $\bdy B$.) 
\end{defn}
\begin{defn}   \label{bdyrecovery.def}
Given a toric fibration $\pi:(M,\omega)\rightarrow (B,\cA,\cS)$, the
{\it boundary recovery} is the unique 
regular Lagrangian fibered manifold $(B\times T^n,\omega_0)$ 
that yields $(M,\omega)$ via boundary reduction.
\end{defn}

\begin{rmk}   \label{mfldwbdy.rmk}
A local model for the neighborhood of a point in the boundary of a 
toric manifold can be extracted from Example~\ref{bdyreduc.ex}.
Simply let $W$ be a top dimensional submanifold with boundary in
$\bR^n$ such that the origin belongs to $\bdy W$ and $\bdy W$ is
transverse to the $\{P_j\}_{j=1,\ldots k}$; then take
$B=W\cap U_k$
and perform the boundary reduction along $\pi^{-1}(\cup_{j=1}^k P_j\cap B)$.
Notice that the image of the boundary of the resulting manifold 
is the closure of $\bdy B-\bdy_R B$ and $B_0=B-\bdy_R B$.
\end{rmk}

The geometric structure of $B$ and the stratification $\cS$
together determine what circles get collapsed
during boundary reduction of $(B\times T^n,\omega_0)$.    
\begin{defn}   \label{collapsing.def}
Let $\pi:(M,\omega)\rightarrow (B,\cA,\cS)$ be a toric fibration and
$\gamma$ a compact embedded curve with one endpoint $b_1$ in the
$(n-1)$-dimensional stratum of $\bdy_R B$ 
and such that $\gamma-\{b_1\}\subset B_0$.  
Let $b_0$ be the other endpoint of $\gamma$. 
The {\it collapsing class}, with respect to $\gamma$, 
for the smooth component of $\bdy_R B$ containing $b_1$ is the 
primitive class $\mathbf a\in H_1(F_{b_0},\bZ)$ that spans the 
kernel of $\iota_*:H_1(F_{b_0},\bZ)\rightarrow
H_1(\pi^{-1}(\gamma),\bZ)$ where $\iota$ is the inclusion map.
Accordingly,  the {\it collapsing covector} with respect to $\gamma$ is
the primitive covector $v^*\in T^*_{b_0}B$  
that determines vectors
$v(x)\in T_x^{\rm vert}M$ for each $x\in\pi^{-1}b_0$ such that the
integral curves of this vector field represent $\mathbf a$.
(Notice that $\mathbf a$ and $v^*$ are well defined up to sign.)
\end{defn}
\begin{rmk}  The collapsing class and covector are independent of $\gamma$
in a toric fibration, but will not be when this definition is extended to
almost toric fibrations.
\end{rmk}
\begin{rmk}   \label{collapsingvector.rmk}
We can also define the collapsing covector in $T^*_{b_1}B$ as the parallel
transport of $v^*$ along $\gamma$.
\end{rmk}

Now we can give a characterization of toric bases that yields a simple
topological construction of all toric manifolds.
\begin{thm} \label{topologicaltoric.thm}
A triple $(B,\cA,\cS)$ is a toric base if and only if there is
an immersion $\Phi:(B,\cA)\rightarrow (B,\cA_0)$ and the
stratification $\cS$ is such that
\begin{itemize}
\item each smooth component of $\bdy_R B$ is affine planar,
\item $(B,\cA)$ is locally convex in a neighborhood of any point
of $\bdy_R B$,
\item each component $K$ of the $k$-stratum of $\bdy_R B$, $k\ge 2$, 
is the intersection
of the closures of $k$ smooth components of the $(n-1)$-stratum,
\item the $k$ primitive integral covectors defining those planes span
a $k$-dimensional sublattice of $\Lambda|_K$ and
\item any point in the boundary of a $k$-stratum of $\bdy_R B$ has
a neighborhood isomorphic to a model neighborhood as described in
Remark~\ref{mfldwbdy.rmk}.
\end{itemize}
%\item every point on the interior of $B$ or the interior of $\bdy_R B$ 
%has a neighborhood locally
%isomorphic to a neighborhood of the origin in Example~\ref{bdyreduc.ex}
%for some $k$ and
%\item each point in the closure of $\bdy B-\bdy_R B$ has a neighborhood
%as described in Remark~\ref{mfldwbdy.rmk}.
%\end{itemize}
Furthermore, the boundary reduction of $(B\times T^n,\omega_0)$ according
to $\cS$ is the unique toric manifold whose base is $(B,\cA,\cS)$
\end{thm}

\begin{rmk}   Bear in mind that in general $B$ is not 
the image of the moment map.
\end{rmk}

\begin{proof}
Given any toric manifold $(M,\omega)$, the moment map $\mu$
factors through the orbit space yielding the map $\Phi:B\rightarrow \bR^n$
and an induced affine structure $\cA=\Phi^*\cA_0$.
We have already noted that Example~\ref{bdyreduc.ex} and 
Remark~\ref{mfldwbdy.rmk} provide normal forms for the fibration;
therefore they provide normal forms for the base.
The listed conditions are all properties implied by these normal forms.

Conversely, the immersion $\Phi$ defines a toric fibration of
$(B\times T^n,\omega_0)$ given by projection to the first factor.
The first hypothesis about the base $(B,\cA,\cS)$ guarantees that we
can perform boundary reduction according
to the stratification $\cS$ and the remaining hypotheses guarantee that
the resulting manifold is smooth.
The uniqueness of the toric manifold fibering over $(B,\cA,\cS)$ follows
from the Boucetta and Molino's classification of locally toric 
fibrations~\cite{BoucettaMolino.fibrations}.
\end{proof}

Since we will be concerned for the most part with bases that are
two-dimensional, to simplify exposition we use the following natural
language:
an {\it edge} is the closure of a component of the
$1$-stratum and a {\it vertex} is a component of the $0$-stratum.

While a toric base $(B,\cA,\cS)$ gives a well-defined toric fibration,
to specify a torus action one needs to choose an affine immersion
$\Phi$, at least up to translation.
The moment map then is $\Delta=\Phi(B)$.
Meanwhile, an immersion $\Phi$ is useful for visualizing $(B,\cA)$
as a subset of $(\bR^n,\cA_0)$.

We will use a bit more information than just the image of $\Phi$:
A {\it base diagram} for a Lagrangian fibered four-manifold
$\pi:(M,\omega)\rightarrow (B,\cA,\cS)$ is a subset 
$\Delta=\Phi(B,\cA)\subset(\bR^2,\cA_0)$ , for some choice of $\Phi$,
together with coding that enables reconstruction of the base
$(B,\cA,\cS)$ from $\Delta$.
For instance, for a toric base we use
heavy and light lines to represent reduced and non-reduced
parts of the boundary respectively.
Then, unless the base 
multiply covers of enough of the image in $(\bR^2,\cA_0)$ so as to obscure 
the structure of the preimage, we can construct the fibered manifold by
boundary reduction.
(If the immersion does obscure the structure of the base 
one would merely have to add extra information to the base diagram.)  
To do this, we need to know the collapsing classes.
They can be determined from the base diagram as follows,
assuming for simplicity of notation that $(B,\cA,\cS)\subset(\bR^2,\cA_0)$
so $\cA=\cA_0$:
\begin{obs}   \label{basecollapsing.obs}
If $v$ is a primitive integral vector that is normal to a smooth component
of the reduced boundary of $(B,\cA_0,\cS)\subset(\bR^2,\cA_0)$ then
$v^*=<v,\cdot>$ is a collapsing covector for that component.
(Here we are using the standard inner product on $\bR^2$ to provide an
isomorphism between its tangent and cotangent bundles.)
\end{obs}

Of course, we can build up new toric bases from existing ones by identifying
open subsets in each that are isomorphic.  Furthermore, it is easy to determine
when bases can be joined along smooth components of their boundaries in such
a way that the affine structures agree.  We are most interested in
when we can \lq\lq see\rq\rq\  the isomorphism, therefore with base
diagrams in mind we note:
\begin{prop} \label{basepasting.prop}
Consider two toric bases $(B_i,\cA_0,\cS_i)\subset (\bR^n,\cA_0)$. 
Suppose there is an element $(A,u)\subset AGL(n,\bZ)$ such that on 
the intersection $B_1\cap (AB_2+u)$ the $n$-stratum is non-empty,
the stratifications $\cS_1,\cS_2$ agree on the intersection
thereby defining $\cS_3$ on $B_1\cup (AB_2+u)$, and  
each connected component of the strata of $\cS_3$ is smooth. 
Then $(B_1\cup (AB_2+u),\cA_0,\cS_3)$ is a toric base.
\end{prop} 
\begin{proof}
The hypotheses of Theorem~\ref{topologicaltoric.thm} remain satisfied because
the union of two affine planar submanifolds of the same dimension
(that have non-empty intersection)
can be smooth only if the union itself is affine planar.
The condition that the $n$-stratum be non-empty then ensures that the
convexity condition is preserved.  
\end{proof} 

\subsection{Examples}
In this section, we present examples of toric fibered four-manifolds 
$(M,\omega)\rightarrow (B,\cA,\cS)$ via base
diagrams corresponding to immersions  $\Phi:(B,\cA)\rightarrow (\bR^2,\cA_0)$.
Accordingly, the preimage (under $\Phi\circ\pi$)
in $M$ of a point in the interior of a light 
edge is a 
$2$-torus while the preimage of a point in the interior of a heavy edge is a 
circle.
Similarly, a vertex can have a preimage in $M$ that is a torus, circle or 
point depending upon whether it is the intersection of
two light lines, one light and one heavy, or two heavy lines.

\begin{figure}
        \psfragscanon
        \psfrag{a}{(a)}
        \psfrag{b}{(b)}
        \includegraphics[scale=.75]{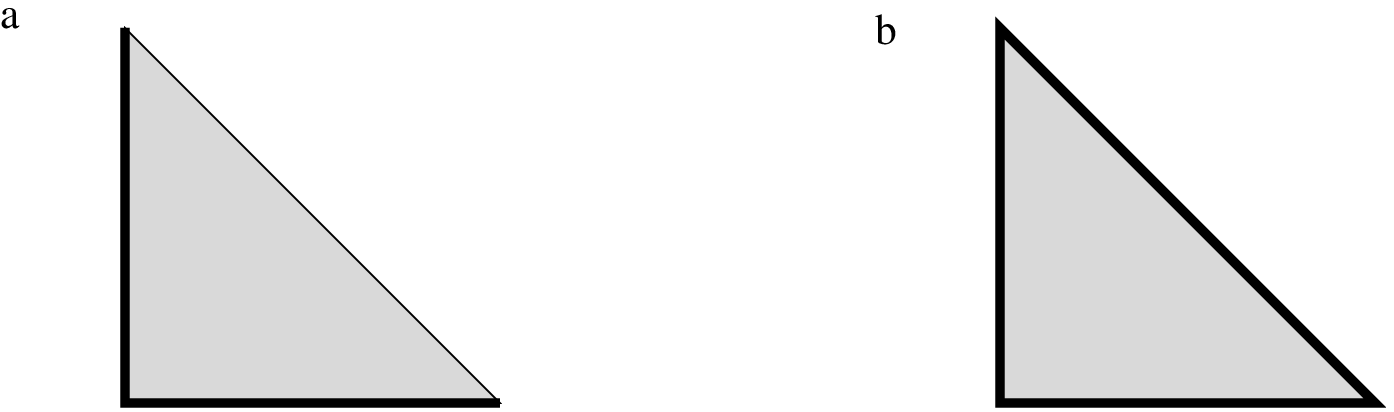}
\caption{(a) $B^4$ \ \ (b) $\bC P^2$}
\label{B4CP2.fig}
\end{figure}

\begin{example}  \label{B4.ex}  
A toric fibration of a four-ball $B^4$ of radius $\sqrt{2}$ in
$(\bR^4,dx\wedge dy)$ is given by
$\pi(x,y)= (\frac{1}{2}(x_1^2+y_1^2),\frac{1}{2}(x_2^2+y_2^2))$ 
whose image is the triangle with vertices $(0,0)$, $(1,0)$ and $(0,1)$.
See Figure~\ref{B4CP2.fig}(a).
Notice that if $B^4$ is identified with a ball in $\bC^2$ by
$z_j=x_j+iy_j$, the preimage of the horizontal and vertical edges are
the disks in complex coordinate axes $z_2=0$ and $z_1=0$ respectively.
\end{example}

\begin{example} \label{CP2.ex}  
The triangle with vertices
$(0,0)$, $(1,0)$ and $(0,1)$ shown in Figure~\ref{B4CP2.fig}(b) represents
a Lagrangian fibered symplectic manifold diffeomorphic to $\bC P^2$.
This reflects the fact that $\bC P^2$ is a boundary reduction of
$(B^4,\omega_0)$ in which the circles of the Hopf fibration are collapsed
to points.
\end{example}

Figure~\ref{B4CP2.fig}(b) also illustrates the standard construction
of $\bC P^2$ from three affine coordinate charts:
\begin{xca}  
Use Example~\ref{B4.ex} and Proposition~\ref{basepasting.prop} to find
the fiber-preserving symplectomorphisms that yield 
the Lagrangian fibered $\bC P^2$ defined by Figure~\ref{B4CP2.fig}(b)
as the union of three open Lagrangian fibered four-balls.
\end{xca}

\begin{figure}
\begin{center}
	\includegraphics[scale=.75]{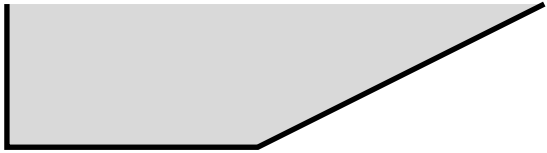}
\end{center}
\caption{Neighborhood of a sphere: $U_{2,1}$}
\label{spherenbhd.fig}
\end{figure}

\begin{example} \label{spherenbhd.ex}  
For each $n\in \bZ$ consider the domain 
\begin{equation} \label{Un1.eqn}
U_{n,1}:=\{(p_1,p_2)|p_1\ge 0, 0\le p_2<t, np_2\ge p_1-s\}.
\end{equation}
for some $s,t>0$.
Figure~\ref{spherenbhd.fig} shows $U_{2,1}$.
Then $U_{n,1}$  defines a toric fibration of a neighborhood
of a sphere -- the sphere that is the preimage of the heavy horizontal
line segment.
The sphere has coordinates $(p_1,q_1)$ on the complement of two
points (the preimages of the vertices), so with respect to the symplectic
form $dp\wedge dq$ it clearly has (symplectic) area $2\pi a$.
Furthermore, the sphere has self-intersection $-n$.
There are several ways to see this.  
For instance, one can check that the fiber-preserving symplectomorphism   
needed to present it as the union of two Lagrangian fibered four-balls
implies it is diffeomorphic to a disk bundle over the sphere with 
Chern class $-n$.
In Section~\ref{selfint.sec} we explain an alternative way to see
this that generalizes beyond spheres that are invariant under a
torus action.
\end{example}

\begin{figure}
\begin{center}
	\includegraphics[scale=.75]{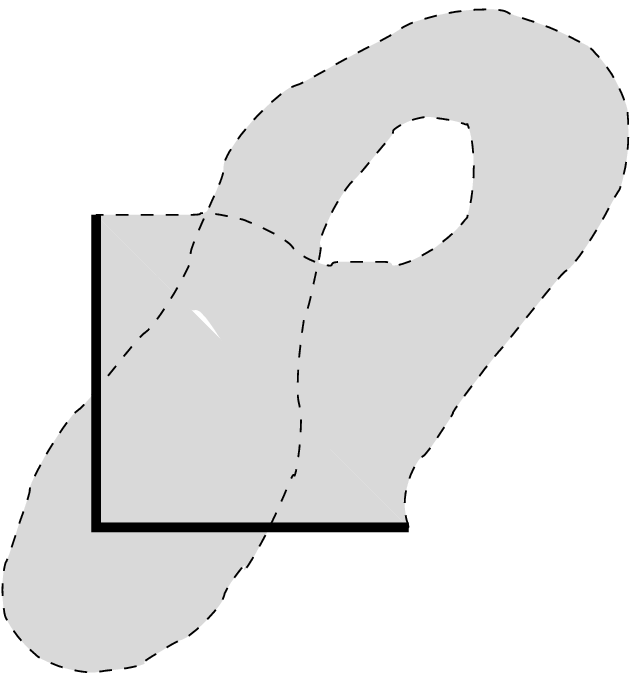}
\end{center}
\caption{Exotic $\bR^4$}
\label{exoticR4.fig}
\end{figure}

\begin{example}   \label{exoticR4.ex}
Figure~\ref{exoticR4.fig} is a base diagram for a toric fibration
of an exotic symplectic $\bR^4$, i.e. $(\bR^4,\omega)$ where $\omega$ 
is exotic in the sense that there
is no symplectic embedding of $(\bR^4,\omega)$  into 
$(\bR^4,\omega_0)$ where $\omega_0$ is the standard symplectic structure.
This example, due to Zung~\cite{Zung.II}, 
generalizes in an obvious way to higher dimensions
yielding an exotic $\bR^{2n}$ for any $n\ge2$.

To see that $\omega$ is exotic observe that
the immersion of the base into $(\bR^2,\cA_0)$ 
gives local coordinates $(p,q)$ on the complement
of the reduced boundary with respect to which the symplectic structure
is $dq\wedge dq$ with primitive $\alpha=pdq$.
Placing the vertex of the base diagram at the origin, the regular torus that
projects to the origin is an exact Lagrangian torus (i.e. $\alpha$ restricted
to this torus is cohomologous to $0$, in this case actually equal to
$0$).  Gromov~\cite{Gromov.pseudohol} proved that there are no exact
Lagrangian tori in $(\bR^{2n},\omega_0)$, thereby establishing that
this $\bR^4$ is symplectically exotic.
\end{example}

\section{Singularities} \label{singularities.sec}
\subsection{Integrable systems and non-degenerate singularities}

When we pass from regular Lagrangian fibrations to those with singularities
we need to restrict the class of allowable singularities in order to
keep control of a relation between the geometry of the
base and the topology of the total space. 
A natural place to turn for a classification of singularities of a Lagrangian
fibration is the theory of integrable systems.  

\begin{defn} \label{intsyst.def}  
An {\it integrable system} on a symplectic manifold $(M^{2n},\omega)$
is a collection of $n$ functionally independent
Poisson commuting functions $F_i:M\rightarrow \bR^n$.
\end{defn}
An equivalent definition is a Lagrangian fibration whose base can be
immersed in $\bR^n$ (ignoring affine structures).

Abusing slightly the integrable systems definition of non-degenerate
singularities, we have

\begin{defn}   \label{nondegen.def}
A critical point of a Lagrangian fibration 
$\pi:(M^{2n},\omega)\rightarrow B^n$ 
is called a {\it non-degenerate singular point of rank $k$}
if it has a Darboux
neighborhood (with symplectic form $dx\wedge dy$) 
in which the projection is 
$\pi=(\pi_1,\ldots \pi_k,\pi_{k+1},\ldots \pi_{n})$
where $\pi_j(x,y)=x_j$ for $j\le k$ and the other components
have some combination of the following forms:

\begin{equation} \label{elliptic.eq}
\pi_j (x,y)=(x_{j}^{2}+y_{j}^{2}) \quad \textit{elliptic}
\end{equation}

\begin{equation} \label{hyperbolic.eq}
\pi_j (x,y)=(x_{j}^{2}-y_{j}^{2}) \quad \textit{hyperbolic}
\end{equation}

\begin{equation} \label{nodal.eq}
(\pi_i,\pi_j) (x,y)= (x_{i}y_{i}+x_{j}y_{j},
         x_{i}y_{j}-x_{j}y_{i}) \quad \textit{nodal, 
{\rm or} focus-focus}.
\end{equation}

\end{defn}
In other words, a non-degenerate singularity splits as a product 
of elliptic, hyperbolic and nodal singularities.
Notice that the elliptic singularities are precisely the singularities
that occur in a toric fibration.

The nodal singularity is a complex hyperbolic singularity:
in terms of complex coordinates
\begin{equation} \label{complex.eq}
{\rm x}=x_1+ix_2, \ {\rm y}=y_1+iy_2, \quad 
\omega={\rm Re}\ {d{\rm \overline x}\wedge d{\rm y}}
\end{equation}
the projection from a four-manifold to $\bC$ is  
\begin{equation} \label{complexnodal.eq}
\pi({\rm x,y)=\overline{x}y}.
\end{equation}
This is the same singularity that appears in {\it Lefschetz fibrations}
of symplectic manifolds (cf.~\cite{GompfStipsicz.4mflds}).

\begin{rmk}
The precise definition of a non-degenerate singular
point is as follows.
Given an integrable system $F=(F_1,\ldots,F_n):(M^{2n},\omega)
\rightarrow \bR^n$, a point $x\in M$ such that 
$DF(x)=0$  is called a {\it non-degenerate fixed point}
if there are
symplectic coordinates near $x$ with respect to which the quadratic
parts of the functions $F_i$ form a Cartan subalgebra of the
Lie algebra of quadratic functions on $\bR^{2n}=T_xM$ with respect to
the Poisson bracket.
To generalize to singular points $x$ of rank $k$ 
(where the image of $DF$ has dimension $k$)
replace $T_{x}M$ with the symplectic vector
space ${\rm Ker}(DF(x))/\cI$ where
$\cI$ is the subspace of $T_xM$ spanned by the Hamiltonian vector
fields $X_i$ of the functions $F_i$ (so $\iota_{X_i}\omega=-dF_i$). 
Eliasson~\cite{Eliasson.normalforms} proved that 
such singularities have neighborhoods that are fiberwise diffeomorphic 
to the singularities of linearized 
integrable systems -- which in turn were classified by 
Williamson~\cite{Williamson.normalforms}.
Since we are interested in questions of global symplectic topology rather
than the intricacies of a given integrable system, our abuse is to assume
the singularity has a neighborhood fiberwise symplectomorphic
to a linearized integrable system.
\end{rmk}

The structure imposed on a fibration by the presence of
a hyperbolic singularity is fundamentally different
from the structure imposed by elliptic or nodal singularities. 
The key features of the fibration near an elliptic or nodal singularity
are that the base of the fibration is a manifold 
and the structure of the singular fiber is very simple.
For instance, an elliptic singular point in a two-manifold has a 
neighborhood that is, 
on the complement of the singular point, fibered by circles
that project to a half-open interval in the base.
In dimension four a nodal singular point 
has a neighborhood that is fibered by cylinders and two 
transversely intersecting disks; since this singular point
is isolated the cylinders belong to tori, thereby forcing
the singular fiber to be a pinched torus (if it has no other
singular points).
If there are $k\ge 2$ nodal singular points on a fiber then
the fiber is a union of $k$ spheres, each of which intersects
two others transversely (and positively) at one point.
In either case, the base of the fibration of the neighborhood is an open disk.

In contrast, the presence of hyperbolic singularities in a fibration
causes the base to become non-smoothable and allows the topology of the
singular fibers to become complicated.
For instance, in dimension two the base for a fibered neighborhood of 
a hyperbolic singularity is an \lq\lq X\rq\rq\ and recovering the 
symplectic manifold from the base becomes
impossible without extra data:
several hyperbolic singular points and the arcs connecting them are
typically mapped to one singular point in the base so the topological
type of the singular fiber would in general have to be stipulated.

Note that in dimension two the only manifolds that admit a
Lagrangian fibration with no hyperbolic singularities are the sphere and
the torus.  
However, all other (orientable) surfaces can be fibered
once one allows hyperbolic singularities. 
This prompts the following:
\begin{question}   \label{fourmfld.quest}
What symplectic four-manifolds admit Lagrangian fibrations with
topologically stable non-degenerate singularities?
\end{question}

In this paper our focus is on using base diagrams to determine symplectic
manifolds.  Therefore we henceforth exclude hyperbolic singularities
and define
\begin{defn}   \label{almosttoric.def}
A non-degenerate Lagrangian fibration $\pi:(M,\omega)\rightarrow
B$ of a symplectic
four-manifold is an {\it almost toric fibration} if it is a
non-degenerate topologically stable fibration with no hyperbolic 
singularities.  A triple $(B,\cA,\cS)$ is an {\it almost toric base}
if it is the base of such a fibration.
A symplectic four-manifold equipped with such a fibration is an
{\it almost toric manifold}.
\end{defn}
In the above definition $\cA$ is the affine structure on $B-\{s_i\}$ where
the $s_i$ are the images of the nodal singularities, and $\cS$ is defined
so that the $s_i$ belong to the top-dimensional stratum.

The definition clearly generalizes to other dimensions,
however we restrict our discussion to dimension four. 
Notice that an almost toric manifold can fail to be toric in one of
two ways: 1) it has nodal fibers or 2) there is no affine immersion of
the base into $(\bR^2,\cA_0)$.

\subsection{Neighborhood of a nodal fiber } \label{nodalnbhd.sec}

\begin{defn}   \label{nodalfiber.def}
A {\it nodal fiber} in a Lagrangian fibered four-manifold
is a singular fiber with nodal singularities and  
a {\it node} is the image of such a fiber under the fibration map.
\end{defn}
Generically, nodal singularities of a Lagrangian fibration of a four-manifold
occur in distinct fibers.  Furthermore,  
one can always vary a fibration locally so that each fiber has only
one singularity.  (See Exercise~\ref{multiplicityknodes.ex}.) 
Therefore,
we in general assume a nodal fiber has just one nodal singularity and
specify that the nodal fiber and the corresponding node in
the base are multiplicity $k$ whenever there are $k\ge 2$ singular points
on the fiber.

Topologically, a nodal fiber is the singular fiber of a Lefschetz
fibration (in which it appears in a family of symplectic fibers).
Specifically, a neighborhood of a nodal fiber is diffeomorphic to
$T^2\times D^2$ with a $-1$-framed $2$-handle attached along a
simple closed curve in $T^2\times \{x\}$ for some $\{x\}\in\bdy D^2$.
This is explained in~\cite{GompfStipsicz.4mflds}, for instance.
It also follows from the techniques of this paper: the core
of the $2$-handle is the vanishing disk (the disk whose boundary represents
the vanishing class, defined in Section~\ref{nodalmonodromy.sec})
and the $-1$-framing is calculated in Lemma~\ref{nodalintersec.lem}.

An example of a nodal fiber is the image of the zero section
of a disk bundle in $T^*S^2$ under a self-plumbing of the disk
bundle.  
To construct this, first rewrite the normal form for the nodal singularity in
terms of complex coordinates ${\rm x}=x_1+ix_2, \ {\rm y}=y_1+iy_2$
with $\omega_0={\rm Re}\ {d\,{\rm \overline x}\wedge d{\rm y}}$.
Thinking of the cotangent bundle as a complex line bundle with 
Chern class $c_1=-2$, we build a disk bundle using two coordinate
charts 
\begin{eqnarray}
U:= & \{({\rm x},{\rm y})|\ \abs{\rm xy}<\epsilon<1, 
\abs{\rm x}<2,\abs{\rm y}<\frac{1}{4} \} \quad and \\ 
V:= & \{({\rm u},{\rm v})|\ \abs{\rm uv}<\epsilon<1,
\abs{\rm u}<2,\abs{\rm v}<\frac{1}{4}\}
\end{eqnarray} 
with ${\rm u}= u_1+iu_2$
and ${\rm v}=v_1+iv_2$.
These are glued together via the map
\begin{equation} \label{transition1.eq}
({\rm u},{\rm v})=\varphi({\rm x},{\rm y})
=\left(\frac{1}{\rm x},\overline{\rm x}^2 {\rm y}\right)
\end{equation}
defined on $U\cap \{\frac{1}{2}<\abs{\rm x}<2\}$.
This map is clearly a symplectomorphism if we equip the second
coordinate chart with the symplectic form 
${\rm Re}(-d\,{\rm \overline u}\wedge d{\rm v})$.

Now construct a self-plumbing of the disk bundle by identifying
neighborhoods of the origin in each coordinate chart, say
$\{(\rm{x},{\rm y})| \ \abs{\rm x}, \abs{\rm y} <\frac{1}{4}\}$ and
$\{(\rm{u},{\rm v})| \ \abs{\rm u}, \abs{\rm v} <\frac{1}{4}\}$,
via the map 
\begin{equation} \label{transition2.eq}
({\rm u},{\rm v})=\psi({\rm x},{\rm y})
=\left(\overline{\rm y},\overline{\rm x}\right)
\end{equation}
which again preserves the symplectic form.

The Lagrangian fibration given by 
$\pi({\rm x},{\rm y})=\overline {\rm x}{\rm y}$,
$\pi({\rm u},{\rm v})=\rm \overline{u}v$ then has one singular point
at $({\rm x},{\rm y})=({\rm u},{\rm v})=(0,0)$.

\begin{xca}   \label{sign.ex}
Check that the singular fiber is a sphere with one positive
self-intersection.  (When checking the positivity of the intersection,
note that the manifold is oriented by  
$\omega^2=dx_1\wedge dy_1\wedge dx_2\wedge dy_2=dx_1\wedge dx_2
\wedge dy_1\wedge (-dy_2)$.)
\end{xca}

\begin{prop} \label{nodalgerm.prop}
The germ of a neighborhood of an isolated nodal fiber in a Lagrangian fibered
symplectic manifold is unique up to symplectomorphism.  Furthermore,
given any two symplectomorphic fibered neighborhoods of a nodal fiber the 
symplectomorphism between them can be chosen to be fiber-preserving on
the complement of smaller fibered neighborhoods.  (See~\cite{Symington.grbd}.) 
\end{prop}
It is not hard to generalize the above proposition to include
multiplicity $k$ nodal fibers (in which case $k$ is of course a topological
invariant of the neighborhood).

\begin{rmk}   \label{San.rmk}
Note that the germ of a fibered neighborhood of a nodal fiber is
not unique up to fiber-preserving symplectomorphism; San~\cite{San.focusfocus}
has identified a 
non-trivial invariant that classifies the germs of such neighborhoods
up to fiberwise symplectomorphism.
\end{rmk}

\subsection{Monodromy around a nodal fiber} \label{nodalmonodromy.sec}

\begin{figure}
\begin{center}
        \psfragscanon
        \psfrag{c1}{$C_1$}
        \psfrag{c2}{$C_2$}
        \psfrag{g1}{$\gamma_1$}
        \psfrag{g2}{$\gamma_2$}
        \psfrag{r=0}{$r=4\delta$}
        \psfrag{r=2}{$r=2$}
        \psfrag{rho=0}{$\rho=4\delta$}
        \psfrag{rho=3}{$\rho=2$}
        \psfrag{theta}{$\theta$}
        \psfrag{alpha}{$\alpha$}
        \psfrag{beta}{$\beta$}
	\includegraphics{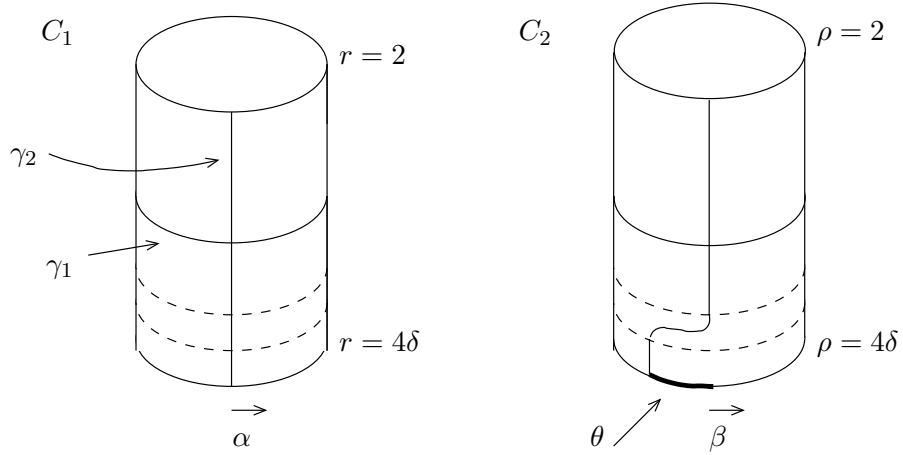}
\end{center}
\caption{Positive Dehn twist as $\theta$ goes from $0$ to $2\pi$}
\label{Dehntwist.fig}
\end{figure}

We can detect the monodromy around a nodal fiber from the local
model in Section~\ref{nodalnbhd.sec}.
Indeed, continuing with the complex variables notation, 
restrict to the torus bundle over a circle of radius $\delta<\epsilon$
in $\bC=\bR^2$ and let $\theta=\arg(z)$ parameterize the circle in the
base, so $\theta=\arg(\overline{\rm x}{\rm y})=\arg(\overline{\rm u}{\rm v})$.
Then the intersection of the torus bundle over the circle 
with the two coordinate charts $U,V$ gives
two cylinder bundles that can be trivialized by coordinates
$(r,\alpha)$ and $(\rho,\beta)$ where
${\rm x}=r e^{ i\alpha}$ and ${\rm u}=\rho e^{i\beta}$.
The restrictions of the maps $\varphi$ and $\psi$ show that these
bundles are glued together by 
$\rho e^{i\beta}=\frac{1}{r} e^{-i\alpha}$ when $\frac{1}{2}<r<2$ 
and
$\rho e^{i\beta}=\frac{ \delta}{ r}  e^{-i(\alpha+\theta)}$ when
$4\delta<r<\frac{1}{4}$.
Figure~\ref{Dehntwist.fig} shows the two cylinders that make up a fiber
torus $F_b$ over $b=z=\delta e^{i\theta}$.  
On this torus two curves are indicated: $\gamma_1$ given by $r=\rho=1$
and $\gamma_2$ given by $\alpha=0$ in $U$.
As $\theta$ varies the curves $\gamma_1,\gamma_2$ give a trivialization
of the torus bundle over the interval $0\le \theta<2\pi$.
Meanwhile, $\gamma_2$ shows
the positive Dehn twist in the torus fiber as $\theta$
varies from $0$ to $2\pi$ and reveals that the 
topological monodromy around a nodal fiber is
\begin{equation} \label{A10.eq}
\Psi(\gamma) = A_{(1,0)}= \left(  
\begin{matrix}
1 & 1  \\  
0 & 1 
\end{matrix} 
\right)
\end{equation}
with respect to a basis $\{[\gamma_1],[\gamma_2]\}$ 
for $H_1(F_b,\bZ)$.

Changing the basis for $H_1(F_b,\bZ)$ causes the image of $\Psi$ 
to be a conjugate of $A_{(1,0)}$.
Such matrices are of the form
\begin{equation} \label{Aab.eq}
A_{(a,b)}=\left(  
\begin{matrix}
1-ab & a^{2} \\  
-b^{2} & 1+ab 
\end{matrix} 
\right)  
\end{equation} 
with eigenvector $(a,b)$.

\begin{def}   \label{vanishing.def}
Consider an almost toric fibration $\pi:(M,\omega)\rightarrow B$ that
has a node at $s$.
Let $\eta$ be an embedded curve with endpoints at $s$ and a regular
point $b\in B_0$ such that $\eta-\{s\}\subset B_0$ contains no other nodes.
Associated to $s$ and  $\eta$ is the (well defined up to sign)
{\it vanishing class} in $H_1(F_b,{\bZ})$, namely 
the class whose representatives bound a disk in $\pi^{-1}(\eta)$.
(With respect to the language for Lefschetz fibrations, 
the vanishing class is the homology class of the
vanishing cycle.)
The {\it vanishing covector} 
$w^*\in T^*_b B$ is the primitive covector (defined up
to sign) that determines vectors
$w(x)\in T_x^{\rm vert}M$ for each $x\in\pi^{-1}b$ such that the
integral curves of this vector field represent vanishing class.
\end{def}

\begin{lemma} \label{vanishingclass.lemma}
Suppose $\gamma$ is
a positively oriented loop based at $b$ that is the boundary of 
a closed neighborhood of $s$ containing $\eta$.
Then the vanishing class is the unique class (up to scale) that is
preserved by the monodromy along $\gamma$.
With respect to a basis for $H_1(F_b,\bZ)$ for which the
monodromy matrix is $A_{(a,b)}$, 
the vanishing class is the class $(a,b)$.
\end{lemma}

\begin{proof}
Without loss of generality we work with the model neighborhood of
Section~\ref{nodalnbhd.sec} and with basis $\{[\gamma_1],[\gamma_2]\}$ 
for $H_1(F_b,\bZ)$, as above, so the topological
monodromy matrix is $A_{(1,0)}$.
Now note that $\gamma_1$ bounds a disk: the
disk $\{({\rm x},{\rm y})| \ {\rm arg}(\overline{\rm x}{\rm y})=0,
|{\rm y}|=|{\rm x}|\le \sqrt{\delta}\}$.
Thus the vanishing class is $(1,0)$, the eigenvector of the topological
monodromy.
\end{proof}

Because a nodal fiber of multiplicity $k$ (a fiber with $k$ nodal
singular points) arises from  coalescing
$k$ nodal fibers of multiplicity $1$, the monodromy around a nodal
fiber of multiplicity $k$ is, up to conjugation,
\begin{equation} \label{A10k.eq}
\Psi(\gamma) = A^k_{(1,0)}= \left(  
\begin{matrix}
1 & k  \\  
0 & 1 
\end{matrix} 
\right).
\end{equation}

\subsection{Affine structure near a node} \label{affinenode.sec}

Suppose $\pi:(M,\omega)\rightarrow B$ 
is a Lagrangian fibered neighborhood
of a nodal fiber in a four-manifold with $s\in B$ the node.
The restriction of the fibration to $B-\{s\}$ is a regular fibration
and hence $B-\{s\}$ has an induced affine structure $\cA$.
Because of the non-trivial affine monodromy around the node (which we
call the {\it nodal monodromy}) there is 
no immersion of $(B-\{s\},\cA)$ into $(\bR^2,\cA_0)$.
However, we can find an affine immersion 
if we remove a curve from $B$ that has an endpoint at the node,
leaving us with a simply connected surface.
Therefore, as a base diagram for the neighborhood of a node, we will
use the image of an affine immersion of the complement of such a curve.

\begin{figure}
\begin{center}
        \psfragscanon
        \psfrag{a}{(a)}
        \psfrag{b}{(b)}
	\includegraphics[scale=.5]{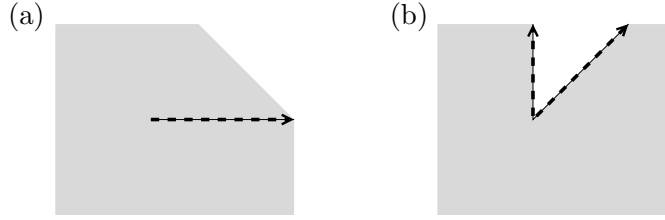}
\end{center}
\caption{Neighborhood of a node}
\label{nodenbhd.fig}
\end{figure}

For convenience we choose, whenever possible, the curve we remove 
to be a ray $R$ based at the node $s$.
Then, up to translation, any immersion 
$\Phi:(B-R,\cA)\rightarrow (\bR^2,\cA_0)$
has an image $\Phi(B-R)$ that is a neighborhood of the origin minus a sector.
Different choices of the ray $R$ cause the internal angle of the sector to vary
in the interval $[0,\pi)$.
Figures~\ref{nodenbhd.fig}(a) and (b) show base diagrams in which
the shaded part of the complement of the heavy dotted lines
is the image $\Phi(B-R)$ for two different choices of $R$.

As we will see below, the sector has angle zero when $R$ belongs to the
one well-defined affine linear submanifold through the node $s$:
\begin{defn}   \label{eigen.def}
The {\it eigenline} $L$ through a node $s$ is the unique maximal affine linear
immersed one-manifold through the node
for which there is a regular point $b\in L$, arbitrarily close to $s$,
such that the affine monodromy along an arbitrarily small loop around $s$ and 
based at $b$ preserves $T_bL\subset T_bB$.
An {\it eigenray} is either of the two maximal affine linear submanifolds
that has an endpoint at the node and is a subset of the eigenline.
\end{defn}
In Figure 4(a) the indicated ray is an eigenray.

We now construct a family of integral affine punctured
planes, $(V_u,\cA_u)$ where $V_u$ is diffeomorphic to $\bR^2-\{0\}$.
Let $(\widetilde W,\cA_0)$ be the universal cover of 
$(\bR^2-\{0\},\cA_0)$ with affine projection $\Psi$.
For each vector $u$ representing a point in $\bR^2$ 
let $\widetilde V_u$ be a sector in 
$\widetilde W$ with interior angle $>\pi$ and $\le 2\pi$ and 
that is bounded by rays $R,R'$ 
consisting of points $x,x'$ such that $\Psi(x)=\lambda u$
and $\Psi(x')=\lambda A_{(1,0)}u$.
(We are ignoring the ambiguity of how $\widetilde V_u$ is affinely embedded
in $\widetilde W$ since this is irrelevant for our construction.)
Then define $(V_u,\cA_u)$ where $V_u:= \widetilde V_u/x\sim x'$ for
$x,x'$ as above and $\cA_u$ is defined by the induced lattice $\Lambda_{u}$ on
$TV_u$:
\begin{equation}
(TV_u,\Lambda_{u})=(\widetilde V_u\times \bR^2,\Lambda_0)/
(x,v)\sim(x',A_{(1,0)}v).
\end{equation}
Of course, $(\widetilde V_u,\cA_u)$ and 
$(\widetilde V_{\lambda u},\cA_{\lambda u})$ are isomorphic for any
$\lambda>0$.
More importantly,

\begin{lemma} \label{nodemodels.lem}
For all $u\in\bR^2$, $(V_u,\cA_u)$ is isomorphic to $(V_{e_1},\cA_{e_1})$.
\end{lemma}
Therefore we simplify notation and call this affine manifold $(V,\cA)$.

\begin{xca}   Prove Lemma~\ref{nodemodels.lem}
\end{xca}

If a node has multiplicity $k$ then local models can be constructed
exactly as above with $A_{(1,0)}$ replaced by
$A^k_{(1,0)}$.  We call such affine punctured planes $(V^k,\cA^k)$.
Note that if $R$ is an eigenray, then the image $\Phi(B-R)$
is independent of $k$.

\begin{prop} \label{nodenbhd.prop}
Any node $s\subset B$ 
of multiplicity $k\ge 1$ has a neighborhood that is affine
isomorphic to a neighborhood of the puncture in $(V^k,\cA^k)$.
\end{prop}

\begin{proof}
Let $U$ be a neighborhood of $s$ and let $R\subset U$ be a ray based at $s$.
Then $(U-R,\cA)$ is homomorphic to a disk and hence there is an affine
immersion $\Phi:(U-R,\cA)\rightarrow (\bR^2,\cA_0)$ such that the
origin is the limit of the images of points in $B$ whose limit is $s$.
Corollary~\ref{onetoone.cor} implies that this is an embedding
 and
hence is isomorphic to a neighborhood of $r=0$ in the interior of 
$(\widetilde V_u,\cA_0)$ for some $u$.
The result then follows from Lemma~\ref{nodemodels.lem} and the
remarks about nodes of multiplicity greater than one.
\end{proof}
                                                                               
Furthermore, for each choice of ray $R\subset B$ based at $s$ there is
an immersion $\Phi:(B-R,\cA)\rightarrow (\bR^2,\Lambda_0)$ and a $u$
such that $\Phi(B-R)$ is the intersection of $\Psi({\rm Int}\widetilde V_u)$
and a neighborhood of the origin in $\bR^2$.
In particular, when $R$ is an eigenray, $\Phi(B-R)$ is the intersection of 
$\Psi({\rm Int}\widetilde V_{\pm e_1})$ and a neighborhood 
of the origin in $\bR^2$, as in Figure~\ref{nodenbhd.fig}(a).

\section{Almost toric manifolds}  \label{almosttoric.sec}

\subsection{Classification in terms of base}
\label{almosttoricfibr.sec}

The base of an almost toric fibration has a geometric structure very
similar to  a toric base.
Indeed, the only difference between an almost toric and a toric base
(on a local level) is the presence of nodes which in turn have 
neighborhoods whose affine structure (on the complement of the node)
depends only on the multiplicity of the node.

\begin{defn}   \label{affinewnodes.def}
An {\it integral affine manifold with nodes} $(B,\cA)$ is a two-manifold 
$B$ equipped with an 
integral affine structure on $B-\{s_i\}$ such that each $s_i$ has
a neighborhood $U_i$ such that $(U_i-s_i,\cA)$ is affine isomorphic
to a neighborhood of the puncture in $(V^k,\cA^k)$ for some $k\ge 1$.
\end{defn}

As for toric manifolds, we are interested in using the base of an
almost-toric fibration to understand topological and symplectic
features of the total space.
The stratification $\cS$ of a toric base is an important invariant of the
fibration and to extend this notion to the almost toric case we
declare that nodes belong to the top-dimensional stratum.
Then triples $(B,\cA,\cS)$ that can appear as bases of almost toric 
fibrations can be characterized as follows:

\begin{thm}   \label{almosttoricbase.thm}
Consider a triple $(B,\cA,\cS)$ such that $(B,\cA)$ is an integral
affine manifold with nodes $\{s_i\}_{i=1}^N$.
Then $(B,\cA,\cS)$ is an almost toric base if and only if 
every point in $B-\{s_i\}_{i=1}^N$ has a neighborhood that is a
toric base.
\end{thm}
\begin{proof} 
An almost toric fibration is by definition locally toric on the complement 
of its nodal fibers.
This implies that the base must, on the complement of the nodes,
have the local structure of a toric base.

Conversely, if the base is locally toric on the complement of the nodes,
we can construct an almost toric fibration
in a similar fashion to a toric fibration: 
Take the Lagrangian fibration of $((B-\{s_i\}_{i=1}^N)\times T^2,\omega_0)$ 
given by projection onto the first factor.
The local toric structure guarantees that boundary reduction to
achieve the stratification $\cS$ yields a toric fibered smooth manifold,
say $\pi':(M',\omega')\rightarrow (B-\{s_i\}_{i=1}^N,\cA,\cS)$.
(Recall that $\pi'$ does not define a toric fibration unless there is an
immersion $(B-\{s_i\}_{i=1}^N,\cA)\rightarrow(\bR^2,\cA_0)$).
Furthermore, the definition of an affine manifold with nodes implies that
each node has a neighborhood
that is the base of a model neighborhood of a nodal singular fiber.
Theorem~\ref{topologicaltoric.thm}
guarantees that the fibration over the regular points of this neighborhood
is unique.  Therefore there is a fiber-preserving symplectomorphism
of the complement of the singular fiber in its model neighborhood and
the restriction of our toric fibration
$\pi':(M',\omega')\rightarrow (B-\{s_i\}_{i=1}^N,\cA,\cS)$
to the preimage of a punctured neighborhood of the node. 
Using such fiber-preserving symplectomorphisms we can glue in nodal fibers
and achieve an almost toric fibration with base $(B,\cA,\cS)$.
\end{proof}

Our definition of almost toric manifolds (Definition~\ref{almosttoric.def})
is motivated by a classification of non-degenerate Lagrangian fibrations due
to Zung~\cite{Zung.II}.
To classify these fibrations up to fiber-preserving symplectomorphism
(or symplectic equivalence) he introduced the following notion: 
two non-degenerate Lagrangian fibrations over the same base $(B,\cA,\cS)$
are {\it roughly symplectically equivalent} if 
their singular fibers have fiberwise symplectomorphic neighborhoods.
Here, $\cA,\cS$ are the affine structure and stratification that Zung
showed exist on the base of a non-degenerate Lagrangian fibration.

Zung also defined an invariant of a non-degenerate Lagrangian fibration,
the {\it Lagrangian Chern class}.
It is an element of the first homology
of the base with values in the sheaf of closed basic $1$-forms
($1$-forms that vanish on vectors tangent to fibers) modulo
those forms that arise from contracting the vector field for
a symplectic fiber-preserving circle action with the symplectic form.
(To be precise, this Chern class is actually a relative class in the sense 
that
it is defined relative to a given reference fibration.)

\begin{thm}[Zung~\cite{Zung.II}] \label{Zung.thm}
Two non-degenerate Lagrangian fibrations over the same base $(B,\cA,\cS)$ 
are symplectically equivalent if and only if they are
roughly symplectically equivalent and have the same 
Lagrangian Chern class.
\end{thm}

In dimension four with the absence of hyperbolic singular points,
we are led to the following corollary:
\begin{cor} \label{almosttoricbase.cor}
Let $(B,\cA,\cS)$ be an almost toric base 
which is either non-compact or has non-empty boundary. 
Then $(B,\cA,\cS)$ is the base of an almost toric fibration of a unique
symplectic manifold.
\end{cor}

\begin{proof}
The hypothesis that the two-dimensional base be either 
non-compact or have non-empty boundary implies that it has the 
homotopy type of a one-dimensional space.  Therefore the Lagrangian
Chern class vanishes and hence all roughly symplectically equivalent
fibrations with base $(B,\cA,\cS)$ are fiberwise symplectomorphic.

Now suppose we have two almost toric fibrations over $(B,\cA,\cS)$.
Because of the hypothesis on the topology of $B$, we know that
there is a unique fibration over $(B-\{s_i\}_{i=1}^N,\cA,\cS)$.
Completing this to a fibration over $(B,\cA,\cS)$ as in the proof of
Theorem~\ref{almosttoricbase.thm} we are not guaranteed of getting
a unique Lagrangian fibration (see San~\cite{San.focusfocus}).  However
the proof of Proposition~\ref{nodalgerm.prop} given in~\cite{Symington.grbd}
shows that the structure of the fibration can be perturbed via a
compactly supported symplectic isotopy, thereby establishing uniqueness
up to symplectomorphism of the total space.
\end{proof}

\subsection{Base diagrams} \label{base.sec}

Our objective now is to draw base diagrams in $\bR^2$ that prescribe the
construction of almost toric bases $(B,\cA,\cS)$.
We restrict ourselves first to the case when $(B-\{R_i\}_{i=1}^N,\cA)$ 
can be immersed in
$(\bR^2,\cA_0)$ where the $R_i$ are a minimal set of disjoint curves, each 
with at an endpoint at the node $s_i$.  
We call the curves $R_i$ {\it branch curves}.
Let $P_i$ be the possible values of $\Phi$ if it were extended, as a
discontinuous map, to $R_i$, i.e. let
$P_i=\{p\in\bR^2 | p=\lim_{x\rightarrow y_i}\Phi(x), y_i\in R_i\}$ and call 
$\cup_{i=1}^N P_i$ the {\it branch locus}.
Take $\Delta=\Phi(B-\{R_i\}_{i=1}^N)\cup \{P_i\}_{i=1}^N$ and
indicate each set $P_i$ by a line with heavy dashes.
Mark the image of a node $s_i$ with the integer $k_i$ 
if the node has multiplicity $k_i\ge 2$ and the branch curve $R_i$ based
at $s_i$ is an eigenray.  (If $R_i$ is an eigenray and the node is
not marked we assume it has multiplicity $1$.)
Now use the conventions of Section~\ref{bdyreduc.sec} to indicate the
stratification.
In the special case that 
a component of the branch locus emanating from a node intersects the 
reduced boundary at a point in the $0$-stratum rather that the $1$-stratum
we would mark intersection point with a heavy dot.

\begin{xca}    Check that the notation for the multiplicity of
a node is sufficient and not redundant.  Specifically, prove 
that if two vectors $u,v\in\bR^2$ satisfy
$A^k_{(a,b)}u=v$ and $u\ne v$, then $u$ and $v$ determine $a,b$ and $k$,
while $u=v$ implies the vector $(a,b)$ is parallel to $u$ and 
$k$ is not determined.
\end{xca}

\begin{example}   
Figures~\ref{nodenbhd.fig}(a) and (b) are both base diagrams 
for the neighborhood of a node.
In the notation of Section~\ref{affinenode.sec} 
they are $\Psi(\widetilde V_{e_1})$ and
$\Psi(\widetilde V_{e_2})$, thereby representing $(V_{e_1},\cA_{e_1})$
and $(V_{e_2},\cA_{e_2})$.
In Figure~\ref{nodenbhd.fig}(b) the map $\Phi$  
extends to a discontinuous map in which the set $P$ (described above) is
the boundary of the missing sector.
In Figure~\ref{nodenbhd.fig}(a) the set $P$ is a ray and the
extension of $\Phi$ is continuous but not differentiable.
\end{example}

\begin{figure}
\begin{center}
	\includegraphics[scale=.5]{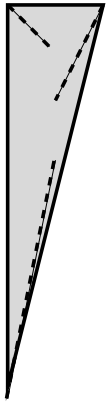}
\end{center}
\caption{Base diagram for an almost-toric fibration of $\bC P^2$.}
\label{funnyCP2.fig}
\end{figure}

\begin{example}   \label{funnyCP2.ex}
Figure~\ref{funnyCP2.fig} shows a base diagram for an 
almost toric fibration with three nodal fibers.
While at first glance one might think that this is a base diagram
for a space with orbifold singularities, the manifold it defines
is in fact smooth.
Sections~\ref{branchmove.sec} and~\ref{nodaltrade.sec} describe the
techniques that allow one to see that the manifold is $\bC P^2$.
\end{example}

To reconstruct the total space from a base diagram we can use
Observation~\ref{basecollapsing.obs} to determine the homology classes
of circles that are collapsed at the reduced boundary and the following 
observation to determine the vanishing class of a nodal fiber.
\begin{obs}   \label{basevanishing.obs}
In a base diagram for an almost toric base, if a primitive integral vector 
$v$ is perpendicular to the eigenline through a node then $v^*=<v,\cdot>$
is a covector that defines the vanishing class of the fiber 
(with respect to the eigenline).
This is an immediate consequence of Proposition~\ref{link.prop}
and the fact that $v^*u=<v,u>$.
\end{obs}

An almost toric fibration can fail to be toric because of the presence
of nodal fibers or because of the lack of an affine immersion of the base
into $\bR^2$. 
To extend to cases where there is no choice of branch curves
$\{R_i\}_{i=1}^N$ such that $B-\{R_i\}_{i=1}^N$ immerses in $\bR^2$ we
follow the usual convention for representing surfaces in $\bR^2$: we 
cut $B$ along more curves and 
use light lines with arrows to indicate components that need to
be identified (with trivial monodromy).
(This same convention could be used to resolve the ambiguities
faced when an immersion into $\bR^2$ obscures the structure of the base
due to multiple covering of portions of the image.)
Proposition~\ref{basepasting.prop} provides criteria for when smooth
components of the boundaries of bases can be identified.

Using these conventions we can reconstruct an almost toric base
from its base diagram.  This is done with the
understanding that smooth components of the branch locus are
identified according to the unique nodal monodromy transformation
(of the correct multiplicity) that maps one onto the other.

\begin{example} \label{S2xT2.ex}  
The manifold $S^2\times T^2$ does not admit a Hamiltonian torus action,
however it is easy to define an almost
toric fibration of it.

Let $S^2\subset \bR^3$ be a sphere of radius $1$ with symplectic structure
$\omega_1$ being the induced area form.
Then the height function gives a Lagrangian fibration 
$\pi_1:(S^2,\omega_1)\rightarrow [-1,1]$ with elliptic singularities at 
the north and south poles.
Equip the torus $T^2$ with angular coordinates $(s,t)$, both of period
$2\pi$, symplectic form $\omega_2=ds\wedge dt$ and Lagrangian fibration
$\pi_2(s,t)=s$.
Then $\pi_1\oplus\pi_2:(S^2\times T^2,\omega_1\oplus\omega_2)
\rightarrow [-1,1]\times S^1$
is an almost toric fibration.
\end{example}

\begin{xca}   \label{draw.ex}
Draw a base diagram for this fibration of $S^2\times T^2$.
\end{xca}

\begin{example} \label{K3.ex1}  
In Section~\ref{K3.sec} we show examples of base diagrams of
almost toric fibrations of the K3 surface.
For now we simply give some background, including a source for
Lagrangian fibrations of the K3.

The smooth K3 surface
is a manifold diffeomorphic to a smooth degree four hypersurface in
$\bC P^3$, i.e. the zero locus of a non-degenerate degree four homogeneous
polynomial in three variables.
While there are many complex K3 surfaces (with different complex structures)
they are all diffeomorphic. 

We start by defining a family of holomorphic fibrations
of K3 surfaces.
Consider a pair of (distinct)
smooth cubic curves in $\bC P^2$ defined by degree three
polynomials $P_1, P_2$; these are   
tori that intersect each other in $9$ points $\{x_1,\ldots, x_9\}$.
The pencil of cubics given by $\{\alpha P_1+\beta P_2 =0\}$ for 
$[\alpha:\beta]\in \bC P^1$ defines a fibration of 
$\bC P^2-\{x_1,\ldots, x_9\}$.
Blowing up those points yields a holomorphic fibration of
$\bC P^2\#9\bC P^1$ over $\bC P^1$;
the generic fiber is a torus (an elliptic curve), 
and for generic choices of $P_1,P_2$ the
singular fibers each have one node on them (and no other singularities).
When equipped with such a holomorphic fibration by elliptic curves, this
manifold is called $E(1)$.
Now take a double branch cover of $E(1)$, branched over a pair of (disjoint)
regular torus fibers.  The result is an elliptically fibered
K3 surface, also known as $E(2)$. 
 
Such a complex K3 surface is equipped with a hyper-K\"ahler structure,
namely a sphere's worth of K\"ahler structures (complex structures
$aI+bJ+ cK$ and Hermitian metrics whose anti-symmetric parts are 
symplectic structures $a\omega_I+b\omega_J+c\omega_K$, 
$a^2+b^2+c^2$).
Given a fibration that is holomorphic with
respect to $J$, it is Lagrangian with respect to $\omega_I$ (or 
$a\omega_I+c\omega_K$, $a^2+c^2=1$).
(In fact it is special Lagrangian, cf.~\cite{GrossWilson.mirrorvia3tori}.)
Generically, a holomorphic fibration has $24$ nodal fibers and hence
so does a Lagrangian fibration.
\end{example}

\subsection{Branch moves} \label{branchmove.sec}

Due to the choices made in the drawing of a base diagram, many base diagrams
correspond to the same almost toric fibration.  The variations resulting
from an application of $AGL(2,\bZ)$ are quite straightforward.  
More involved changes occur when we alter how we cut the base before 
projecting it to $(\bR^2,\cA_0)$.
In particular, a change in a choice of the branch curves $R_i$ constitutes
a {\it branch move}.
For instance, Figures~\ref{nodenbhd.fig} (a) and (b) differ by a branch
move. 

\begin{figure}
\begin{center}
        \psfragscanon
        \psfrag{a}{(a)}
        \psfrag{b}{(b)}
	\includegraphics[scale=.75]{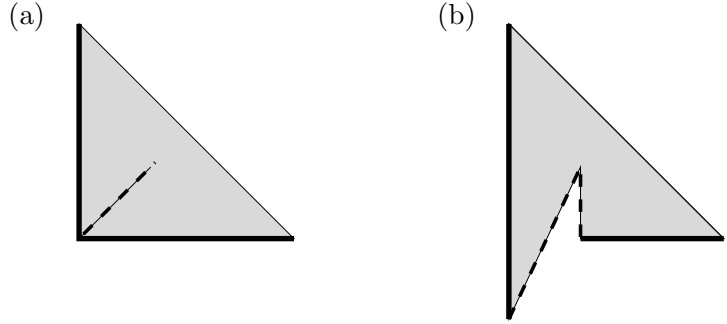}
\end{center}
\caption{Two base diagrams of the same fibration.}
\label{branchmove.fig}
\end{figure}

\begin{xca}   Verify that Figure~\ref{branchmove.fig} 
shows two different
base diagrams for one almost toric fibration and that they 
differ by a branch move.
In both diagrams the   
coordinates of the node are $(1,1)$ and the coordinates of the endpoints
of the line segment with slope $-1$ are $(3,0)$ and $(0,3)$.
\end{xca}

\begin{figure}
\begin{center}
        \psfragscanon
        \psfrag{a}{(a)}
        \psfrag{b}{(b)}
	\includegraphics[scale=.5]{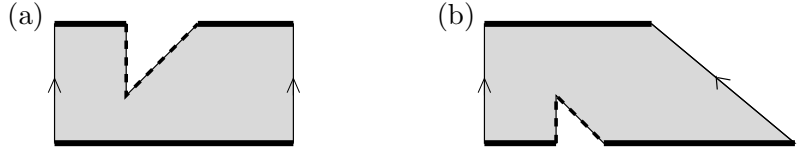}
\end{center}
\caption{Base diagrams for one fibration of 
$(S^2\times T^2) \#\overline{\bC P}^2=
(S^2\tilde\times T^2) \#\overline{\bC P}^2$.} 
\label{diffequiv.fig}
\end{figure}

\begin{example} \label{S2xT2blowup.ex}  
Figure~\ref{diffequiv.fig}(a) shows a base diagram for 
$(S^2\times T^2)\#\overline{\bC P}^2$.  
This follows from Exercise~\ref{draw.ex} 
and the almost toric blow-up described
in Section~\ref{blowup.sec}. 
(Filling in the missing wedge gives the answer to Exercise~\ref{draw.ex}.)
Notice that the eigenline of the node is parallel to the reduced boundary.
\end{example}

\begin{xca} \label{twistedS2xT2blowup.ex}  
Verify that Figure~\ref{diffequiv.fig}(b) is a base diagram for
a fibration of 
$(S^2\tilde\times T^2) \#\overline{\bC P}^2$.
Hint: 
Fill in the missing wedge to blow down 
(as per Section~\ref{blowup.sec}) and then convince yourself that
this is the base diagram for a sphere bundle over a torus.
Use the techniques of Section~\ref{selfint.sec} to calculate
that the self-intersections of the tori represented by the 
connected components of the reduced boundary are $\pm 1$.
\end{xca}

As an application of branch moves we can use base diagrams to see that
\begin{prop}
The manifolds $(S^2\times T^2)\#\overline{\bC P}^2$ and 
$(S^2\tilde\times T^2) \#\overline{\bC P}^2$ are diffeomorphic.
\end{prop}
\begin{proof}
As already noted, Figures~\ref{diffequiv.fig}(a)
and (b)  are base diagrams for  $(T^2\times S^2)\#\overline{\bC P}^2$ 
and $(S^2\tilde\times T^2) \#\overline{\bC P}^2$ respectively.
In  Figure~\ref{diffequiv.fig}(a) the monodromy across the vertical
upward ray is  $A_{(1,0)}$ and across the downward vertical ray (which
is not drawn) is the identity.
To get from the first to the second by a branch move cut
Figure~\ref{diffequiv.fig}(a) along the downward vertical and apply
the linear map $A^{-1}_{(1,0)}$ to the right half
of the image (with the origin at the node); this makes the
monodromy across the upward vertical ray the identity and the monodromy
across the downward vertical ray $A^{-1}_{(1,0)}$.
Reglue along the upward vertical ray, 
thereby obtaining Figure~\ref{diffequiv.fig}(b).
\end{proof}

\begin{xca}   \label{diffequiv.ex}
Adapt this argument to show that 
$(S^2\times \Sigma)\#\overline{\bC P}^2$ and 
$(S^2\tilde\times \Sigma) \#\overline{\bC P}^2$ are diffeomorphic for
any surface $\Sigma$.  (The case of $\Sigma=S^2$ is discussed in the
appendix.)
Hint: Blow up inside a toric fibration of $S^2\times D^2$ where $D^2$ is
a disk in $\Sigma$.  Choose the symplectic structure and the disk
so that the area
of the fiber $S^2$ is smaller than the area of the $D^2$.
\end{xca}

\subsection{Blowing up and down} \label{blowup.sec}

Blowing up is a simple surgery which, on a four-manifold,
amounts to replacing a ball with a neighborhood of a sphere of
self-intersection $-1$.
It is well known how, in the toric category, to blow up at a point that
is itself a fiber (i.e. a rank zero elliptic singular point).
After reviewing this we explain how, in the almost toric category, to
blow up at a point that lies on a one-dimensional fiber at the expense of
introducing a nodal fiber.
Blowing down is just the reverse operation to blowing up.  

To blow up in the symplectic category one proceeds as follows:
find a symplectic embedding $\varphi:(B^4(r),\omega_0)\rightarrow (M,\omega)$ 
of the  closed
ball of radius $r$ with the standard symplectic structure, consider
$M-{\rm Int}(\varphi(B^4(r)))$ and perform a boundary reduction along 
$\varphi(\bdy B^4(r))$, thereby collapsing the fibers of
the Hopf fibration.  

\begin{figure}
\begin{center}
        \psfragscanon
        \psfrag{a}{(a)}
        \psfrag{b}{(b)}
	\includegraphics[scale=.75]{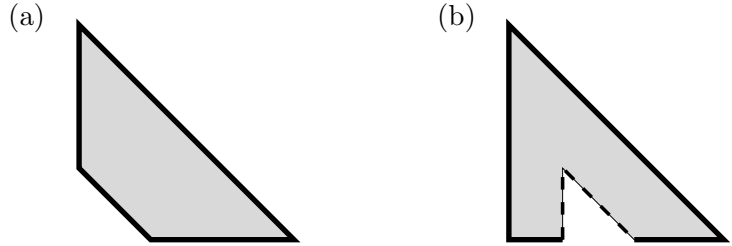}
\end{center}
\caption{Two fibrations of $\bC P^2 \#\overline{\bC P}^2$}
\label{blowup.fig}
\end{figure}

In general the symplectic embedding $\varphi$ could be very complicated.
However, if $(M,\omega)$ is equipped with a Lagrangian fibration that has 
a rank zero elliptic singularity, then for sufficiently small $r$ 
one can choose  the embedding $\varphi$ so that
the origin gets mapped to the rank zero singular point and it is
fiber preserving with respect to the standard toric fibration 
of $B^4(r)\subset(\bR^4,\omega_0)$ (as in Example~\ref{B4.ex}).
Then the surgery can be described in terms of a surgery on the base:
remove from the base a triangular neighborhood of the vertex that is
the image of the rank zero singular point; do this in such a way that
the edges of the triangle meeting at that vertex have length
$\frac{1}{2}r^2$ times the length of a primitive integral tangent vector; 
then connect the two resulting vertices of the
base with a heavy line.  (See Figure~\ref{blowup.fig}(a).)
The resulting base reveals the change in topology: the new edge
represents 
a symplectic $2$-sphere of area $\pi r^2$ whose self-intersection
is $-1$. (See Sections~\ref{area.sec} and~\ref{selfint.sec}.)

With the allowance of nodal fibers there is another \lq\lq fibration
compatible\rq\rq\  way to blow up, namely an {\it almost toric blow-up}.  
We describe the surgery on the
base in terms of the base diagram 
and then explain how it corresponds to blowing up a point.
(This operation was initially observed by Zung~\cite{Zung.II}.)

Consider an almost toric 
fibration $\pi:(M,\omega)\rightarrow (B,\cA,\cS)$ such that 
the $1$-stratum of $B$ is non-empty.
(Recall that the $1$-stratum must belong to $\bdy_R B$.)
Choose a base diagram for $B$ such that for some $a,\varepsilon>0$
the set
$\{(p_1,p_2)| p_1>-\varepsilon , p_2\ge 0, p_1+p_2<a+\varepsilon\}$
with the boundary components marked by heavy lines represents a fibered
neighborhood of a point in the $1$-stratum.
Remove the triangle with vertices $(0,0)$, $(a,0)$, $(0,a)$ and insert
line segments with heavy dashes so that their endpoints are
$(0,a),(0,0)$ and $(0,a),(a,0)$.
Viewing $(0,a)$ as the origin, $A_{(1,0)}$ maps the vertical dashed segment
to the diagonal one and preserves horizontal vectors.
Therefore, this is a base diagram for an almost toric fibration and with
respect to the induced standard coordinates the nodal fiber has vanishing
class $[\frac{\del}{\del q_2}]$.
Notice that
this is the same as the collapsing class of the component of the 
$1$-stratum in question.
See Figure~\ref{blowup.fig} for an example.

Theorem~\ref{almosttoricbase.cor} guarantees that 
the described operation on the base gives a well-defined surgery
on the manifold $(M,\omega)$.
To see that this surgery amounts to blowing up a point notice that
we start with a neighborhood diffeomorphic to $S^1\times D^3$; on a
topological level, the effect of the surgery is to add a vanishing
disk (to create the nodal fiber) whose boundary represents the same
class as the collapsing class of the almost toric fibered $S^1\times D^3$,
i.e. to adding a $-1$-framed $2$-handle along a contractible circle.
But this is the same thing, topologically, as making a connect sum
with $\overline{\bC P}^2$.

\begin{xca}   
Use the techniques of Section~\ref{readingI.sec} 
to verify that the preimage of either segment marked by heavy dashes
is a visible sphere of self-intersection $-1$.
(Both segments yield the same visible sphere.)
This is the exceptional sphere introduced by the almost toric blow-up.
\end{xca}

\begin{rmk}  
In the toric case -- and in the complex category -- it is clear what
point is being blown up; this is not so in either a symplectic blow-up
or an almost toric blow-up.
\end{rmk}

\subsection{Closed almost toric manifolds} \label{closedalmosttoric.sec}

As mentioned in the introduction, the set of closed manifolds that
admit almost toric fibrations is a modest expansion over those that
admit toric fibrations.  More precisely, we can list them:

\begin{thm}[Leung and Symington~\cite{LeungSymington.almosttoric}]
\label{closedalmosttoric.thm}
If $M$ is a closed almost toric manifold then $M$ is diffeomorphic
to $S^2\times S^2$, $S^2\times T^2$, $N\#n\overline{\bC P}^2$ with
$N=\bC P^2$ or $S^2\tilde\times T^2$ and $n\ge 0$, 
a torus bundle over a torus with monodromy
\begin{equation}
\left\{\left(  
\begin{matrix}
1 & \lambda \\  
0 & 1 
\end{matrix} 
\right), Id\right\} \quad {\rm with\ } \lambda\in\bZ,
\end{equation} 
a torus bundle over the Klein bottle that is double covered by one of
the above bundles over a torus, 
the K3 surface or the Enriques surface.
\end{thm}
In this section we give a brief sketch of the proof,
referring the reader to~\cite{LeungSymington.almosttoric} for details.

There are two parts to the proof of Theorem~\ref{closedalmosttoric.thm}:
showing that the aforementioned manifolds admit almost toric fibrations
and showing that there are no other almost toric manifolds.

As to the first part, notice that $S^2\times S^2$ and 
$\bC P^2\#n\overline{\bC P}^2$ admit toric fibrations and hence
almost toric fibrations. 
Example~\ref{S2xT2.ex} gives an almost toric fibration of $S^2\times T^2$.
A base defining an almost toric fibration of $S^2\tilde\times T^2$
was constructed in Exercise~\ref{twistedS2xT2blowup.ex}.
Applying the almost toric blow-up of Section~\ref{blowup.sec} 
to these examples yields almost toric fibrations
of $S^2\times T^2\#n\overline{\bC P}^2$ and 
$S^2\tilde\times T^2\#n\overline{\bC P}^2$ (with the case $n=1$ shown in
Figures~\ref{diffequiv.fig}(a) and (b)).
Note that we do not need to list separately the blow-ups of these
trivial and non-trivial bundles as they are diffeomorphic.
(See Section~\ref{diffequiv.sec}.)
Any fibration over the Mo\"ebius band is a $\bZ_2$ quotient of a 
fibration over the cylinder and hence can be constructed explicitly
in terms of a $\bZ_2$ quotient of the cylinder with an integral affine
structure.

Meanwhile, Geiges' classification of torus bundles over a torus that admit 
a symplectic structure~\cite{Geiges.torusbundles} implies
that a such a bundle can have Lagrangian fibers
if and only if the manifold is
diffeomorphic to a bundle with the monodromy given in the statement of
Theorem~\ref{closedalmosttoric.thm} .
Taking $\bZ_2$ quotients of the base, when possible, yields the 
total spaces of all the regular
torus bundles over the Klein bottle that admit Lagrangian fibers.

Explicit almost toric fibrations of the K3 surface are given in 
Section~\ref{K3.sec}.
The first of these, built by taking symplectic sums of eight copies of the
fibration defined by Figure~\ref{octant.fig}, has lots of symmetry.
Taking a fiber-preserving $\bZ_2$ quotient of this fibration yields 
a fibration with $12$ singular fibers and base $\bR P^2$.  
One way to see that this is the Enriques surface is to notice that the
resulting manifold is diffeomorphic to the one obtained by symplectic
summing four copies of the fibration defined by
Figure~\ref{octant.fig} to give a fibration of 
$E(1)=\bC P^2\#9\overline{\bC P}^2$ and then performing two smooth 
multiplicity two 
log transforms along tori isotopic to the preimage of the reduced boundary.
Since this latter construction is known to yield the Enriques surface
(cf. \cite{GompfStipsicz.4mflds}) we are done.

The proof  that the set of manifolds 
mentioned in Theorem~\ref{closedalmosttoric.thm} 
contains all four-manifolds that admit an almost toric fibration 
is outlined by the following lemmas which are proved 
in~\cite{LeungSymington.almosttoric}.

\begin{lemma} \label{almosttoricbase.lem}
If $(B,\cA,\cS)$ is a compact almost toric base such that $\bdy_R B=\bdy B$,
then $B$ must have non-negative
Euler characteristic and hence must be homeomorphic to a disk,
a cylinder, a Mo\"ebius band, a sphere, $\bR P^2$, a torus or a Klein bottle.
Furthermore, if $B$ is a closed manifold then it has $12\chi$ nodes
(counted with multiplicity) where $\chi$ is the Euler characteristic
of $B$.
\end{lemma}
(This lemma, in a more general form, is stated but not proven 
in~\cite{Zung.II}.)  
The proof is an application of the Gauss-Bonnet
theorem.
We separate the arguments for when $B$ has boundary and when it does not.
Notice that the hypotheses that $B$ is compact and $\bdy_R B=\bdy B$ imply
that $B$ is the base of a fibration of a closed manifold.

If $B$ is a closed manifold then there must be $12k$ nodes for some
integer $k\ge 0$ since any product of matrices conjugate to $A_{(1,0)}$
that equals the identity is a product of $12k$ such matrices 
(cf.~\cite{GompfStipsicz.4mflds} or~\cite{FriedmanMorgan.cmplxsurf}).
Choosing $k$ disjoint open (topological) disks that each contain $12$ nodes
around which the monodromy is trivial, put a flat metric on the
complement of the union of the disks and then extend the metric over the
rest of the base $B$.  The fact that the sphere has an affine structure
with $24$ nodes then implies that each collection of $12$ nodes (with
trivial monodromy) contributes $2\pi$ to the total curvature, or $1$ to
the Euler characteristic.  Therefore a closed base can only be a
sphere with $24$ nodes, $\bR P^2$ with $12$ nodes, or $T^2$ or the Klein
bottle with no nodes.

If $B$ has boundary then
heuristically, $B$ has positive Euler characteristic because its boundary
is locally convex and each node contributes non-negative curvature.
To calculate the Euler characteristic we can appeal to
the Gauss-Bonnet formula for polygons -- that the
sum of the total curvature of the polygon, the geodesic curvature of
its boundary and the discontinuities at the vertices equals $2\pi$
(cf.~\cite{Spivak.diffgeomIII}).
Triangulating $B$ and summing over the triangles we 
find the Euler characteristic as $\frac{1}{2\pi}$
times the sum of the contributions from the triangles.
Since an affine manifold with nodes has no metric that is compatible with the
affine structure we work instead with a compact manifold
$(B',\cA')$ that is homeomorphic to $(B,\cA)$, has no nodes, and 
such that $(B'-P',\cA')$ is isomorphic to $(B-P,\cA)$ where $P'$ is a subset
of $\bdy B'$ and $P$ is a union of disjoint curves, each with one endpoint
at a node in $(B,\cA)$.  Then the affine structure $\cA'$ induces a flat
metric on $B'$, allowing us to carry through the Gauss-Bonnet argument.

Now we can consider what manifolds can have an almost toric fibration
over each of the bases listed in Lemma~\ref{almosttoricbase.lem}.  
\begin{lemma} \label{diskbase.lem}
If there is an almost toric fibration $\pi:(M,\omega)\rightarrow (B,\cA,\cS)$
where $B$ is homeomorphic to a disk and $M$ is a closed manifold, 
then there is a toric fibration $\pi':(M,\omega')\rightarrow (B,\cA')$.
\end{lemma}
The proof of this lemma relies on branch moves to change the fibration
to a new one in which the base is equivalent, after making nodal trades, to
a toric base.  We do not include the details here but merely note that
to make the branch moves it may be necessary
to deform the symplectic structure (as in Section~\ref{deformations.sec}),
hence the indicated change in the symplectic structure from $\omega$
to $\omega'$.

\begin{lemma} \label{cylinder.lem}
Let $(B,\cA,\cS)$ be an almost toric base such that $\bdy_R B=\bdy B$.
If $B$ is a closed cylinder then two boundary components are parallel.
If  $B$ is either a cylinder or a Mo\"ebius band, then the eigenlines of any
nodes must be parallel to $\bdy B$.
\end{lemma}
The essence of this lemma can be seen in Figure~\ref{diffequiv.fig}(a)
where, if 
the reduced boundary components were not parallel, there would be no
integral affine transformation that could glue the collar neighborhoods of the
non-reduced boundary components to obtain a cylinder.  
Furthermore, the introduction of a node whose eigenline is not parallel
to the boundary would have the effect of making the boundary components
converge.  
With these restrictions, the only closed manifolds that have almost toric
fibrations over a cylinder or a Mo\"ebius band are sphere bundles over
tori and their blow-ups, and we have constructed fibrations of all of
them. 

We already saw above that Geiges' classification implies that any almost
toric fibration over a torus or a Klein bottle is diffeomorphic to one
of the manifolds listed in the theorem.

If the base is a sphere then it must have $24$ nodes.  Then it is the base
of an almost toric fibration with $24$ nodal fibers.  With respect to a 
different symplectic structure, the manifold is the total space of a
Lefschetz fibration with $24$ singular fibers and hence is the
K3 surface (cf.~\cite{GompfStipsicz.4mflds}).
Finally, an almost toric fibration over $\bR P^2$ must be double covered
by an almost toric fibration over $S^2$, namely the K3 surface.
Therefore, up to diffeomorphism we have only one such manifold.

\section{Modifying fibrations}

\subsection{Nodal slides and nodal trades} \label{nodaltrade.sec}

The first operation we introduce in this section is a nodal slide:

\begin{defn}   \label{nodalslide.def}
Two almost toric bases $(B,\cA_i,\cS_i)$, $i=1,2$, are related by
a {\it nodal slide} if there is a curve $\gamma\subset B$
such that $(B-\gamma,\cA_1,\cS_1)$ and $(B-\gamma,\cA_2,\cS_2)$ 
are isomorphic and $\gamma$ contains one node and belongs
to the eigenline $L_i$ through that node.
\end{defn}

Any two bases that are related by a nodal slide belong to a family
of bases in which the node \lq\lq slides\rq\rq\ along the eigenline.
(The affine structures on the complement of the eigenlines are not
affected by the position of the node.)
These define a family of almost toric fibrations of one 
manifold and hence a family of symplectic structures on that manifold.
Furthermore, the manifolds are (fiberwise)
symplectomorphic on the complement of
a compact set belonging to the preimage of the eigenline. 
Therefore we can work in a contractible neighborhood of that preimage where
the symplectic structures are exact (since the symplectic structure
in a neighborhood of a regular Lagrangian fiber is exact).
Invoking Moser's argument~\cite{McDuffSalamon.intro} we can conclude that the
there is an isotopy, compactly supported in that neighborhood, that defines
a symplectomorphism between the almost toric manifolds fibering over
the two bases.
Therefore,

\begin{prop} \label{isotopic.prop} 
If two bases are related by a nodal slide then they define the
same manifold equipped with isotopic symplectic structures.
\end{prop}

Under certain circumstances, sliding a node right into an edge in $\bdy_R B$
and replacing the limit point in the edge with a vertex yields an
almost toric base.
If so, then the new base again determines a manifold symplectomorphic
to the original one.
For instance, Figure~\ref{B4CP2.fig}(a) and
Figure~\ref{branchmove.fig}(a) are different fibrations of the same
symplectic four-ball.
Specifically,

\begin{lemma} \label{nodaltrade.lem}
Given an almost toric base $(B,\cA,\cS)$ let $R$ be an embedded eigenray 
that connects a node $s$ with a point $b$ in an edge $E\subset\bdy_R B$
such that there are no other nodes on $R$.
Let $v^*,w^*\in T_b^* B$ be the vanishing and collapsing covectors.
If $v^*$ and $w^*$ span $\Lambda_b^*$ then there is an almost toric
base $(B,\cA',\cS')$ such that
\begin{enumerate}
\item $(B,\cA',\cS')$ contains one fewer node than $(B,\cA,\cS)$,
\item $(B-R,\cA',\cS')$ is isomorphic to $(B-R,\cA,\cS)$ and
\item in $(B,\cA',\cS')$ the intersection of $R$ and 
$\bdy_R B$ is a vertex.
\end{enumerate}
\end{lemma}

\begin{proof}
Without loss of generality we can work locally in a neighborhood of
$R$.
The condition that $v^*$ and $w^*$ span $\Lambda_b^*B$ implies that
we can choose a base diagram in which $v^*=(-1,1)$ and $w^*=(0,1)$.
Then, possibly after rescaling, a neighborhood of the dashed line in 
Figure~\ref{branchmove.fig}(a) provides a base diagram for 
$(B,\cA,\cS)$ with the dashed line representing $R$.
Simply erasing the dashed line reveals a neighborhood in an appropriately
scaled copy of Figure~\ref{B4CP2.fig}(a) that serves as $(B,\cA',\cS')$,
clearly satisfying all of the conditions listed in the lemma.
This suffices to prove the lemma.
Note that the condition that $v^*$ and $w^*$ span $\Lambda_b^*$ is
required to prevent the vertex from being the image of an orbifold
singular point.
\end{proof}

\begin{rmk}  
If the edge in $B$ and an eigenray emanating from the 
node do not intersect inside
$B$, the symplectic manifold defined by $B$ is still a ball, but the ball
does not have a toric fibration that agrees with the almost toric fibration
near the boundary.
\end{rmk}

\begin{def}   \label{nodaltrade.def}
Two almost toric bases $(B,\cA,\cS)$ and $(B,\cA',\cS')$ 
are related by a {\it nodal trade} if they satisfy the hypotheses of
Lemma~\ref{nodaltrade.lem}.
\end{def}

\begin{thm} \label{nodaltrade.thm}
Two almost toric bases that are related by a nodal trade are symplectomorphic;
in fact their symplectic structures are isotopic.
\end{thm}

\begin{proof}
It suffices to prove that the almost toric bases
$(U,\cA,\cS)$ and $(U,\cA',\cS')$ in the proof of Lemma~\ref{nodaltrade.lem}
define diffeomorphic four-manifolds with symplectic forms that are
isotopic via an isotopy that is compactly supported in a smaller neighborhood.
We already know that $(U,\cA,\cS)$ is the base of a fibration of an
open four-ball.
To see that $(U,\cA',\cS')$ also defines a four-ball, notice that it can
be viewed as $D^4$ with a $1$-handle attached (to get $S^1\times D^3$)
and then a $-1$-framed $2$-handle (a thickening of the vanishing disk)
attached along a circle that runs once over the $1$-handle
(since $v^*,w^*$ span $\Lambda_b^*$).
Such a $2$-handle cancels the $1$-handle so the manifold is $D^4$. 

Now repeat the Moser argument proving Proposition~\ref{isotopic.prop} 
to complete the proof of the theorem.
\end{proof}

\begin{rmk}  
As observed by Zung~\cite{Zung.II}, the one parameter
family of fibrations corresponding to sliding a node into an edge 
and achieving a nodal trade is, in the language of integrable systems,
a Hamiltonian Hopf bifurcation that occurs in many
Hamiltonian systems~\cite{vanderMeer.HamiltHopf}.
\end{rmk}
\begin{xca}   \label{multiplicityknodes.ex}
Suppose $(B,\cA,\cS)$
has a multiplicity $k$ node at $s$ and $(B,\cA',\cS')$ has
$k$ distinct multiplicity $1$ nodes on a segment $L$ of the eigenline
through these nodes.
Use the same type of arguments to prove that if 
$(B-L,\cA,\cS)$ and $(B-L,\cA',\cS')$ are isomorphic, then 
$(B,\cA,\cS)$ and $(B,\cA',\cS')$ define the same four-manifold with
isotopic symplectic structures.
\end{xca}

\subsection{Connecting toric fibrations} \label{diffequiv.sec}

It is well known that $(S^2\times S^2,\omega\oplus\lambda\omega)$
can admit distinct toric fibrations and that the number of such
fibrations depends on $\lambda$.
If, without loss of generality, we assume $\lambda\ge 1$, then
the number of such fibrations is $\lceil \lambda \rceil$, the
smallest integer greater than or equal to $\lambda$.
The number of these fibrations is interesting in that it equals
the number of non-conjugate tori in the group of Hamiltonian
symplectomorphisms of 
$(S^2\times S^2,\omega\oplus\lambda\omega)$~\cite{Karshon.maximaltori}.
Furthermore, each time $\lambda$ passes an integer the topology of
the full symplectomorphism group changes~\cite{Gromov.pseudohol,
Abreu.symplectom, AbreuMcDuff.symplectom, Anjos.symplectom}.

Using nodal slides and branch moves the different toric fibrations
of $(S^2\times S^2,\omega\oplus\lambda\omega)$ can be connected
by a path of almost toric fibrations.
This provides a simple geometric proof of the 
the diffeomorphism equivalence of the manifolds defined by the 
$\lceil \lambda \rceil$ different toric bases and also, invoking
a Moser argument, the isotopy equivalence of the induced symplectic
forms. 

\begin{figure}
\begin{center}
        \psfragscanon
        \psfrag{a}{(a)}
        \psfrag{b}{(b)}
	\includegraphics[scale=.75]{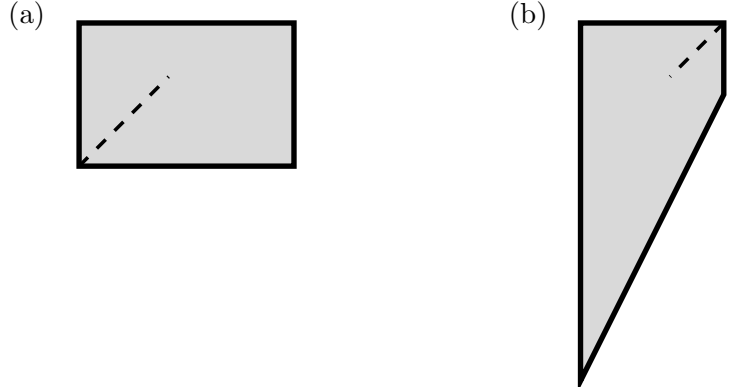}
\end{center}
\caption{Branch move for one fibration of $(S^2\times S^2,\sigma\oplus
\frac{3}{2}\sigma)$}
\label{Hirzebruch.fig}
\end{figure}

We illustrate this in Figure~\ref{Hirzebruch.fig} in which
case $\lambda=\frac{3}{2}$.
Figures~\ref{Hirzebruch.fig}(a) and (b) show fibrations that clearly
differ from toric ones by a nodal trade
and the two figures differ from each other by a branch move.

\begin{conj} \label{connecting.conj}
Any two toric fibrations of $(M,\omega)$ can be connected by a path
of almost toric fibrations of $(M,\omega)$.  (Note that $\omega$ is
fixed throughout this path.)
\end{conj}

\subsection{Deformations} \label{deformations.sec}

A {\it deformation} of a symplectic manifold is a one-parameter
family of symplectic structures.  In general the cohomology
class of the symplectic structure varies during the deformation.
Any one-parameter family of almost toric bases that defines the same
smooth manifold defines a deformation of the symplectic structure.

One obvious deformation arises from rescaling an almost toric base,
thereby causing a rescaling of the symplectic structure.
Also, if the base is non-compact or has boundary then removing portions of the
top stratum of the base or adding to it so as
to yield a new base homeomorphic to the first will give a deformation.

Two other useful variations of the base are \lq\lq thickening\rq\rq\
and \lq\lq thinning\rq\rq\ along the reduced boundary.  
(These deformations were originally introduced 
in~\cite{McDuffSymington.assoc} where they were used to prove
the diffeomorphism equivalence of certain four-manifolds.)
Specifically, an edge in the reduced boundary
has a collar neighborhood 
fibered by curves that are affine linear.
When the edge contains two vertices the collar neighborhood
has a base diagram as shown in Figure~\ref{spherenbhd.fig} for
some choice of length $a$ of the heavy horizontal line and slope $\frac{1}{n}$,
$n\in\bZ$, for the rightmost heavy line.
The thickness of the neighborhood is the length of the vertical heavy line,
say $t$.  Call such a neighborhood $(N_{a,n,t},\cA_0,\cS)$.
Notice that for fixed lengths $a,t$ and fixed $n$ the base
$(N_{a,n,t},\cA_0,\cS)$ belongs to a family 
$(N_{a-n\tau,n,t+\tau},\cA_0,\cS)$ 
where $\tau\in(-t,\frac{a}{n})$ if $n> 0$ and $\tau\in(-t,\infty)$ if 
$n\le 0$.
(The upper bound on $\tau$ when $n>0$ merely guarantees that $a-n\tau$
stays positive.)
Then to {\it thicken} or {\it thin} $(B,\cA,\cS)$ 
along an edge containing two vertices means to
replace a collar neighborhood of it that is isomorphic to
$(N_{a,n,t},\cA_0,\cS)$ with one isomorphic to 
$(N_{a-n\tau,n,t+\tau},\cA_0,\cS)$ for some $\tau>0$ or $\tau<0$ respectively.

\begin{xca}    Extend this definition of thickening and thinning to
the case when the $1$-stratum is compact.  Hint:
Examples~\ref{S2xT2.ex} and \ref{funnyCP2.ex} contain base diagrams
in which the $1$-stratum is compact.
\end{xca}

The effect on the cohomology class of the symplectic structure is evident
from the changes in the areas of the spheres (or torus) that are the
preimage of the closed edges of $B$.
This is discussed again near the end of Section~\ref{selfint.sec}.

\section{Reading the base I: surfaces } \label{readingI.sec}

\subsection{Visible surfaces} \label{visible.sec}

Some surfaces embedded in a toric fibered manifold are easily reconstructed
from their images in the base, such as the symplectic spheres
that are preimages of edges in the base of a toric fibration of
a closed manifold. 
In this section we explain how a curve in the base of an almost toric 
fibration and the choice of
elements of $H_1$ of regular fibers of the curve together
determine surfaces in the total space.
We also show how this data determines whether or not the surface is
symplectic or Lagrangian.

Throughout this section, given an immersed curve 
$\gamma:I\rightarrow (B,\cA,\cS)$, we use $\{\gamma_i\}_{i=1}^k$ to denote 
the continuous components of 
$\gamma|_{\gamma^{-1}(B-\bdy_R B)}$.
We assume that if $\gamma$ approaches and then coincides with a
component of the $1$-stratum of $\bdy_R B$ (on a set with non-empty 
interior) that the tangency of
$\gamma$ with $\bdy_R B$ at a point of first contact is sufficiently
high so as to ensure that $\Sigma_\gamma $ is non-singular.  
Specifically, if there are local standard coordinates in which $\gamma$ is
given by $(p_1,f(p_1))$  then at such a point of first tangency
$\sqrt{f}$ is differentiable.

\begin{defn}   \label{affinecircle.def}
A circle in a regular fiber $F_b$ of a Lagrangian fibration is an 
{\it affine circle} if it has a tangent vector field belonging to
$\Lambda^{\rm vert}$.
\end{defn} 

\begin{defn}   \label{visiblesurface.def}
A {\it visible surface} $\Sigma_\gamma$ in an almost toric fibered manifold
$\pi:(M,\omega)\rightarrow (B,\cA,\cS)$
is a an immersed surface whose image is an immersed curve
$\gamma$ with transverse self-intersections such that 
$\pi|_{\Sigma_\gamma\cap \pi^{-1}(B_0)}$ is a submersion onto
$\gamma\cap B_0$,  
any non-empty intersection of $\Sigma_\gamma$ with a regular fiber 
is a union of affine circles, and no component of $\bdy\Sigma_\gamma$ projects
to a node.
\end{defn}
Notice that a visible surface must be a sphere, a disk, a cylinder or torus.

\begin{defn}   \label{compatible.def}
Given an immersed curve $\gamma:I\rightarrow (B,\cA,\cS)$,
let $\{\gamma_i\}_{i=1}^k$ be the continuous components of 
$\gamma|_{\gamma^{-1}(B-\bdy_R B)}$.
A primitive class ${\mathbf a}_i$ in $H_1(\pi^{-1}(\gamma_i),\bZ)$ 
(such that $\pi_* {\mathbf a}_i=0$ if $\gamma_i$ is a loop) is
{\it compatible} with $\gamma$ if all of the following are satisfied:
\begin{enumerate}
\item ${\mathbf a}_i$ is the 
vanishing class of every node in $\gamma_i$,
\item 
$|{\mathbf a}_i\cdot{\mathbf c}|\in\{0,1\}$ for each ${\mathbf c}$ that 
is the collapsing class, with respect to $\overline \gamma_i$, 
for a component of
the $1$-stratum of $\bdy_R B$ that intersects $\overline \gamma_i$,
\item
$|{\mathbf a}_i\cdot{\mathbf c}|=1$ if $\overline \gamma_i$ 
intersects the $1$-stratum non-transversely,
\item
$|{\mathbf a}_i\cdot{\mathbf d}|=1$ for each ${\mathbf d}$ that is one of the 
two collapsing classes at a vertex contained in the closure of
$\gamma_i$.
(Here $\cdot$ is the intersection pairing in $H_1(\pi^{-1}(\gamma_i),\bZ)$
and $\overline\gamma_i$ is the closure of $\gamma_i$.)
\end{enumerate}
\end{defn}
Note that we defined the collapsing class to be in the first homology of
a regular fiber.  But since the collapsing class in $H_1(F_b,\bZ)$,
$F_b$ being any fiber in $\pi^{-1}(\gamma_i)$, has the same image
in $H_1(\pi^{-1}(\gamma_i),\bZ)$ induced by the inclusion map, it
makes sense to view the collapsing class as living in
$H_1(\pi^{-1}(\gamma_i),\bZ)$.

\begin{thm} \label{visiblesurfaces.thm}
Suppose $(B,\cA,\cS)$ is an almost toric base such that each node has
multiplicity one.
An immersed curve $\gamma:I\rightarrow (B,\cA,\cS)$ with transverse 
self-intersections and a set of compatible classes $\{{\mathbf a}_i\}_{i=1}^k$
together determine a visible surface $\Sigma_\gamma$ such that for each
$b\in\gamma_i$, 
\begin{equation} \label{compat.eq}
\iota_*[\Sigma_\gamma\cap F_b]={\mathbf a}_i
\end{equation}
where $\iota: F_b\rightarrow \pi^{-1}(\gamma_i)$ is the inclusion map.
The surface $\Sigma_\gamma$ is unique up to 
isotopy among visible surfaces in the preimage of $\gamma$ that satisfy
Equation~\ref{compat.eq}.
Furthermore, no such surface exists if the classes ${\mathbf a}_i$
are not compatible with $\gamma$.
\end{thm}

\begin{proof}

If $\gamma\subset\bdy_R B$ then $\gamma$ defines a symplectic submanifold
with an induced almost toric fibration over $\gamma$.
In this case $\Sigma_\gamma$ is a sphere or disk if 
$\gamma$ has two endpoints or one endpoint, respectively, 
in the $0$-stratum; it is a torus if $\gamma$ is compact and without
boundary; otherwise $\Sigma_\gamma$ is a cylinder. 
Tori projecting to $\gamma\subset\bdy_R B$ 
can be found in Examples~\ref{funnyCP2.ex} 
and~\ref{S2xT2blowup.ex},
among others.

Now suppose $\gamma\cap(B-\bdy_R B)\ne \emptyset$.  In that case we
certainly need to specify classes ${\mathbf a}_i$ in order to obtain
a visible surface $\Sigma_\gamma$ that is well defined up to isotopy.
Since the compatibility conditions for the ${\mathbf a}_i$ deal only with
local phenomena (near a given fiber), we assume for
simplicity of exposition that $\gamma$ is embedded.

First suppose that $\gamma_i$ contains no nodes.
Then $\pi^{-1}(\gamma_i)$ is a trivial torus bundle over $\gamma_i$.
Choose a trivialization $T^2\times\gamma_i$ of this torus bundle.
Then for any ${\mathbf a}_i$ and any section 
$\sigma:\gamma_i\rightarrow \gamma_i\times T^2$ there is a unique  visible surface
that is a cylinder which contains the image of the section and
whose intersection with any fiber represents ${\mathbf a}_i$.
Since all sections are isotopic, $\Sigma_{\gamma_i}$ is unique up to
isotopy among visible surfaces over $\gamma_i$.

If $\gamma_i$ (and hence $\gamma$) has a node $s$ as an endpoint, then
the requirement that no component of $\bdy\Sigma_\gamma$ project to a
node implies $\pi^{-1}(s)\cap\Sigma_\gamma$ is the nodal singular point. 
Furthermore, the smoothness of $\Sigma_\gamma$ implies that
${\mathbf a}_i$ is the vanishing class of that node.
If a node $s$ is on the interior of $\gamma_i$, then because the vanishing
class is the only well defined class in a neighborhood of the node,
in order for
the class of the affine circles in $\Sigma_\gamma$ to be well defined as
$\gamma$ passed through a node we must have the ${\mathbf a}_i$ is
the vanishing class for that node. 

Next suppose $\overline \gamma_i$ intersects at $b$ a component of
the $1$-stratum of $\bdy_R B$ that has collapsing class ${\mathbf c}$
with respect to $\overline \gamma_i$.  
If the preimage of $b$ in $\Sigma_\gamma$ is a point (in which case
$b$ must be an endpoint of $\gamma$), then the smoothness of 
$\Sigma_\gamma$ implies that ${\mathbf a}_i$ must be, up to sign,  
the collapsing
class ${\mathbf c}$, i.e. ${\mathbf a}_i\cdot{\mathbf c}=0$.
Furthermore, as can be checked in local coordinates,
transversality of the intersection of $\gamma$ and the reduced boundary
at $b$ guarantees smoothness of $\Sigma_\gamma$.
If however the preimage of $b$ in $\Sigma_\gamma$ 
is a circle $C$ that represents the generator of
$H_1(F_b,\bZ)=\bZ$, 
for $C$ to be a limit of representatives of ${\mathbf a}_i$
we must have ${\mathbf a}_i={\mathbf d}+k{\mathbf c}$ for some $k\in\bZ$
where ${\mathbf d}\in H_1(\pi^{-1}(\gamma_i),\bZ)$ is such that
${\mathbf c}$ and ${\mathbf d}$ generate $H_1(\pi^{-1}(\gamma_i),\bZ)$
(or a two dimensional subgroup of it if $\gamma_i$ is a loop).
This is equivalent to having $|{\mathbf a}_i\cdot {\mathbf c}|=1$.

If the closure of $\gamma_i$ contains a vertex then an affine circle
representing a class ${\mathbf a_i}$ is a torus knot in an embedded
$S^3$ bounding a neighborhood of the preimage of the vertex.
In fact, it is a torus knot $T(k,l)$ where 
${\mathbf a_i}=k{\mathbf b}_1+l{\mathbf b}_2$ where 
${\mathbf b}_1,{\mathbf b}_2$ are the vanishing classes of the two
edges that meet at the vertex.
But a torus knot $T(k,l)$ is trivial if and only if $k$ or $l$ equals $\pm 1$
and $|{\mathbf b}_i\cdot {\mathbf b_j}|$ equals $0$ if $i=j$ and
$1$ if $i\ne j$.
Therefore $k$ or $l$ equals $\pm 1$ if and only if ${\mathbf a}\cdot 
{\mathbf b}_i=\pm 1$ for $i=1$ or $2$.

The sufficiency of the compatibility conditions for the existence of
$\Sigma_\gamma$ can be checked directly in the 
the local models for neighborhoods of the singular fibers.
Because the family of affine circles in a nodal fiber that represent 
the vanishing class is connected 
and because the node is the only isolated point that is the limit of
affine circles in neighboring regular fibers,
the presence of a node in $\gamma$ does not disrupt the
surface $\Sigma_\gamma$ being well-defined up to the stated fiberwise isotopy.
Furthermore, any fiberwise isotopy of $\Sigma_{\gamma_i}$ extends to its
closure in an obvious way.
\end{proof}

In light of our goal of reading topological and symplectic properties of
a manifold from a base diagram, we want to use
Theorem~\ref{visiblesurfaces.thm} to
reconstruct a visible surfaces from a curve in a base diagram.  
Again, for simplicity of exposition we assume
that our base is a subset of $(\bR^2,\cA_0)$.

Recall that a class in the first homology of a regular fiber $F_b$ is
determined by a covector in $T_b^* B$.  
Given a base diagram we can then use the standard inner product on
$TB\subset T\bR^2$ to represent the covector by its dual in $T_b B$.
Therefore, to record the choices of compatible classes ${\mathbf a}_i$
we assign a {\it compatible vector} $v_i\in\bR^2$ to $\gamma_i$ such that the integral
curves of the vector field $v_i\frac{\del}{\del q}\subset \Lambda^{\rm vert}$ 
represent ${\mathbf a}_i$.
If $\gamma_i$ contains a node then $v_i$ is forced on us so we do
not need to label $\gamma_i$ with that vector.

For reading base diagrams, and in particular for determining whether or 
not vectors $v_i$ are compatible, the following direct calculation
is helpful:
\begin{prop} \label{vectorcompat.prop}
Suppose $F_b$ is a regular fiber in an almost toric fibration defined
by $(B,\cA_0,\cS)\subset (\bR^2,\cA_0)$.
Given two elements ${\mathbf a},{\mathbf b}\in H_1(F_b,\bZ)$
and vectors $v,w\in \bR^2$ such that 
 $[v\frac{\del}{\del q}]={\mathbf a}$ and $[w\frac{\del}{\del q}]
={\mathbf b}$ we
have $|{\mathbf a}\cdot{\mathbf b}|=|\det(v w)|=|v\times w|$.
\end{prop}
Also recall Observations~\ref{basecollapsing.obs} and~\ref{basevanishing.obs}.

Of course, if the image of $\gamma_i$ in a base diagram intersects
the branch locus, the compatible vector $v_i=(k_i,l_i)$ 
will not stay constant as $\gamma$
crosses this locus (unless it is the eigenvector of the monodromy for
the corresponding node).
However, the change in the vector is determined by the base diagram
and therefore it suffices to indicate $v_i$ at just one point of $\gamma_i$.

\begin{rmk} The above restrictions 
that all nodes have multiplicity one and that no boundary component of
$\bdy\Sigma_\gamma$ project to a node were made only to simplify 
reconstruction of $\Sigma_\gamma$ from a base diagram.
With greater willingness to add data to base diagrams these could
be eliminated.
\end{rmk}

\subsection{Symplectic area} \label{area.sec}
It is easy to read the (signed) 
symplectic area of a visible surface off of its image in a base diagram.  
The toric case is well known: If $\gamma$ is an edge that is a smooth
component of the reduced boundary, then the area of $\pi^{-1}(\gamma)$ is
$2\pi|a|$ if the edge is $av$ for some primitive integral vector $v$. 
Indeed, by applying an element of $GL(2,\bZ)$ to the base we can
assume that a neighborhood of $\gamma$ is as in Example~\ref{spherenbhd.ex}
in which case the area of $\pi^{-1}(\gamma)$ equals $2\pi|a|$
where $|a|$ is the length of $\gamma$.

In general we need a way to keep track of orientations.
\begin{defn}   \label{coorient.def}
Consider a parameterized curve $\gamma:I\rightarrow (B,\cA,\cS)$, its
image in a base diagram and a 
set of compatible vectors $\{v_i\}$ that define 
a visible surface $\Sigma_\gamma$.
Let $b=\gamma(t)$ be a point on the interior of $\gamma_i$ and
$x\in M$ such that $\pi(x)=b$.
Let $\tilde u,\tilde v_i\in T_x\Sigma_\gamma\subset T_xM$ 
be as in Theorem~\ref{affine.thm}
where $u=\gamma'(t)\in T_b\bR^2$.
Then the $\{v_i\}_{i=1}^N$ are {\it co-oriented} 
if each ordered pair $(\tilde u,\tilde v_i)$
defines the chosen orientation of $\Sigma_\gamma$.
\end{defn}

Then we have:
\begin{prop} \label{area.prop}
Let $\gamma:I\rightarrow (B,\cA,\cS)$ be a parameterized
immersed curve and $\{v_i\}_{i=1}^N$ a set of 
co-oriented compatible vectors in a base diagram
that define an oriented surface $\Sigma_\gamma$.
Then the (signed) symplectic area of $\Sigma_\gamma$ is
\begin{equation}
 {\rm Area}(\Sigma_\gamma) = 2\pi\int_0^1 \gamma'(t)\cdot v(t) dt 
\end{equation}
where $v(t)=v_i$ if $\gamma(t)\in \gamma_i$ and for other values of
$t$ (when $\gamma\subset\bdy_R B$) $v(t)$ is in integral vector
such that $u(t)\times v(t)=1$ for some integral vector 
$u(t)=\lambda \gamma'(t)$, $\lambda>0$.
\end{prop}
Notice that 
when $\gamma(t)$ belongs to a component of $\bdy_R B$ that intersects
$\overline\gamma_i$ we can choose $v(t)$ to be equal to $v_i$.

\begin{proof}
With respect to standard coordinates we can parameterize
the interior of each $\Sigma_{\gamma_i}$ as 
$p=\gamma_i(t)$, $q=v_i s$ 
for $0\le \alpha_i<t<\beta_i\le 1$ and $0\le s\le 1$.
Then since $\omega=dq\wedge dq$ it is immediate that
${\rm Area}(\Sigma_{\gamma_i})=\int_{\alpha_i}^{\beta_i}\int_0^1 \gamma'_i(t)\cdot v_i dt ds=\int_{\alpha_i}^{\beta_i}\gamma'_i(t)\cdot v_i dt$.
As for the components of $\gamma$ that lie in a component of the
$1$-stratum of $\bdy_R B$, each one
lifts to a visible surface in the boundary recovery $B\times T^2$.
Supposing that the collapsing class of this component of the stratum is
${\mathbf b}$, the visible surface is well defined up to isotopy and 
choice of class ${\mathbf a}$ such that $|{\mathbf a}\cdot {\mathbf b}|=1$.
With $v\frac{\del}{\del q}$ being 
the tangent vector to integral curves representing
${\mathbf a}$, the inner product $\gamma'(t)\cdot v$ is independent
of our choice of ${\mathbf a}$ since for two choices ${\mathbf a},
{\mathbf a'}$, $({\mathbf a}-{\mathbf a'})\cdot {\mathbf b} = 
(v-v')\cdot \gamma'(t)=0.$
\end{proof}

\begin{cor}
If a visible surface $\Sigma_\gamma$ is Lagrangian then $\gamma$ is
an affine linear one-dimensional submanifold 
that contains no open subsets of the $1$-stratum of $\bdy_R B$.
\end{cor}

\begin{proof}
For $\Sigma_\gamma$ to be Lagrangian we must have that $\gamma'(t)\cdot v(t)=0$
for all $t\in I$.  Since $v(t)$ is an integral vector valued function it
is piecewise constant, therefore, because $\gamma$ is smooth we have that 
$\gamma$ is linear.  The constraint $\gamma'(t)\cdot v(t)=0$ and the
integrality of $v(t)$ then force $\gamma$ to be affine linear.
Furthermore, $\gamma$ cannot contain an open subset of the $1$-stratum of 
$\bdy_R B$ since for such points we would have $\gamma'(t)\cdot v(t)>0$.
\end{proof}

\subsection{Intersections numbers} \label{selfint.sec}

In this section we explain how to calculate the intersections and
self-intersections of visible surfaces.
The latter are of particular importance from the symplectic standpoint
because  the symplectic neighborhood
theorem (cf. \cite{McDuffSalamon.intro}) implies that the germ of a
symplectic neighborhood of a symplectic surface in a symplectic four-manifold
is determined by its self-intersection and area.

We start by calculating, separately, 
the signed transverse geometric intersection of visible
surfaces in the preimage of a regular point or a node.

\begin{lemma} \label{regularintersec.lem}
Consider two oriented visible surfaces $\Sigma_{\gamma_i}$, $i=1,2$,
defined by parameterized arcs $\gamma_i:I\rightarrow B_0\subset
(B,\cA,\cS)$
and co-oriented compatible vectors $v_i$ in a base diagram. 
Suppose $\gamma_1$ and $\gamma_2$ 
intersect transversely at one point $b\in B_0$ and
assume that $\Sigma_{\gamma_1}$ and $\Sigma_{\gamma_2}$ intersect
transversely in $F_b$.
Let $u_i$ be a tangent vector at $b$ that is positively oriented
with respect to the parameterization of $\gamma_i$. 
Then $\Sigma_{\gamma_1}$ intersects $\Sigma_{\gamma_2}$ in 
$|\det(v_1 v_2)|$ points in $F_b$ 
and the sign of all of the intersections equals
the sign of the product $\det(u_1 u_2)\det(v_1 v_2)$.
\end{lemma}

\begin{proof}
The vectors $v_i$ determine homology classes $[v_i\frac{\del}{\del q}]$ 
in $H_1(F_b,\bZ)$
whose algebraic intersection in $F_b$ is 
$[v_1\frac{\del}{\del q}]\cdot[v_2\frac{\del}{\del q}]=v_1\times v_2=
\det (v_1v_2)$.
Since the representatives of these classes are affine circles the
algebraic intersection equals the geometric intersection.
To figure the sign of the intersection of the two surfaces in the
four-manifold from the intersection of curves in the fiber torus,
note that the sign of the intersection equals the sign of
$\det(\tilde u_1 \tilde v_1 \tilde u_2 \tilde v_2)\
= \det(u_1 u_2)\det(v_1 v_2)$.
\end{proof}

Notice that two visible surfaces can have a transverse intersection
at a node only if they share an endpoint at that node.  Indeed, as we
saw in the proof of Theorem~\ref{visiblesurfaces.thm}, if $\gamma$
has a node on its interior then the preimage of that point in $\Sigma_\gamma$
is an affine circle representing the vanishing class.  The only transverse
intersection two such circles,
or one such and the node, can have is an empty intersection.
\begin{lemma} \label{nodalintersec.lem} 
The signed intersection of a pair of vanishing disks that intersect once
transversely at a nodal singular point is $-1$.
\end{lemma}

\begin{proof}
To calculate this intersection we refer to the local model in 
Section~\ref{nodalnbhd.sec}.
There we have used complex coordinates ${\mathbf x}=x_1+ix_2$,
${\mathbf y}=y_1+iy_2$ on a neighborhood of the nodal singularity (at
the origin) and the total space has the orientation
given by $({\mathbf x},{\mathbf \overline y})$.
Since the neighborhood of the nodal fiber supports the 
$S^1$ action defined by 
$t\cdot({\mathbf x},{\mathbf y})=
(e^{2\pi i t}{\mathbf x},e^{2\pi i t}{\mathbf y})$(and its
continuation into the $({\mathbf u},{\mathbf v})$ coordinate chart),
the vanishing cycles can be chosen to be $T(1,1)$ torus knots in
an $S^3$ that is the boundary of a neighborhood of the nodal singularity.
These intersect each other positively once with respect to the usual
orientation on $\bC^2$, but since we have the reverse orientation 
the intersection is $-1$.
\end{proof}

\begin{xca}   Check that the only transverse intersections of
visible surfaces that lie in the preimage of the $1$-stratum are empty.
\end{xca}

We now consider how to calculate the self-intersection of a visible
sphere or torus $\Sigma_\gamma$ from its image $\gamma$ in a base diagram.
First, if possible, find
an isotopic visible surface $\Sigma_\eta$ such that
the intersection of $\Sigma_\gamma$ and $\Sigma_\eta$
(hence the self-intersection of $\Sigma_\gamma$) can be
calculated using Lemmas~\ref{regularintersec.lem}
and \ref{nodalintersec.lem}. 
If such a $\Sigma_\eta$ does not exist (as when $\gamma=\bdy B=\bdy_R B$)
then calculate the
self-intersection by applying the above lemmas to the intersections
of two surfaces $\Sigma_{\gamma_1}$, $\Sigma_{\gamma_2}$ that are 
isotopic to $\Sigma_{\gamma}$.   

We make the calculation explicit for the two most useful cases. 
\begin{prop} \label{sphereint.prop}
Suppose $\Sigma_\gamma$ is a sphere such that $\gamma\subset\bdy_R B$
(so $\gamma$ is an edge).
Let $v_1,v_2$ be the inward pointing collapsing vectors for the two components
of the $1$-stratum whose closures intersect the endpoints of $\gamma$.
Let $v_\gamma$ be the inward pointing collapsing vectors for the component
of the $1$-stratum contained in $\gamma$.
Then the self-intersection of the sphere $\Sigma_\gamma$ equals
$v_1\times v_2$ where the vectors are indexed so that $v_1\times v_\gamma>0$.
\end{prop}

\begin{figure}
\begin{center}
        \psfragscanon
        \psfrag{gamma}{$\gamma$}
        \psfrag{g1}{$\gamma_1$}
        \psfrag{g2}{$\gamma_2$}
        \psfrag{10}{$\left(\begin{smallmatrix} 1 \\ 0 
                                 \end{smallmatrix}\right)$}
        \psfrag{1mk}{$\left(\begin{smallmatrix} 1 \\ 
                                          -k \end{smallmatrix}\right)$}
	\includegraphics[scale=.75]{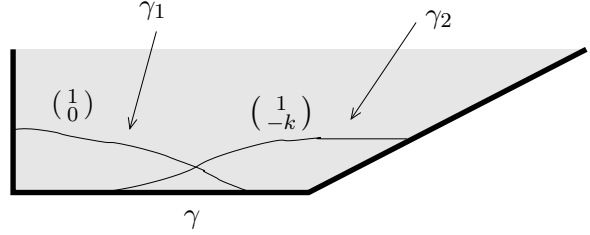}
\end{center}
\caption{Isotopic spheres $\Sigma_{\gamma_1}$ and $\Sigma_{\gamma_2}$.}
\label{sphereisotopy.fig}
\end{figure}

\begin{proof}
Figure~\ref{sphereisotopy.fig} shows two surfaces $\Sigma_{\gamma_i}$,
$i=1,2$ that are isotopic to $\Sigma_\gamma$.  For each the compatible
vector is the collapsing vector for the $1$-stratum that the curve
intersects transversely.
Applying Lemma~\ref{regularintersec.lem} then gives the formula.
\end{proof}

\begin{prop} \label{torusintersec.prop}
Suppose $\Sigma_\gamma$ is a torus with $\gamma\subset\bdy_R B$, so
$\gamma=\bdy_R B=\bdy B$.  Let 
$v$ be an inward pointing collapsing vector 
for the reduced boundary.
Then the self-intersection of $\Sigma_\gamma$ equals $v\times v'$ where $v'$
is the parallel transport of $v$ clockwise along $\gamma$.
\end{prop}

\begin{figure}
\begin{center}
        \psfragscanon
        \psfrag{g1}{$\gamma_1$}
        \psfrag{g2}{$\gamma_2$}
        \psfrag{10}{$\left(\begin{smallmatrix} 1 \\ 0 
                      \end{smallmatrix}\right)$ }
        \psfrag{1m7}{$\left( \begin{smallmatrix} 1 \\ -7 
                      \end{smallmatrix}\right)$ }
	\includegraphics[scale=.75]{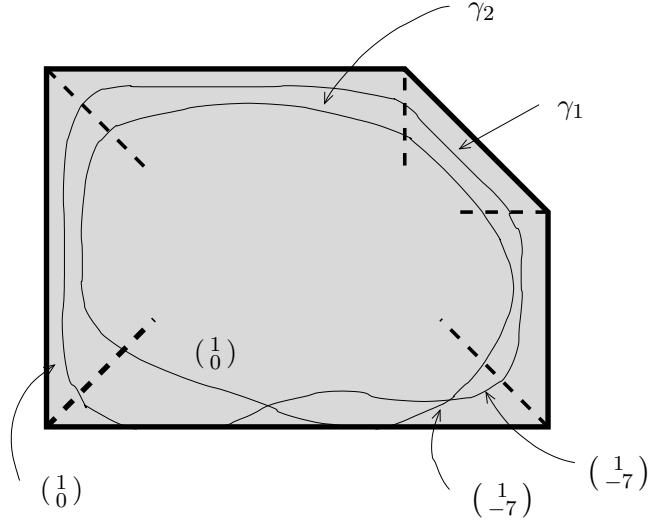}
\end{center}
\caption{Isotopic tori $\Sigma_{\gamma_1}$, $\Sigma_{\gamma_2}$
in $S^2\times S^2\#\overline{\bC P}^2$.}
\label{torusint.fig}
\end{figure}

\begin{proof}
In a base diagram $\gamma$ appears as a union of $k\ge 3$ 
line segments with 
oriented integral tangent vectors $\{u_1,\ldots, u_k\}$ 
such that $u_{j+1} = A_j u_j$ where $A_j$ is conjugate to 
$A_{(1,0)}$ and the vectors are indexed so that $j$ increases as one
moves clockwise along the boundary.
Without loss of generality let 
$v$ be a primitive integral vector normal to the edge with tangent
vector $u_1$.
Now choose two curves $\gamma_1,\gamma_2$ isotopic to $\gamma$ in
a small collar neighborhood of $\gamma$ that contains no nodes.
Choose them so that $\gamma_i\cap\gamma$ for $i=1,2$ each consists of
a closed subset of \lq\lq first edge\rq\rq\ -- 
the segment with tangent vector $u_1$.
Choose compatible vectors $v_1,v_2$ for 
$\Sigma_{\gamma_1},\Sigma_{\gamma_2}$ to equal $u_1$
just to the left of the intersection of $\gamma_i$ with the first edge.
(See Figure~\ref{torusint.fig} for an example.)
Then on the other side of the intersection the compatible vector representing
the same homology class is
$A_k^{-T}\cdots A_2^{-T}A_1^{-T}u_1$.
This equals  $A_k\cdots A_2A_1v$, which is the parallel transport of $v$
along the curves $\gamma_i$, because $v=Ju_1$ and $A_j J=J A_j^{-T}$.
Applying Lemma~\ref{regularintersec.lem} we get the result.
\end{proof}

When the image $\gamma$ of a visible sphere or torus belongs to
$\bdy_R B$ the surface is necessarily symplectic.
A collar neighborhood of $\gamma$ is fibered by parallel affine linear
submanifolds $\eta_t$ 
and the rate of change of the affine length of these submanifolds
measures the self-intersection of the surface, i.e. the first
Chern class of the normal bundle.
Indeed, the pre-image of $\eta_t$ is a circle bundle over the surface
whose first Chern class equals the self-intersection of $\Sigma_\gamma$.
Furthermore, the rate of change of the affine length of $\eta_t$ measures
the Chern class of the normal bundle of the surface (or 
self-intersection of the surface).

Specifically, let $u$ be a primitive integral vector tangent to $\gamma$
and $w$ be an inward pointing integral vector such that $|n\times w|=1$.
Choose the parameter $t$ so that $\eta_t$ passes through points
$b+tw$ for $b\in\gamma$.  Letting $l(t)$ be the affine length of $\gamma$,
$\frac{dl}{dt}$ is a constant and is equal to the self-intersection
of $\Sigma_\gamma$.
We leave it to the reader to verify this alternative formula for
the self-intersection.  It is a special case of the following:
\begin{prop}(cf.~\cite{McDuffSalamon.intro})
Given a circle bundle $\pi:P\rightarrow X$ with first Chern class
$c_1$, then there is an $S^1$ invariant symplectic form $\omega$ 
on the manifold $P\times I$ with Hamiltonian function $H$ equal to the 
projection $P\times I\rightarrow I$ and with reduced spaces $(X,\tau_\lambda)$
where $[\tau_\lambda]=[\tau_\mu]+(\lambda-\mu)c_1$.
\end{prop}

\section{Reading the base II: total space} \label{readingII.sec}

\subsection{Euler characteristic and fundamental group} \label{euler.sec}
Given an almost toric fibration 
we can calculate the Euler characteristic of the total space 
$M$ very easily from the base:
\begin{equation}
\chi(M)=V+\sum_i k_i
\end{equation}
where $V$ is the number of vertices and $k_i$ is the multiplicity of
each node $s_i$.

This follows from the additivity property of the Euler characteristic
and the Euler characteristics of almost toric fibered
neighborhoods of the different types of fibers. 
Each neighborhood has the homotopy type of the singular fiber
so the Euler characteristic of the neighborhood of a point, circle, 
or torus fiber is $1,0$ and $0$ respectively.
Accordingly, the neighborhood of a multiplicity $k$ nodal fiber has
the homotopy type of either a pinched torus (if $k=1$) or a cyclic
chain of $k$ spheres, each of which intersects two others once.
Such a fiber has Euler characteristic $k$.
We can choose to cover the manifold with fibered neighborhoods so that
the intersections of the neighborhoods are all either empty or 
diffeomorphic to a torus bundle or $S^1\times D^3$ and hence always
have Euler characteristic equal to $0$.   

\begin{xca} Note that by encoding the collapsing and vanishing classes,
the base of an almost toric fibration reveals the fundamental
group of the total space.
Verify that any closed toric manifold is simply connected.
\end{xca}

\subsection{First Chern class}

Given a symplectic manifold $(M,\omega)$ its first Chern class
$c_1(M)$ is $c_1(TM)$ where $TM$ is a complex vector bundle
with respect to a compatible almost complex structure on $M$.
(An almost complex structure $J$ is {\it compatible}
if $\omega(u,Ju)>0$ for all $u\ne 0$ and $\omega(Ju,Jv)
=\omega(u,v)$ for all $u,v$.)  
Recall that the set of compatible almost complex structures 
is contractible, so in particular $c_1(M)$ is well-defined.

One way to define the first Chern class of a complex vector bundle of
rank $k$ is as the Poincar\'e dual to the homology class of the set $D$
of points where a generic $k$-frame becomes degenerate (i.e. drops in rank);
in doing so the sign of the homology class is chosen with regard to 
orientations.
For a Lagrangian fibered symplectic manifold we 
bypass choosing an almost complex structure: a representative of
the homology class Poincar\'e dual to $c_1(M)$, up to sign, is given by
the set of points where a generic Lagrangian $k$-frame of $TM$ degenerates.

To get the orientations right we can proceed as follows:
Choose  a generic $n$-tuple of sections $\eta_1,\ldots,\eta_n$ 
that span a Lagrangian subspace of each tangent plane in which they
are linearly independent.
Let $D\subset M$ be the set of points where they are linearly dependent,
and $D'\subset D$ the set of points where they span an $n-1$-dimensional
subspace.
In a neighborhood of a point $x\in D'$
choose two other sections $\mu_1,\mu_2$ such that 
$\omega(\mu_1,\mu_2)>0$ and $\{\eta_j\}_{j<n}\cup \{\mu_i\}$ 
defines a Lagrangian frame for either $i=1$ or $2$.
Then $\eta_n=\sum_{j<n} f_j \eta_j + g \mu_1 + h\mu_2$ for smooth
functions $f_j,g,h$.
Now let $\Omega$ be an orientation for the component of $D'$ containing
$x$ so that $\Omega\wedge dg\wedge dh$ gives the orientation of $M$.
Then the homology class represented by $D$ with the orientations given
by $\Omega_x$ for $x\in D'$ gives the
Poincar\'e dual to $c_1(M)$.
This calculation follows from using an almost complex structure to complete
the Lagrangian frame to a complex frame and then carrying out the usual 
calculation of the orientations (cf.~\cite{GriffithsHarris.alggeom}).
The answer is independent of whatever local compatible almost 
complex structure we choose.

\begin{prop} \label{chernclass.prop}
Consider an almost toric fibration $\pi:(M,\omega)\rightarrow (B,\cA,\cS)$
of a closed manifold $M$.
The  Poincar\'e dual to $c_1(M)$ is
given by the homology class of $\pi^{-1}(\bdy_R B)$.
\end{prop}

\begin{proof}
This is especially easy for a toric manifold: we have standard coordinates
$(p,q)$ defined on $\pi^{-1}(B_0)$.  
There, we can take sections $\eta_i=\frac{\del}{\del q_i}$, $i=1,\ldots, n$.
These become linearly dependent exactly on the preimage of the reduced
boundary and span an $n-1$-dimensional subspace on the preimage of the
$n-1$-stratum of the reduced boundary.
It suffices to consider the case when the collapsing class is 
$[\frac{\del}{\del q_n}]$ in which case we have local coordinates
$(p_1,q_1,\ldots p_{n-1},q_{n-1},x_n,y_n)$.
Taking $\mu_1=\frac{\del}{\del x_n}$, $\mu_2=\frac{\del}{\del y_n}$ we
have $\frac{\del}{\del q_n}=
-y_n\frac{\del}{\del x_n}+x_n\frac{\del}{\del y_n}$ so 
$dg\wedge dh=dx_n\wedge dy_n$.
Therefore the smooth components of the 
preimage of the reduced boundary contribute with their 
orientations as  symplectic manifolds.

Restricting to dimension four,
we do not necessarily have global sections 
$\frac{\del}{\del q_i}$ in an almost toric fibration.
However, on the complement of nodal singularities
$\frac{\del}{\del q_1}\wedge\frac{\del}{\del q_2}$ is well
defined since it is invariant under the action of $GL(2,\bZ)$.
Therefore the Chern class can be calculated using local sections
that give this same section of $\wedge^2 TM$.
Since the nodal singularities are codimension four they do not
contribute to the Chern class calculation.
\end{proof}

\subsection{Volume} \label{volume.sec}
Given an almost toric fibration $\pi:(M,\omega)\rightarrow (B,\cA,\cS)$,
the symplectic volume of $M$ is equal to
the symplectic volume of $\pi^{-1}(B_0)$ where
$B_0$ is the set of regular values.
This is a regular Lagrangian fibration with $\omega=dp\wedge dq$ with
respect to standard local coordinates (which are compatible with the 
fibration).
Therefore the symplectic volume of $M$, $\int_M \omega^n$,
equals $n!\int_{\pi^{-1}(B_0)} dp_1\wedge dq_1\wedge \cdots \wedge dp_n
\wedge dq_n=n!(2\pi)^n \int_{B_0} dp_1\wedge \cdots \wedge dp_n$
since $q$ is defined mod $2\pi$.
(Note that $dp_1\wedge \cdots \wedge dp^n$ gives a well defined volume
form on the base because elements of $AGL(n,\bZ)$ preserve volume.) 
For instance, the symplectic volume of $(\bC P^2,\omega)$ 
where $\omega$ is determined by an almost toric fibration over the
base in Figure~\ref{B4CP2.fig}(b), is $4\pi^2$.

\section{Lens spaces and manifolds they bound}

\subsection{The cone on a lens space}
Two common definitions of the lens space $L(n,m)$ where
$n,m$ are relatively prime integers, $n\ge1$, are:
\begin{itemize}
\item The quotient of the sphere 
$S^3=\{(z_1,z_2)|\abs{z_1}^2+\abs{z_2}^2=1\}\subset C^2$ 
by the free $\bZ_n$ action generated by 
$\alpha\cdot(z_1,z_2)=(\alpha^m z_1,\alpha z_2)$  
where $\alpha=e^{\frac{2\pi i}{n}}$.
\item The space obtained by gluing two solid tori along their boundaries
$T_1,T_2$ as follows:
choose generators $(a_i,b_i)$, $i=1,2$, for $H_1(T_i,\bZ)$
such that $a_i$ is the class of a meridinal circle
in $T_i$ and $a_i\cdot b_i=1$; 
then glue the two solid tori together by a diffeomorphism
$\phi:T_1\rightarrow T_2$ such that $\phi_*a_1=-m a_2+n b_2$.
\end{itemize}
 
For the first definition to make sense, one certainly must assume
that $n\ne 0$, hence $L(0,1)$ is not a lens space. 
However, $n=0$ does make sense in the second definition and one can
easily see that $L(0,1)=S^2\times S^1$.
The orientation of the lens space $L(n,m)$ is the one it inherits from $S^3$
in the first definition and from the solid tori oriented by the 
frame 
$\{\frac{\del}{\del r_i},\frac{\del}{\del \theta_i},\frac{\del}{\del \phi_i}\}$
where $r_i$ is a radial coordinate and $\theta_i,\phi_i$ are cyclic
coordinates such that the meridians $a_i$ are represented
by curves in the boundaries $T_i$ given by 
$\phi_i={\rm const.}$

For our purposes, the following rephrasing of the second definition is
most natural:
\begin{defn} \label{lens.def}   
Let $(q_1, q_2)$ be coordinates on $T^2$, $t$ a coordinate
on $I$ and orient $T^2\times I$ by the frame
$\{\frac{\del}{\del t},\frac{\del}{\del q_1},\frac{\del}{\del q_2}\}$.
The {\it lens space} $L(n,m)$ is the quotient of
$T^2\times I/\sim$ where $\sim$ collapses the integral curves
of $\frac{\del}{\del q_1} $ on $T^2\times \{0\}$
and collapses the integral curves of 
$-m\frac{\del}{\del q_1}+n\frac{\del}{\del q_2}$
on $T^2\times \{1\}$.
\end{defn}
Using Definition~\ref{lens.def} it is very easy to construct a toric fibration
of the cone on a lens space:

\begin{example} \label{Lnm.ex}  
Let $u=(n,m)$ and $v=(0, 1)$ be two vectors in $(\bR^2,\cA_0)$ and
define the domain
$V_{n,m}:=\{xu \pm yv\in \bR^2|x,y\ge 0 \}-\{(0,0)\}$ and its closure
$\overline{V_{n,m}}$.
Then $(V_{n,m},\cA_0,\cS)$ where $\cS$ is specified by 
heavy lines on the two edges of $V_{n,m}$, 
is the base or a toric fibration of $(L(n,m)\times \bR_+,\omega_{n,m})$.
To see this, let $p_2$ be the coordinate on the vertical axis in $\bR^2$;
then Definition~\ref{lens.def} implies immediately
that $\pi^{-1}(V_{n,m}\cap\{p_2=s\})$, $s>0$ is the lens space $L(n,m)$.
Notice that given any embedded curve $\gamma\subset V_{n,m}$ 
which intersects the boundary of $V_{n,m}$ once transversely on each edge,
$\pi^{-1}(\gamma)=L(n,m)$.
\end{example}

\begin{rmk}  \label{orbifold.rmk}
Notice that in general the four-dimensional space defined 
by $(\overline{V_{n,m}},\cA_0,\cS)$ is not smooth.  Indeed, the cone on a lens
space $L(n,m)$ is smooth only when $n=1$, i.e. when the
lens space is $S^3$.  Otherwise it has an {\it orbifold singularity}.
\end{rmk}

The following theorem, which records the redundancy in the above definitions
is well known (cf.~\cite{HirzebruchNeumannKoh.quadraticforms}).

\begin{thm} \label{lens.thm}
Two lens spaces $L(n,m)$ and $L(n,m')$ are orientation preserving
diffeomorphic if and only if
\begin{eqnarray}
m &\equiv&  m' \ ({\rm mod}\ n) \qquad{\rm or} \label{equiv1.eq} \\
mm' &\equiv& 1 \ ({\rm mod}\ n). \label{equiv2.eq} 
\end{eqnarray}
\end{thm}
Note that since the fundamental group of $L(n,m)$ is $\bZ_n$, 
$L(n,m)$ and $L(n',m')$ can only be diffeomorphic if $n=n'$.

We can use base diagrams to \lq\lq see\rq\rq\ these diffeomorphism 
equivalences.
We merely need to apply linear maps $A\in GL(2,\bZ)$ to the base diagram for 
$(V_{n,m},\cA_0,\cS)$ with the vertex at the origin
(e.g. Example~\ref{Lnm.ex}).
If we choose $A$ so that it leaves $v=(0,1)$ fixed then 
$ A=\left(  
\begin{smallmatrix} 
1 & 0 \\  
k &  1
\end{smallmatrix} 
\right).$  
so the image of $V_{n,m}$ is $V_{n,m'}$ where $m'=kn+ m$.
To see the second equivalence, choose $A$ so that it sends $(n,m)$ to
$(0,1)$, i.e.
$ A=\left(  
\begin{smallmatrix} 
-m  & n  \\  
c &  m'
\end{smallmatrix} 
\right)$ with $ \det A=-1$.

\begin{xca}   \label{oreintrev.ex}
Two lens spaces $L(n,m)$ and $L(n,m')$ are orientation reversing
diffeomorphic if and only if  
$m \equiv -m'\ ({\rm mod}\ n)$ or $mm' \equiv -1\ ({\rm mod}\ n)$.
Use base diagrams to verify that such lens spaces are diffeomorphic.
How does the reversal of the orientations show up?
\end{xca}

Toric fibrations also give a way to see the equivalence of
the two definitions of the lens space $L(n,m)$, (1) as a
quotient of $S^3$ and (2) as a union of solid tori.

To explain this, we present a symplectic $n$-fold covering of 
$L(n,m)\times\bR_+$ by $L(1,0)\times\bR_+=S^3\times\bR_+$ 
from which we can extract the quotient
construction.
Indeed, define a linear map $B:\bR^2\rightarrow\bR^2$ by
\begin{equation} 
B=\left(  
\begin{matrix}
 1 & 0 \\  
\frac{m}{n} &  \frac{1}{n} 
\end{matrix} 
\right)  
\end{equation} 
which maps $V_{1,0}$ onto $V_{n,m}$.
The map $(B,B^{-T}): (V_{1,0}\times T^2,\omega_0)\rightarrow
(V_{n,m}\times T^2,\omega_0)$ 
is then a symplectic $n$-fold cover that descends under boundary
reduction to a symplectic $n$-fold cover of $L(n,m)\times \bR_+$ 
by $L(1,0)\times\bR_+$ that preserves the almost toric fibrations
defined by $V_{1,0}$ and $V_{n,m}$.
Restricting this $n$-fold cover to a sphere $S^3$ that projects to an
embedded curve $\gamma\subset V_{1,0}$, we see $L(n,m)$ as the stated quotient
of $S^3$.
Because Definition~\ref{lens.def} is just a rewording of (2),
we are done.

\subsection{An induced contact structure} \label{contact.sec}

A {\it contact structure} $\xi$ on a three-manifold $N$
is a completely non-integrable plane field.  Locally (and often 
globally) such a plane
field can be prescribed as the kernel of a $1$-form $\alpha$ 
such that 
\begin{equation}\label{contact.def}
\alpha\wedge d\alpha\ne0.
\end{equation}
Conversely, given any $1$-form $\alpha$ satisfying Equation~\ref{contact.def}, 
$\xi={\rm ker}\ \alpha$ is a contact structure.

Contact structures are the cousins of symplectic structures that
live on odd dimensional manifolds.
In particular, contact manifolds arise naturally as hypersurfaces in 
symplectic manifolds.

\begin{defn}   \label{expanding.def}
A vector field $X$ on a symplectic manifold $(M,\omega)$ is {\it expanding}
if it satisfies $\cL_X\omega=\omega$.
\end{defn}
Recall that $\cL_X\omega:=d(\iota_X\omega)+\iota_X d\omega$.
In our case, since $d\omega=0$ only the first term is non-zero.

A basic result from contact geometry is:
\begin{prop} \label{contact.prop}
Let $Y\subset (M,\omega)$ be a hypersurface.  If there is an expanding
vector field $X$ defined on a neighborhood of $Y$ that is 
transverse to $Y$,
then $\alpha:=\iota_X \omega$ defines a contact structure on $Y$.
\end{prop}
                                                                       
An almost 
toric fibration of a symplectic manifold allows for a simple criterion
to determine the existence of an induced contact structure on certain 
lens space hypersurfaces.

\begin{defn}   \label{lens.obs}
A {\it visible lens space} in a symplectic manifold with
almost toric fibration $\pi:(M,\omega)\rightarrow (B,\cA,\cS)$
is the preimage of an embedded curve $\gamma\subset (B,\cA,\cS)$ that is
disjoint from nodes and the $0$-stratum  and intersects
the $1$-stratum of $\bdy B$ transversely at its endpoints.
\end{defn}
That $\pi^{-1}(\gamma)$ is a lens space follows directly from 
Definition~\ref{lens.def}.
Notice that a neighborhood of $\gamma$, chosen small enough to
contain no nodes, necessarily admits an integral affine immersion
to $\bR^2$ (since it is simply connected).

\begin{prop} \label{visiblelens.prop}
Let $\gamma\subset(B,\cA,\cS)$ be an embedded curve that defines a
visible lens space and let $\cN(\gamma)$ be a neighborhood of $\gamma$.
Consider an immersion 
$\Psi:(\cN(\gamma),\cA)\rightarrow (\bR^2,\cA_0)$.
If there is a point $x_0\in\bR^2$ such that $\Psi(\gamma)$ is transverse to 
all lines through $x_0$ and those lines through $x_0$ that intersect
the endpoints of
$\Psi(\gamma)$ are tangent to the image of $\bdy_R B$
at those intersections, 
then the lens space has an induced contact structure.
\end{prop}

\begin{proof} 
Let $\Phi_n:(\widetilde W_n,\cA_n)\rightarrow (\bR^2-\{0\},\cA_0)$
be an $n$-fold cover equipped with local coordinates $p=(p_1,p_2)$ 
induced from $\bR^2$.
This defines a Lagrangian fibration
$\pi_n:(\cN(\gamma)\times T^2,\omega_n)\rightarrow (W_n,\cA_n)$
where $\omega=dp\wedge dq$ with respect to local standard coordinates.
Without loss of generality take $x_0$ to be the origin in $\bR^2$.
The transversality assumption on $\gamma$ implies that if $\cN(\gamma)$
is chosen small enough, the immersion $\Psi$ factors through
an embedding 
$\Psi_n:(\cN(\gamma),\cA)\rightarrow (\widetilde W_n,\cA_n)$ 
so $\Psi=\Phi_n\circ \Psi_n$.
Then the vector field $p\frac{\del}{\del p}$  on $(\bR^2-\{0\},\cA_0)$
defines a vector field of the same form on $(\widetilde W_n,\cA_n)$, 
and likewise a vector field $X=p\frac{\del}{\del p}$ on 
$\cN(\gamma)\times T^2\subset (\widetilde W_n\times T^2,\omega_n)$.
Note that $X$ is transverse to $\pi_n^{-1}(\Psi_n(\gamma))$.
Therefore, since $d(\iota_X dp\wedge dq)=d( p dq)=dp\wedge dq$,
$X$ defines a contact structure on 
$\gamma\times T^2\subset (N(\gamma)\times T^2,\omega_n)$.
The assumed tangency at the endpoints of $\gamma$ implies 
this vector field extends to a well defined vector field
on the boundary reduction of $\cN(\gamma)\times T^2$ along
$\pi^{-1}_n \Psi_n(\bdy_R B\cap \cN(\gamma))$.
\end{proof}

In contact geometry on three-manifolds
there is a dichotomy of fundamental importance:
a contact structure can be tight or overtwisted, and from a topological
standpoint the tight structures are the interesting ones.

\begin{defn} \label{tight.def}  
A contact structure $\xi$ on a three-manifold
$Y$ is {\it overtwisted} if there is a 
an embedded (two-dimensional)
disk $D\subset Y$ such $\xi$ is tangent to $D$ at all 
points of $\bdy D$.
It is {\it tight} if there is no such disk.
\end{defn}

For visible lens spaces, there is an easy way
detect overtwistedness of a contact structure.  
The following construction of
overtwisted contact structures was first made by 
Lerman~\cite{Lerman.contactcuts} to illustrate the existence of
a countable number of non-conjugate tori of dimension two in
the group of contactomorphisms of these contact structures on
lens spaces.

\begin{prop} \label{overtwisted.prop}
Suppose $\gamma\subset (B,\cA,\cS)$ determines a visible lens space.
Consider the image of $\gamma$ in $(\bR^2,\cA_0)$ and a point $x_0$ in
$\bR^2$ as in Proposition~\ref{visiblelens.prop}.  If the image of
${\rm Int}(\gamma)$ intersects
any line through $x_0$ more than once, then the induced contact
structure on the lens space is overtwisted.
\end{prop}

\begin{proof}
View $\gamma$ as an immersion $\gamma:[0,1]\rightarrow (\bR^2,\cA_0)$
with tangent vectors $\gamma'(t)$.
We can assume, without loss of generality, that $x_0$ is the origin
in $\bR^2$, that $\gamma(0)$ is on the positive
vertical axis and $\gamma'(0)$ has a positive horizontal component.  
Proposition~\ref{visiblelens.prop} implies the existence of a contact
structure $\xi$ that in terms of local coordinates can be
given as $\xi={\rm ker}\ pdq$.

Consider the visible disk $D=\Sigma_{\gamma([0,\tau])}$, $0<\tau<1$, 
with compatible vector $(1,0)$.  (See Section~\ref{visible.sec} for
definitions.)
The boundary of the disk consists of points that project to $\gamma(\tau)$.
In terms of local coordinates, 
the tangent plane to the disk at these points is spanned by the 
vectors $\gamma'(\tau)=v\frac{\del}{\del p}$, $v\in\bR^2$, 
and $\frac{\del}{\del q_1}$.
Meanwhile, the contact planes are spanned by
$\gamma'(\tau)$ and $(Jp)\frac{\del}{\del q}$.
Here,   
$J=\left(  
\begin{smallmatrix} 0 & -1 \\  1 & 0
\end{smallmatrix} 
\right).$
Thus tangency of $D$ and the contact planes occurs all along the boundary
if we can choose $\tau$ so that $\gamma(\tau)$ lies on the vertical
axis.
Definition~\ref{tight.def} implies that under these conditions the
induced contact structure is overtwisted.
\end{proof}
Note that in the above construction it does not suffice for 
the endpoints of $\gamma$ to lie on a line.
If they do, and ${\rm Int}(\gamma)$ belongs to an open half-space, then
$\gamma([0,\tau])=\gamma([0,1])$ defines a sphere -- and hence is not a
candidate for an overtwisted disk.

\begin{xca}   \label{overtwisted.ex}
Check that the characteristic foliation
(defined in~\cite{Etnyre.contactintro}) 
on the overtwisted disk $D$ in the preceding proof
is exactly the degenerate characteristic
foliation shown in Figure 3 of~\cite{Etnyre.contactintro}.
\end{xca}

Determining tightness of a contact structure is a difficult problem in
general.  One of the most useful criteria is that of fillability:

\begin{defn}   \label{convex.def}
Consider a symplectic manifold $(M,\omega)$ with boundary $Y=\bdy M$.
The boundary is {\it $\omega$-convex} if a collar neighborhood of $Y$ admits
an outward pointing expanding vector field transverse to $N$.
Then the induced contact structure $\xi$ on $Y$ is
{\it (strongly symplectically) fillable} and $(M,\omega)$ is 
a {\it (strong symplectic) filling} of $(Y,\xi)$.
\end{defn}

\begin{thm}[Gromov, Eliashberg; cf.~\cite{Etnyre.convexity}]
\label{fillable.thm}
If a contact three-manifold $(Y,\xi)$ is fillable then the contact structure
$\xi$ is tight.
\end{thm}
There are several notions of fillability of contact manifolds; any one of
them suffices in the above theorem.  For more details 
see~\cite{EtnyreHonda.nofilling}.
\begin{cor} \label{convex.cor}
If the boundary of a symplectic manifold $(M,\omega)$ is $\omega$-convex,
then the induced contact structure is tight.
\end{cor}

\begin{cor} \label{onetoone.cor}
Let $\pi:(M,\omega)\rightarrow(B,\cA,\cS)$ be an almost toric fibration of a
neighborhood of a nodal fiber, $R\subset B$ an eigenray based at the node, and 
$\Phi:(B-R,\cA)\rightarrow(\bR^2,\cA_0)$ an immersion.
Then $\Phi$ is an embedding.
\end{cor}

\begin{proof}
Without loss of generality suppose that the image (under $\Phi$) of the
node is the origin in $\bR^2$.
Let $\gamma\subset (B,\cA,\cS)$ be a loop around the node, the image
of $\bdy M'$ where $M'\subset M$ is a closed almost toric neighborhood
of the nodal fiber.
Suppose without loss of generality that the image of $\gamma\cap(B-R)$
in $(\bR^2,\cA_0)$ is transverse to the vector field $p\frac{\del}{\del p}$.
With respect to standard coordinates on the preimage of $B-R$, 
$p\frac{\del}{\del p}$ is a vector field that extends to 
an expanding vector field $X$ on  a full collar neighborhood of
$\bdy M'$.
Since $X$ is transverse to $\bdy M'$, Theorem~\ref{fillable.thm} implies
the induced contact 
structure $\xi={\rm ker}(\iota_X \omega)$ is tight.
But Proposition~\ref{overtwisted.prop} implies this can happen only
if $\Phi$ is one to one. 
\end{proof}

\subsection{Rational balls} \label{rationalballs.sec}

\begin{figure}
\begin{center}
	\psfragscanon
        \psfrag{k}{$l$}
        \psfrag{l}{$k$}
        \psfrag{l2}{$k^2$}
        \psfrag{klm1}{$kl-1$}
	\includegraphics[scale=.5]{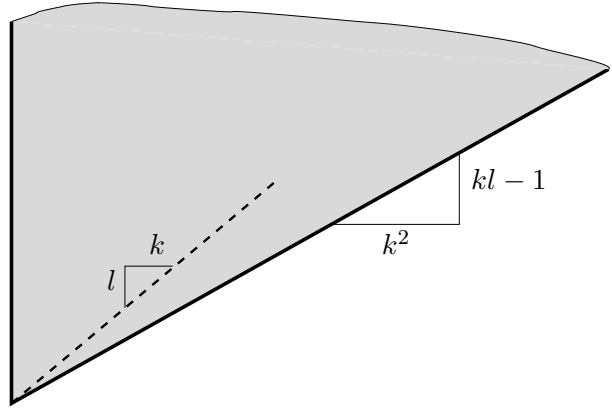}
\end{center}
\caption{Rational ball with boundary $L(k^2,kl-1)$.}
\label{rationalball.fig}
\end{figure}

The base $\overline V_{n,m}$ defines a toric fibration of the cone on
$L(n,m)$ and has an  orbifold singularity whenever $n\ne 1$.
When $n=1$ we have seen in Section~\ref{nodaltrade.sec}
that an another almost-toric fibration (of $\bR^4$,
the cone on $S^3$) can be obtained by trading the vertex for a node.
We can in fact mimic this construction whenever $n=k^2, m=kl-1$ for  
relatively prime integers $k,l$ because the linear map $A_{l,k}$ (defined
in Equation (\ref{Aab.eq})) maps the vector $(k^2,kl-1)$ to $(0,-1)$.
Specifically, modify the base diagram for $\overline V_{n,m}$
by removing the vertex and including a ray $R$ marked by heavy dashes,
based at $(tk,tl)$ for some $t>0$ and passing through the origin.
The resulting base $\widetilde V_{n,m}$ is the base of an almost toric
fibration of a smooth manifold.
(See Figure~\ref{rationalball.fig}.) 

What is the manifold?  The structure of the base implies that
this manifold, say $M$, is obtained from $T^2\times D^2$ by forcing
$(-l,k)$ to be a vanishing class (thereby introducing a nodal fiber)
and having $(1,0)$ be a collapsing class in a boundary reduction.
Since these two vectors span $\bQ^2$, we have $H_*(M,\bQ)=H_*(\bR^4,\bQ)$. 
The construction also reveals that $H_1(M,\bZ)=\bZ_k$. 

Now consider any embedded arc $\gamma$ with vertices on the two edges of the
base diagram for $\widetilde V_{n,m}$.
If $\gamma$ does not intersect the ray $R$ then it defines a visible
lens space $L(k^2,kl-1)$ that bounds a compact manifold with boundary
whose homology with rational coefficients is that of a ball;  this
manifold with boundary, denoted here as $B_{k,l}$, 
is called a {\it rational ball}.  
Invoking Proposition~\ref{visiblelens.prop}, the
boundary is equipped with an induced tight contact structure if
$\gamma$ is transverse to the radial vector field $p\frac{\del}{\del p}$.

\begin{xca}   \label{S1xS2.ex}
Check that if $\gamma$ intersects the ray $R$ once (away from the node)
then it defines a visible surface that is diffeomorphic to $S^1\times S^2$.
Does the transversality of $\gamma$ with the vector field
$p\frac{\del}{\del p}$ ensure an induced tight contact structure?
\end{xca}

\subsection{Chains of spheres} \label{chains.sec}

The rational balls of the previous section are almost toric 
manifolds bounded by the lens spaces $L(k^2,kl-1)$.
In fact, every lens space bounds a smooth manifold that can be equipped
with a symplectic structure and a toric fibration.

\begin{defn} \label{chain.def}  
A {\it chain of spheres} in a smooth four-manifold is a collection
of embedded spheres $S_i$, $i=1,\ldots k$ such that
$S_i\cdot S_{i+1}=1$ for $i=1,\ldots, k-1$, $S_i\cdot S_j=0$ for 
$j\ne i, i+1$, and all intersections are transverse.
\end{defn}

A model for a neighborhood of a chain of spheres in a four-manifold is 
a linear plumbing of disk bundles:

\begin{defn}   \label{plumbing.def}
A {\it plumbing} of two disk bundles over spheres
$\pi_i:X_i\rightarrow \Sigma_i$ is the space $(X_1\cup X_2)/\varphi$
where $\varphi:\pi_2^{-1}(D_2)\rightarrow \pi_1^{-1}(D_1)$
is an orientation preserving diffeomorphism that maps the disk 
$D_2\subset \Sigma_2$ to a fiber of $\pi_1^{-1}(D_1)$ and maps a fiber
of $\pi_2^{-1}(D_2)$ to the disk $D_1\subset \Sigma_1$. 
A {\it linear plumbing} of disk bundles $\pi_i:X_i\rightarrow \Sigma_i$,
$i=1,\ldots k$, is the space $\sqcup_i X_i/{\varphi_1, \ldots\varphi_{k-1}}$ 
where
$\varphi_i:\pi_{i+1}^{-1}(D_{i+1})\rightarrow \pi_{i}^{-1}(D'_{i})$ 
is as above and
for each $i\in\{2,\ldots, k-1\}$,  
$\pi_i^{-1}(D_i)$ and $\pi_i^{-1}(D'_i)$ are disjoint. 
\end{defn}

The following theorem is well known 
(cf.~\cite{HirzebruchNeumannKoh.quadraticforms}).  
Here we give a proof in terms of base diagrams.
\begin{thm}
Let $N$ be  a closed 
neighborhood of a chain of spheres $\{S_i\}_{i=1,\ldots k}$ in a four-manifold
such that the spheres $S_i$ have 
self-intersections $-b_i$.  Then $\bdy N=L(n,m)$ where 
\begin{equation}
\frac{n}{m}=[b_1,b_2,\ldots b_k]:=b_1-\frac{1}{b_2-\cdots\frac{1}{b_k}}.
\end{equation}
\end{thm}

\begin{proof}
Let $U_i$ be the quadrilateral in $\bR^2$ that is the intersection
of the open half-plane $\{p_2<\varepsilon\}$ and
$U_{n,1}$ of Example~\ref{spherenbhd.ex} with $n=b_i$ and $s=a_i$.
Thus the base diagram $U_i$ defines
a model neighborhood of a sphere with self-intersection $-b_i$.
Now \lq\lq plumb\rq\rq\ the base diagrams $U_i$, $i=1,\ldots k$.
Specifically, let $u_i=(a_i,0)$ and 
$A_i=\left(  
\begin{smallmatrix} 
b_i & -1 \\  
1 & 0 
\end{smallmatrix}
\right)$.
Let $D_i,D_i'$ be two dimensional polydisks such as 
shown in Figure~\ref{polydisks.fig}
and then construct
$B=\sqcup_{i=1}^k U_i/{\phi_1,\ldots \phi_{k-1}}$ 
where $\phi_i:D_{i+1}\rightarrow D'_i$
is defined by $\phi_i(x)=A_ix+u_i$.
The parameter $\varepsilon$ is assumed to be chosen small enough that the
$D_i,D'_i$ are disjoint.
Then Proposition~\ref{basepasting.prop} implies we have constructed
a base that defines a plumbing of the disk bundles defined by the
$U_i$.
An example is shown in Figure~\ref{chain.fig} with $k=3$, $b_1=2$,
$b_2=1$, and $b_3=3$.
Finally, an embedded smooth curve $\gamma\subset B$ that projects to
a curve with endpoints on the \lq\lq initial\rq\rq\ and \lq\lq final\rq\rq\ 
edges (as in Figure~\ref{chain.fig})
defines a visible lens space
that is the boundary of a closed neighborhood of the
spheres.

\begin{figure}
\begin{center}
	\psfragscanon
        \psfrag{Di}{$D_i$}
        \psfrag{Dip}{$D'_i$}	
	\includegraphics[scale=.75]{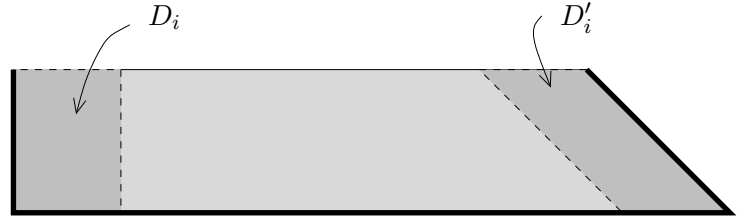}
\end{center}
\caption{Polydisks $D_i,D'_i$.}
\label{polydisks.fig}
\end{figure}

\begin{figure}
\begin{center}
	\psfragscanon
        \psfrag{g}{$\gamma$}
	\includegraphics[scale=.75]{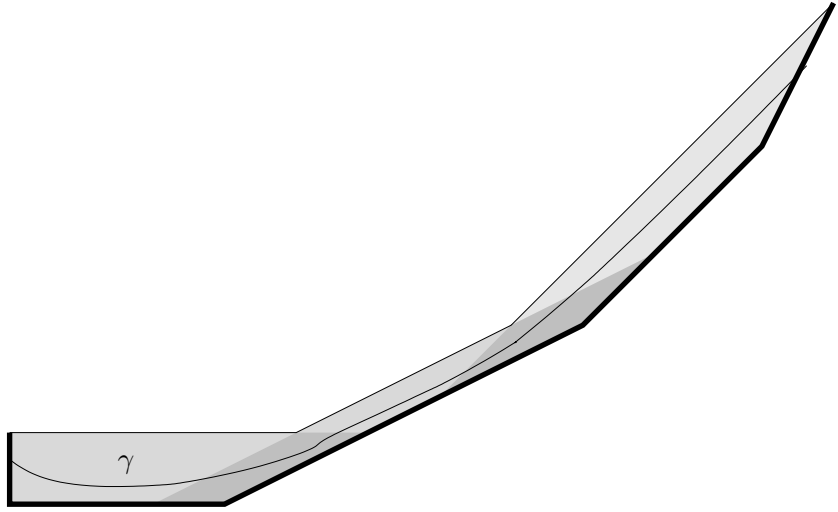}
\end{center}
\caption{Chain of spheres in a linear plumbing.}
\label{chain.fig}
\end{figure}

To calculate which lens space bounds this neighborhood, 
choose an immersion $\Phi:(B,\cA)\rightarrow (\bR^2,\cA_0)$
such that $\Phi(U_1)=U_{b_1,1}$.
Then  notice that the image (under $\Phi$) of the final vector is 
\begin{equation} 
 A_1\cdots A_{k-1} \left(  
\begin{matrix}
b_k  \\  
1  
\end{matrix} 
\right).  
\end{equation} 
Now, 
\begin{equation} 
 A_j \left(  
\begin{matrix} 
x \\  
y  
\end{matrix} 
\right)=
\left(  
\begin{matrix} 
b_jx-y \\  
x
\end{matrix} 
\right)
\end{equation} 
which changes the ratio $\frac{x}{y}$ to $b_j-1/\frac{x}{y}$.
The conclusion of the theorem follows immediately.
\end{proof}

\begin{xca}   \label{blowdown.ex}
Use base diagrams to check that if $b_i=-1$ for some $1\le i\le k$
then $[b_1,\ldots, b_{i-1},b_i,b_{i+1},\ldots,b_k]=
[b_1,\ldots, b_{i-1},b_{i+1},\ldots,b_k]$.
\end{xca}

\begin{xca}   \label{convexity.ex}
Consider a linear plumbing that is a model neighborhood for a chain
of spheres $\{S_i\}_{i=1,\ldots k}$ with self-intersections
$b_i\ge 2$.
Show that inside the linear plumbing 
there is an arbitrarily small closed toric neighborhood of the
chain of spheres
whose boundary is $\omega$-convex.  (See Definition~\ref{convex.def}.)
\end{xca}

\section{Surgeries} \label{surgeries.sec}

\subsection{Rational blowdown} \label{blowdown.sec}

Almost toric fibrations of open manifolds and manifolds with boundary
are very useful for understanding surgeries of symplectic manifolds.
The key is that the gluing locus for such surgeries often admit an 
almost toric fibration.

Blowing up and down  are the simplest examples of this.  
One need not have an almost toric fibration of a symplectic manifold
$(M,\omega)$ in order to \lq\lq see\rq\rq\ a blow-up or blow-down;
one only needs a toric fibration of an embedded symplectic ball
or a neighborhood of a symplectic $-1$-sphere.  
In the first case one can push forward the fibration of the standard
ball that is being embedded; in the second case the symplectic neighborhood
theorem guarantees that the $-1$ sphere has a neighborhood symplectomorphic
to a toric fibered model.
The fact that balls and spheres of the appropriate sizes have fibrations
whose bases are equivalent near their boundaries gives a \lq\lq toric
proof\rq\rq\ that this surgery can be achieved. 

More general than blowing down, in which one replaces a neighborhood
of a $-1$-sphere by a ball, one can replace certain chains of
spheres (Section~\ref{chains.sec}) by a rational ball.
This surgery is called a (generalized)  rational blowdown:

\begin{defn}   \label{grbd.def}
If a smooth four-manifold $M$ contains a chain of spheres
such that the boundary of a closed neighborhood $N$ of the spheres is the lens
space $L(n^2,nm-1)$,
then the {\it generalized rational blowdown} of $M$ 
along the spheres is
\begin{equation}
\wtilde M :=(M-{\rm Int} (N))\cup_\varphi B_{n,m}
\end{equation}
 where $\varphi$ is
an orientation preserving diffeomorphism from the boundary of $N$
to the boundary of the rational ball $B_{n,m}$.
\end{defn}

This smooth surgery (in the case $m=1$) was used by Fintushel
and Stern to prove the diffeomorphism classification of elliptic surfaces
and to construct interesting smooth four-manifolds including
an infinite family of simply connected four-manifolds not homotopic
to a complex surface.
They also used it to give an alterative construction of the
exotic K3 surfaces of Gompf and Mrowka (manifolds homeomorphic but not
diffeomorphic to the K3 surface).

Almost toric fibrations yield a virtually pictorial proof 
of

\begin{thm}[Symington~\cite{Symington.grbd}] \label{grbd.thm}
If a symplectic four-manifold $(M,\omega)$ 
contains a chain of spheres, each of which
is symplectic, such that the boundary of
a closed neighborhood of the spheres is the lens
space $L(k^2,kl-1)$, $k\ge 2$ and $l\ge 1$,
then the {\it generalized rational blowdown} of $(M,\omega)$ 
along the spheres admits a symplectic structure induced by $\omega$.
\end{thm}

Indeed, we constructed in Section~\ref{chains.sec} a toric 
neighborhood of a chain of spheres whose boundary is $L(n,m)$.
Taking $n=k^2$, $m=kl-1$,
a slight generalization of the symplectic neighborhood theorem implies
that the chain of spheres in the theorem have a neighborhood 
symplectomorphic to our model toric one (for some curve $\gamma$).
We also constructed, in Section~\ref{rationalballs.sec}, an almost toric
fibration of a rational ball.
The structure of the bases of these fibrations makes it clear that whenever
a rational ball with boundary $L(k^2,kl-1)$ can replace a neighborhood of
a chain of spheres smoothly then this can be done symplectically.

As a consequence of Theorem~\ref{grbd.thm}, the smooth manifolds
mentioned above (constructed by Fintushel and Stern) all admit
symplectic structures.

\subsection{Symplectic sum}

The symplectic sum was defined in~\cite{Gromov.pdr} by
Gromov and used to great effect (mostly in
dimension four) by Gompf~\cite{Gompf.anewconstr} 
and Fintushel and Stern~\cite{FintushelStern.knotslinks}.

\begin{defn} \label{fibersum.def}  
Consider two symplectic manifolds $(M_i,\omega_i)$, $i=1,2$,
that contain symplectomorphic codimension two submanifolds $\Sigma_i$ with 
anti-isomorphic normal bundles.
Let $\overline{M_i-\Sigma_i}$ be the manifold with boundary such
that $\overline{M_i-\Sigma_i}-Y_i$ is symplectomorphic to 
$(M_i-\Sigma_i,\omega)$ with $Y_i\subset \bdy\overline{M_i-\Sigma_i}$
being diffeomorphic to a circle bundle over $\Sigma_i$.
The {\it fiber sum} $M_1\#_{\Sigma_1=\Sigma_2}M_2$ is then
$\overline {M_1-\Sigma_1}\cup_\varphi \overline {M_2-\Sigma_2}$
with the induced symplectic structure, where 
$\varphi:Y_1\rightarrow Y_2$ is an orientation reversing diffeomorphism.
\end{defn}

In dimension four, the normal bundles of surfaces $\Sigma_i$ are
anti-isomorphic when their self-intersections sum to zero and
the surfaces themselves are symplectomorphic if they have
the same genus and area.

When the surfaces are spheres base diagrams give a way to picture 
the symplectic sum.
Specifically, if the surfaces are spheres of self-intersection 
$n$ and $-n$, 
invoking the symplectic neighborhood theorem we know that we can
find neighborhoods of the $\Sigma_i$ that admit toric fibrations with bases
$U_{n,1}$, $U_{-n,1}$ of Example~\ref{spherenbhd.ex}
for some choice of $s$.
Change the stratification in each base by making the interior of the edge
that is the image of $\Sigma_i$ part of the non-reduced boundary.
This gives  base diagrams that define toric fibrations of  collar
neighborhoods of the $Y_i$ in $(\overline{M_i-\Sigma_i},\omega)$.
Under the hypotheses of the theorem,  Proposition~\ref{basepasting.prop}
implies the bases for the collar neighborhoods of the $Y_i$ can 
be joined along the images of the $Y_i$ to yield a base
for a toric fibration of the gluing locus in the summed manifold.

\begin{rmk}   \label{symplsum.rmk}
When performing the symplectic sum of symplectic four-manifolds it is
not necessary for the surfaces $\Sigma_1,\Sigma_2$ to be symplectomorphic.
Unless they have self-intersection zero, it suffices that they have the same 
genus and opposite self-intersection numbers, and that
one be willing to thicken or thin one of the $(M_i,\omega_i)$, $i=1$ or $2$,
along $\Sigma_i$.
\end{rmk}

If $\Sigma_1,\Sigma_2$ are symplectic 
tori then they have almost toric neighborhoods.
Modifying $U_{n,1}, U_{-n,1}$
by replacing the left and right edges by light
lines marked with arrows to indicate that they are identified
yields base diagrams for such  neighborhoods.
Then, as above, carrying out the symplectic sum along these surfaces
corresponds to replacing the reduced boundary component in each
base by a non-reduced boundary component and then identifying these two
boundary components by an affine isomorphism.

Of course, if an almost toric fibration is given and the $\Sigma_i$
are spheres or tori
belonging to the preimage of the reduced boundary in each $M_i$, then the
above procedure shows the symplectic sum can be performed
within the category of almost toric fibered manifolds.

\begin{xca}   \label{symplsum.ex}
If the surfaces $\Sigma_i$ have genus $\ge 2$ then they do not have
an almost toric neighborhood.  Nonetheless, one can still draw
base diagrams that encode the sum.
We leave it to the interested reader to sort this out.
Hint: Decompose each surface as $\Sigma_i = D_i\cup \Sigma'_i$
where $\Sigma'_i=\Sigma_i-D'_i$ and $D'_i\subset D_i$ are disks such
that the areas of $\Sigma_i-D'_i$ are equal and so are the areas
of the $D_i-D'_i$.
Find disk bundle neighborhoods of the $\Sigma_i$ and local coordinates
that trivialize the bundles over the $\Sigma-D'_i$ with compatible almost toric
fibrations of the bundles over the $D_i-D'_i$ and the $D_i$.
The base diagrams for the fibrations of the disk bundles over the $D_i$
then reveal the essence of the symplectic sum.
\end{xca}

\section{Images of the K3 surface} \label{K3.sec}

The base of a Lagrangian fibered K3 surface is a sphere so there are
no base diagrams per se. 
However, using symplectic sums we can give 
presentations of the K3 as a union of fibrations over topological
disks for which we can draw base diagrams.

\begin{figure}
\begin{center}
	\includegraphics[scale=.75]{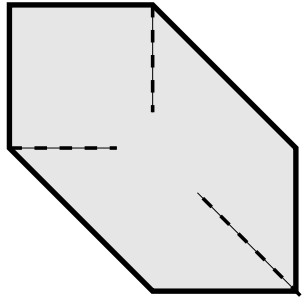}
\end{center}
\caption{Fibration of $\bC P^2\#3\overline{\bC P}^2$.}
\label{octant.fig}
\end{figure}

Our first construction was initially observed by Zung~\cite{Zung.II}.
Consider the base diagram shown in Figure~\ref{octant.fig} that defines
a Lagrangian fibration  $\pi:(M,\omega)\rightarrow (B,\cA,\cS)$ where 
$M=\bC P^2 \#3\overline{\bC P}^2$
The essence is to perform successive
symplectic sums of eight copies of this manifold by summing eight copies of
$(B,\cA,\cS)$, i.e. by performing symplectic sums that are compatible
with the Lagrangian fibrations.
Specifically, the reduced boundary shown in
Figure~\ref{octant.fig} is a union of three rational line segments that
are the images of three symplectic
spheres of self-intersection zero, all having equal area.
If we sum two copies of the base 
we get a second base that defined a fibered manifold with six nodal
fibers and two spheres of self-intersection zero that project to the
boundary of the base.
(The new spheres are connect sums of two pairs of the original ones.)
Performing a symplectic sum of two copies of this second base
yields a third; it defines a fibration over a disk with $12$ nodes and one
smooth component to the reduced boundary.  
(As a shortcut, one can sum four copies of the original manifold at once
using the four-fold sum~\cite{McDuffSymington.assoc}, also in terms
of gluing the bases together.)
Lemma~\ref{diskbase.lem} and Theorem~\ref{toricfourmanifolds.thm}
together imply that the resulting manifold is diffeomorphic to E(1).
The preimage of the boundary is a symplectic torus of self-intersection
zero so we can perform a symplectic sum of two copies
to get a fibration over the sphere with $24$ nodal
fibers; by the proof of
Theorem~\ref{closedalmosttoric.thm} this is a K3 surface.

A standard fact about the K3 surface is that its intersection form,
with respect to an appropriate basis for $H_2(K3,\bZ)$,
is $-E_8\oplus -E_8\oplus 3
\left(\begin{smallmatrix} 0 & 1 \\  1 & -2 \end{smallmatrix}\right)$.
Our goal now is to present an almost toric fibration of the
K3 surface so that, via base diagrams, we can identify surfaces that represent
such a basis for $H_2(K3,\bZ)$.

We do this by presenting the K3 surface via an almost toric
fibration of E(1)$=\bC P^2 \#9\overline{\bC P}^2$ with two key
properties:
\begin{itemize}
\item the reduced boundary has one smooth component that is the whole
boundary and
\item we can identify visible surfaces, disjoint form the preimage of
the reduced boundary, that 
represent a basis for $H_2(E(1),\bZ)$ with respect to which 
the intersection form is $-E_8\oplus \left( 
\begin{smallmatrix} 0 & 1 \\  1 & -1 \end{smallmatrix}\right)$.
\end{itemize}
The first property allows us to define an almost toric fibration of 
the K3 surface by identifying two copies of this base diagram along
their boundaries.  
This amounts to performing 
a symplectic sum of two copies of E(1) along the tori that are the
preimage of the reduced boundary.
The second property allows us to immediately see most of the 
surfaces we are looking for in the K3 surface.
We  explain how to identify all $22$ of them.

\begin{figure}
\begin{center}
        \psfragscanon
        \psfrag{eta10}{$\eta_0$}
        \psfrag{eta1}{$\eta_1$}
        \psfrag{eta2}{$\eta_2$}
        \psfrag{eta3}{$\eta_3$}
        \psfrag{eta4}{$\eta_4$}
        \psfrag{eta5}{$\eta_5$}
        \psfrag{eta6}{$\eta_6$}
        \psfrag{eta7}{$\eta_7$}
        \psfrag{eta8}{$\eta_8$}
        \psfrag{eta9}{$\eta_9$}
	\includegraphics[scale=.5]{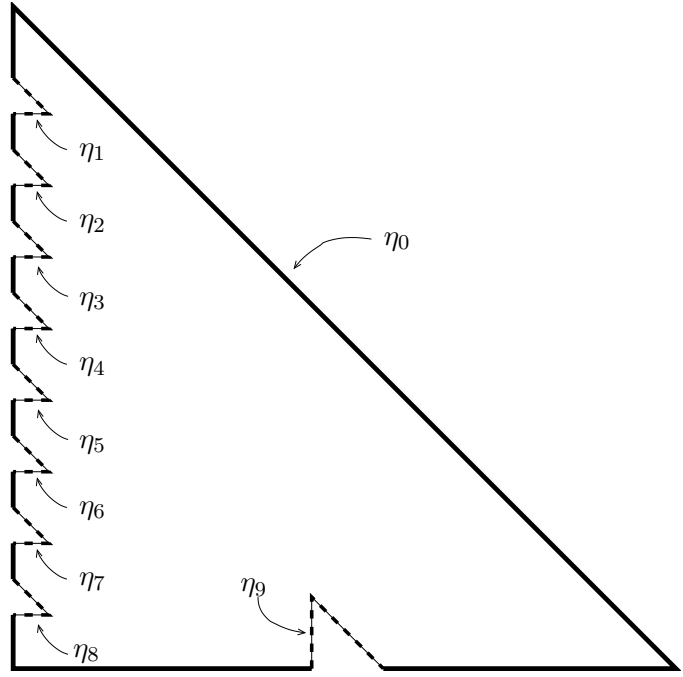}
\end{center}
\caption{Visible surfaces showing basis for $H_2(E(1),\bZ)$.}
\label{E1.fig}
\end{figure}

As preparation we start with a base diagram (Figure~\ref{E1.fig})
for an almost toric fibration of
E(1)$=\bC P^2 \# 9\overline{\bC P}^2$ in which we can see visible surfaces
$\Sigma_{\eta_0},\Sigma_{\eta_1},\ldots \Sigma_{\eta_9}$ where
each $\Sigma_{\eta_i}$, $i\ge 1$ is a $-1$-sphere coming from a blow-up.
Letting $h=[\Sigma_{\eta_0}]$ and $e_i=[\Sigma_{\eta_i}]$ for $i\ge 1$
we get a basis for the second homology in which the intersection form
is $<1>\oplus 9<-1>$.

\begin{figure}
\begin{center}
        \psfragscanon
        \psfrag{g1}{$\gamma_1$}
        \psfrag{g2}{$\gamma_2$}
        \psfrag{g3}{$\gamma_3$}
        \psfrag{g4}{$\gamma_4$}
        \psfrag{g5}{$\gamma_5$}
        \psfrag{g6}{$\gamma_6$}
        \psfrag{g7}{$\gamma_7$}
        \psfrag{g8}{$\gamma_8$}
        \psfrag{g9}{$\gamma_9$}
        \psfrag{g10}{$\gamma_{10}$}
        \psfrag{g11}{$\gamma_{11}$}
	\includegraphics[scale=.5]{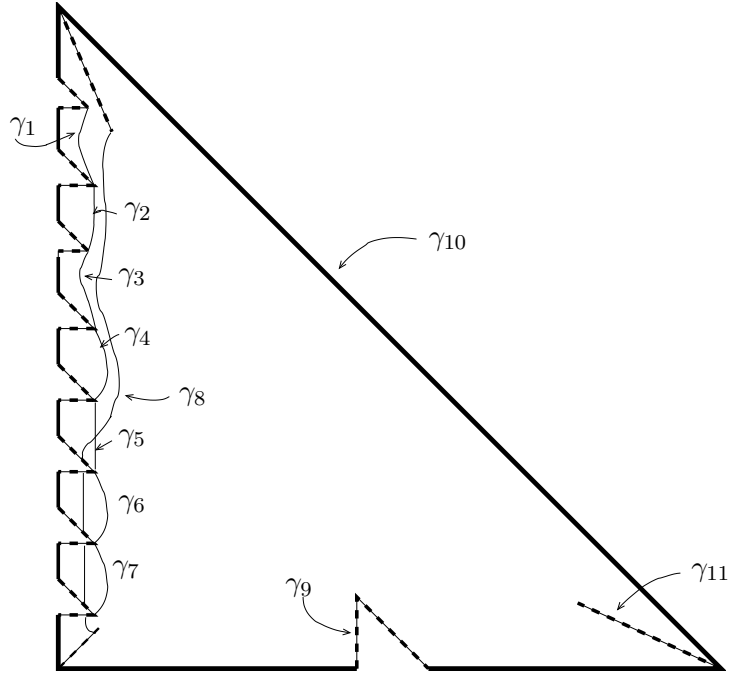}
\end{center}
\caption{Another basis for $H_2(E(1),\bZ)$.}
\label{E8.fig}
\end{figure}

Now making nodal trades at each vertex 
we get the fibration defined by
the base diagram in Figure~\ref{E8.fig}.
The visible surfaces $\Sigma_{\gamma_i}$, $i=1,\ldots 10$ represent
a basis for the second homology with respect to which we get the
intersection form $-E_8\oplus \left(  
\begin{smallmatrix} 0 & 1 \\  1 & -1 \end{smallmatrix}\right)$.
Indeed, $[\Sigma_{\gamma_j}]=e_j-e_{j+1}$, $j=1,\ldots 7$, 
$[\Sigma_{\gamma_8}]= e_6+e_7+e_8-h$, $[\Sigma_{\gamma_9}]=e_9$, and 
$[\Sigma_{\gamma_{10}}] =3h-\sum_{i=1}^9 e_i$
when the appropriate orientations are chosen.
Notice that $\Sigma_{\gamma_8}$ does define a smooth sphere:
With respect to the induced standard coordinates the compatible vector
must be the vanishing covector at each endpoint, namely $(-1,1)$ and
$(2,1)$; each time $\gamma_8$ crosses a wedge corresponding to a
blow-up the coordinates of the vector change by an application of
$\left(\begin{smallmatrix} 
0 & 1 \\  1 & 1 \end{smallmatrix}\right)$;
since 
$\left(\begin{smallmatrix} 
0 & 1 \\  1 & 1 \end{smallmatrix}\right)^3
\left(\begin{smallmatrix} -1 \\  1 \end{smallmatrix}\right)
= \left(\begin{smallmatrix} 2 \\  1 \end{smallmatrix}\right) $ the compatible
vector is well defined.

Using the techniques developed in Section~\ref{visible.sec}
the homology classes represented by the surfaces $\Sigma_{\gamma_i}$
can be read off the base diagram (for which the previous fibration of
E(1) may be helpful) and it 
is straightforward to check the intersection pairing of
the surfaces.

Now perform the symplectic sum of the two copies of E(1):
remove the torus $\Sigma_{\gamma_{10}}$ from each E(1) and replace
it with a $3$-torus (i.e. replace the almost toric fibered E(1) with
its boundary recovery, Definition~\ref{bdyrecovery.def}) and then
identify the two resulting four-manifolds $M_1,M_2$ along their boundaries.
In doing so, we can lift $\Sigma_{\gamma_{10}}$ to a torus inside
the boundary of each $M_i$ and can identify $\bdy M_1=\bdy M_2$ so
that the lifted two-tori in the two boundaries are identified.
Call the homology class of this torus in K3 surface $T_1$.

Because each $\gamma_i$, $i=1,\ldots 8$ is disjoint from the reduced
boundary, the spheres $\Sigma_{\gamma_i}$, $i=1,\ldots, 8$ embed into
the K3 surface after we have performed the symplectic sum.
Furthermore, the induced map on homology is injective on the subspace
spanned by the corresponding homology classes so the
images of the $[\Sigma_{\gamma_j}]$, $j=1,\ldots, 7$ 
and of $[\Sigma_{\gamma_8}]$ represent
distinct homology classes that we name $V_j^i$ and $W^i$ respectively,
where $i$ indicates which copy of E(1) the surface came from.

Meanwhile, the surfaces $\Sigma_{\gamma_9}$ in each copy of E(1) lift to 
closed disks in $M_1,M_2$ and the visible disks $\Sigma_{\gamma_{11}}$
in each $E(1)$ remain visible disks in $M_1,M_2$.
It is easy to perform the symplectic sum in such a way that the
boundaries of these corresponding visible disks in each $M_i$ are identified.
The result will be two spheres of self-intersection $-2$, the first
of which intersects $T_1$ transversely in one point.  
Call the homology classes of these spheres $S_1$ and $S_2$ respectively.
Let $S_3$ be the homology class of a section of the almost toric
fibration of the K3.

Inside the $3$-torus boundary of $M_1$ we can find two
other tori that are products of affine circles so that together with
the lift of $\Sigma_{\gamma_{10}}$ they span $H_2(\bdy M_1,\bZ)$.
Let $T_2,T_3$ be the homology classes of the images in the K3 surface 
(via inclusion) of these tori, with $T_3$ corresponding to the one that
is a Lagrangian fiber in $M_1$.

Then it is straightforward to check that
$V_j^i,W^i,T_k,S_k$, $i=1,2$, $j=1,\ldots 7$, $k=1,2,3$ is
a basis for the second homology of the K3 with respect to which 
the intersection form is $-E_8\oplus -E_8\oplus 3
\left(\begin{smallmatrix} 0 & 1 \\  1 & -2 \end{smallmatrix}\right)$.

\appendix
\section{Closed toric four-manifolds} \label{closedtoric.sec}

The following basic theorem is well known but direct proofs in
the symplectic category are difficult to find.

\begin{thm} \label{toricfourmanifolds.thm}
Suppose $M$ is a closed four-manifold that admits a symplectic structure
and toric fibration
$\pi:(M,\omega)\rightarrow (B,\cA,\cS)$.  Let $k$ be the number of
vertices of the base.
Then $M$ is diffeomorphic to $\bC P^2\#(k-3)\overline{\bC P}^2$ 
if $k\ne 4$ and
is diffeomorphic to  $\bC P^2\#\overline{\bC P}^2$ 
or $S^2\times S^2$ if $k=4$.
The latter two cases can be distinguished by whether or not (respectively)
a pair of non-adjacent edges has odd cross product (thereby
revealing a sphere with odd self-intersection).
\end{thm}

Delzant~\cite{Delzant.moment}
proved that a closed symplectic toric manifold is diffeomorphic to a
smooth compact toric variety.  (See also the appendix 
of~\cite{Audin.torus} for the dimension four case.)
Then one can appeal to the classification of smooth compact toric
complex surfaces, cf.~\cite{Fulton.toric}.
Here we provide a simple proof that is in keeping 
with the perspective of
this paper, namely we classify the bases $B$ that
can define a toric fibration of a closed four-manifold.
When there are few vertices in the base (three or four) it is easy to 
determine
the possible bases and the corresponding manifolds.
The main issue is to prove that if there are more than four vertices
then at least one of the edges must be the image of a $-1$-sphere that
can be blown down (Proposition~\ref{-1spheres.prop}).

First we set some notation to be used throughout this section.  
Since the base of a toric fibration of a closed 
manifold is isomorphic to 
the image of the moment map 
for a torus action $(B,\cA)$ must be isomorphic to a convex polygon 
$(B,\cA_0)\subset(\bR^2,\cA_0)$.
Label vertices of the polygon in
$\bR^2$ clockwise as $x_i$, $i=1,\ldots n$.
Let $e_i$ be the primitive integral vector that points from $x_{i+1}$ to $x_i$
(where $x_{n+1}=x_1$) and let $S_i$ be the sphere whose image is
the edge with vertices $x_i,x_{i+1}$.

The fact that $M$ is a smooth closed manifold implies that 
at each vertex $x_i$ we have $\abs{e_{i+1}\times e_i}=1$.
The clockwise numbering of the vertices and the convexity condition then imply
the refinement
\begin{equation} \label{smooth.eq}
e_{i+1}\times e_i=1.
\end{equation}
Furthermore,
\begin{lemma} \label{edgeint.lem}
The sphere $S_j$ has self-intersection $e_{j-1}\times e_{j+1}$.
\end{lemma}
\begin{xca}  
Prove this as a corollary of Proposition~\ref{sphereint.prop}
\end{xca}

\begin{lemma} \label{threefour.lem}
Let $\pi:(M,\omega)\rightarrow (B,\cA_0,\cS)$ be a toric fibration of a closed
manifold $M$.
If the base has
three vertices then $M$ must be diffeomorphic to $\bC P^2$ and
there is a unique fibration (up to scaling).
If the base has four vertices then $M$ is diffeomorphic to 
$S^2\times S^2$ or the non-trivial $S^2$ bundle over $S^2$ (or, 
equivalently, $\bC P^2\#\overline{\bC P}^2$).
\end{lemma}

\begin{proof}
Assume without loss of generality that $x_1=(a_1,0)$ for
some $a_1>0$, $x_2=(0,0)$ and $x_3=(0,a_2)$.
This implies $e_1=(1,0)$, $a_2>0$ and $e_2=(0,-1)$.
Suppose $(B,\cA_0,\cS)$ has three vertices.
The smoothness and convexity constraint (Equation~\ref{smooth.eq}) 
at $x_1$ and 
$x_3$ forces $e_3=(-1,1)$ and $a_1=a_2$.
This shows that up to scaling and integral affine isomorphism there is only
one base $(B,\cA_0,\cS)$ with three vertices.
Consulting Figure~\ref{B4CP2.fig}(b) we see that this base defines
$\bC P^2$ and that, up to scaling of the total volume
(Section~\ref{volume.sec}) all Lagrangian fibrations of $\bC P^2$ are
equivalent.

If $(B,\cA_0,\cS)$ has four vertices, the
smoothness and convexity constraint Equation~\ref{smooth.eq} 
at $x_3$ and $x_1$ forces $e_3=(-1,n_3)$ and
$e_4=(m_4,1)$.  
Applying the constraint at $x_4$ implies either $n_3=0$ or $m_4=0$.  
Our choice of
which edge to put on the horizontal axis allows us to
stipulate, without loss of generality, that $n_3=0$ and $m_4=n\ge 0$.
Appealing to Example~\ref{spherenbhd.ex} and our 
calculations of self-intersections in Section~\ref{selfint.sec}
we see that $M$ is the union of two disk bundles over the sphere glued
along their boundaries and hence is a sphere bundle over a sphere.
Because $\pi_1({\rm Diff}^+ S^2)=\pi_1(SO(3))=\bZ_2$ there are only two
possibilities: a trivial bundle $S^2\times S^2$ and a non-trivial
bundle.
Taking $S_1,S_2$ as a basis for the second homology the
intersection forms is
$
Q_{M}=\left(  
\begin{smallmatrix} 
-n & 1 \\  
1 & 0
\end{smallmatrix} 
\right)$.
This implies $M$ is diffeomorphic to  $M_n=S^2\times S^2$ when $n$ is even and
$M_n=S^2\tilde\times S^2$ when $n$ is odd.
To see that the non-trivial sphere bundle is diffeomorphic to 
$\bC P^2\#\overline\bC P^2$, notice that when $n=1$, $S_1$ is
the image of a $-1$-sphere and that blowing
this sphere down (Section~\ref{blowdown.sec}) yields $\bC P^2$.
\end{proof}

It is not so convenient to enumerate the possible bases $(B,\cA_0,\cS)$ that
have five or more vertices (and satisfy  $\bdy_R B=\bdy B$ so that they define 
a closed manifold).
However,

\begin{prop} \label{-1spheres.prop}
If $\pi:(M,\omega)\rightarrow (B,\cA_0,\cS)$ is a toric fibration of a closed 
manifold $M$ such that the base has at least five vertices,
then some edge in $\bdy B$ is the image of a sphere
with self-intersection $-1$.
\end{prop}

\begin{lemma} \label{nonnegspheres.lem}
If $(B,\cA_0,\cS)$ has  at least five vertices, then 
there are at most two edges in $\bdy B$ that are the 
image of a sphere with non-negative self-intersection and
these two edges must be adjacent.
Furthermore, only one of these two spheres can have positive self-intersection.
\end{lemma}

\begin{proof}

Position the vertices $e_i$ so that their bases are all at the origin
and let $\alpha_i$ be the angle between $e_{i}$ and $e_{i+1}$ and
observe that $\sum_i\alpha_i=2\pi$.
(Note that having done so one sees the {\it fan} corresponding to
the complex toric variety diffeomorphic to the manifold defined by 
$(B,\cA_0,\cS)$).
Lemma~\ref{edgeint.lem} implies that 
an edge of $B$ with tangent vector $e_i$ 
is the image of a sphere $S_i$ of non-negative
self-intersection if and only if $\alpha_{i-1}+\alpha_i\ge \pi$.

Suppose that three of the spheres have non-negative self-intersection.  
Then the condition on the
angles implies the spheres are adjacent and there are at most four
angles, and hence at most four edges contradicting our assumption of
at least five edges.

Knowing now that there are at most two such spheres, suppose there
are two and that they are not adjacent.  Let $S_2$ be one of
the two spheres with non-negative self-intersection and, without
loss of generality, let $e_1=(1,0)$.
Then $\alpha_{1}+\alpha_2\ge \pi$ implies $e_3=(a,b)$ with $b\le 0$.
Since $\sum_i\alpha_i=2\pi$, the spheres of non-negative 
self-intersection being adjacent forces the other such sphere
to be $S_4$, $e_3=(-1,0)$.
Hence both $S_2,S_4$ must have self-intersection zero and
there are only four spheres, again a contradiction.

Finally, suppose there is a pair of adjacent spheres having
positive self-intersection.
Refer again to the base in $(\bR^2,\cA_0)$ with the vectors $e_1,e_2$
chosen as in the proof of Lemma~\ref{threefour.lem}. 
Smoothness at $x_{1}$ and $x_{3}$ implies that
$e_{n}=(-k,1)$ and $e_{3}=(-1,l)$ for some $k,l>0$.  But then convexity
of the image implies that $k=l$ and $e_{3}=e_{n}$ so there
are only three edges, a contradiction.
\end{proof}

\begin{proof}[Proof of Proposition~\ref{-1spheres.prop}]
Assume that no $S_i$ is a $-1$-sphere.

By Lemma~\ref{nonnegspheres.lem} we can assume that either
\begin{enumerate}
\item $S_i\cdot S_i\le-2$
for all $i$, or
\item $S_n\cdot S_n\ge 0$, $S_{n-1}\cdot S_{n-1}=0$ and $S_i\cdot S_i\le-2$
for all other $i$, or
\item $S_n\cdot S_n\ge 0$ and $S_i\cdot S_i\le-2$
for all other $i$.
\end{enumerate}

Appealing to Exercise~\ref{convexity.ex} 
we see right away that the
first case is impossible since $(B,\cA_0)$ is a closed polygon.
In the second case, since $S_{n-1}\cdot S_{n-1}=0$,
we can assume that  $e_{n-2}=(1,0)$,
$e_{n-1}=(0,-1)$ and $e_{n}=(-1,0)$.
Then $S_n\cdot S_n\ge 0$
implies $e_1=(m,1)$ with $m\ge 0$ while $S_i\cdot S_i\le-2$
for $1\le i\le n-2$ implies $e_1=(k,l)$ for some relatively prime
$k,l>0$ such that $\frac{k}{l}=[b_{n-2},\ldots,b_2]$ with $b_{j}\ge 2$
for all $2\le j\le n-2$.   
For these constraints on $e_1$ to agree we must have $l=1$ 
which forces $n-2=2$ in which case there are only four edges, 
a contradiction.
In the third case we can assume that $e_n=(0,-1)$ and $e_{n-1}=(1,0)$.
Then $S_n\cdot S_n\ge 0$ implies $e_1=(-1,k)$, $k\ge 0$ and $S_i\cdot S_i\le-2$
for $2\le i\le n-1$ implies  $e_1=(k,l)$ as in the previous case; but this
is a contradiction since $k,l>0$.
\end{proof}

\begin{proof}[Proof of Theorem~\ref{toricfourmanifolds.thm}]
We now know that we can blow down any closed toric manifold fibering over
a polygon with at least five sides to yield a toric fibered closed
manifold whose base has one less side.  It follows that any such manifold
is the blow-up of $S^2\times S^2$ or $\bC P^2\#\overline{\bC P}^2$.
But we can refine this statement a bit thanks to the following fact:

\begin{xca}   \label{minimalmodel.ex}
Use base diagrams (moment map images) for toric fibrations of 
$S^2\times S^2$ and $\bC P^2\#\overline{\bC P}^2$ to show that
$(S^2\times S^2)\#\overline{\bC P}^2$ and 
$\bC P^2\#2\overline{\bC P}^2$ are diffeomorphic.
\end{xca}

This completes the proof of Theorem~\ref{toricfourmanifolds.thm}.
\end{proof}

\bibliography{refs}

\end{document}